\documentclass[11pt]{amsart}
\usepackage{amsmath,amsxtra,amssymb}

\numberwithin{equation}{section}
\newtheorem{lemma}{Lemma}[section]
\newtheorem{prop}[lemma]{Proposition}
\newtheorem{theorem}[lemma]{Theorem}
\newtheorem{cor}[lemma]{Corollary}

\newtheorem{rem}[lemma]{Remark}

\newcommand{\re}{\begin{rem}\rm}
  \newcommand{\mar}{\end{rem}}
\newtheorem{exam}[lemma]{Example}

\newcommand{\kla}{\left ( }
\newcommand{\mer}{\right ) }

\renewcommand{\for}{\begin{eqnarray*}}

\newcommand{\mel}{\end{eqnarray*}}

\newcommand{\kl}{\pl \le \pl}
\newcommand{\gl}{\pl \ge \pl}

\newcommand{\lel}{\pl = \pl}

\newcommand{\nz}{{\mathbb N}}
\newcommand{\nen}{n \in \nz}

\newcommand{\rz}{{\mathbb R}}
\newcommand{\zz}{{\mathbb Z}}

\newcommand{\cz}{{\mathbb C}}
\newcommand{\kz}{{\rm  I\! K}}

\newcommand{\ten}{\otimes}

\newcommand{\wet}{\stackrel{\wedge}{\otimes}}

\newcommand{\p}{\hspace{.05cm}}
\newcommand{\pl}{\hspace{.1cm}}
\newcommand{\pll}{\hspace{.3cm}}

\newcommand{\hz}{\vspace{0.5cm}}
\newcommand{\lz}{\vspace{0.2cm}}

\newcommand{\qd}{\end{proof}\vspace{0.5ex}}

\newcommand{\circled}[1]{\overset{\circ}{#1}}
\newcommand{\Om}{\Omega}
\newcommand{\om}{\omega}

\newcommand{\al}{\alpha}
\newcommand{\si}{\sigma}
\newcommand{\Si}{\Sigma}

\newcommand{\la}{\lambda}
\newcommand{\eps}{\varepsilon}

\renewcommand{\L}{{\mathcal L}}
\newcommand{\F}{{\mathcal F}}
\newcommand{\E}{{\mathcal E}}
\newcommand{\A}{{\mathcal A}}

\newcommand{\K}{{\mathcal K}}

\newcommand{\Pa}{{\mathcal P}}
\newcommand{\Ha}{{\mathcal H}}

\newcommand{\M}{{\mathcal M}}

\newcommand{\N}{{\mathcal N}}

\newcommand{\U}{{\mathcal U}}
\newcommand{\B}{{\mathcal B}}
\newcommand{\Hh}{{\mathcal H}}
\newcommand{\ff}{{\mathbb  F}}

\newcommand{\noo}{\left \|}
\newcommand{\rrm}{\right \|}

\newcommand{\bet}{\left |}
\newcommand{\rag}{\right |}

\newcommand{\intt}{\int\limits}
\newcommand{\summ}{\sum\limits}

\newcommand{\prodd}{\prod\nolimits}

\newcommand{\lb}{\langle \langle }
\newcommand{\rb}{\rangle \rangle }

\newcommand{\8}{\infty}

\DeclareMathOperator{\krn}{ker}

\allowdisplaybreaks
\begin{document}

\title[Embedding of OH and logarithmic `little Grothendieck inequality']
{\bf \large Embedding of the operator space OH and the logarithmic
`little Grothendieck inequality'}


\author[Marius Junge]{Marius Junge$^\dag$}
\address{Department of Mathematics\\
University of Illinois at  Urbana-Champaign, Urbana, IL 61801}
\email{junge@math.uiuc.edu}

\thanks{$^\dag$Partially supported by
National Science Foundation grant DMS-0301116}

\begin{abstract}We  use Voiculescu's concept of free probability to   construct a
completely isomorphic embedding of the operator space OH in the
predual of a von Neumann algebra. We analyze the properties of
this embedding and determine the operator space projection
constant of OH$_n$:
 \[  \frac{1}{108} \pl \sqrt{\frac{n}{1+\ln n}} \kl \inf_{P:B(\ell_2)\to OH_n, P^2=P}
  \noo P\rrm_{cb} \kl 288  \pi  \sqrt{\frac{2n}{1+\ln n}} \pl .\]
The lower estimate is a recent result of Pisier and Shlyakhtenko
that improves an estimate of order $1/(1+\ln n)$ of the author.
The additional factor $1/\sqrt{1+\ln n}$ indicates  that the
operator space OH$_n$ behaves differently than its classical
counterpart $\ell_2^n$. We give an application of this formula to
positive sesquilinear forms on $\B(H) $. This leads to logarithmic
characterization of $C^*$-algebras with the weak expectation
property introduced  by Lance.
 \end{abstract}
\maketitle

\newcommand{\llz}{\vspace{0.15cm}}

{\bf Contents:}
\begin{enumerate}
 \item[0.] Introduction and notation.\llz
 \item[1.] Preliminaries.\llz
 \item[2.] A logarithmic characterization of C$^*$-algebras with WEP.\llz
 \item[3.] The Pusz/Woronowicz formula and the operator space OH. \llz
 \item[4.] The projection constant of the operator space OH$_n$. \llz
 \item[5.] Norm calculations in a quotient space.\llz
 \item[6.] $K$-functionals.\llz
 \item[7.] Sums of free mean zero variables. \llz
\end{enumerate}

\newpage

\setcounter{section}{-1}
\section{Introduction and notation}
\setcounter{section}{0} \setcounter{equation}{0}

Probabilistic techniques  and concepts play an important role in
the theory of Banach spaces and operator algebras. For example,
the Khintchine inequality and the `little Grothendieck inequality'
are fundamental tools in Banach space theory. In the theory of
operator algebras probabilistic concepts are important in
Takesaki's proof of Sakai's theorem (see \cite{TAK}) and in
Connes' characterization of injective von Neumann algebras (see
\cite{Co}). Pisier/Haagerup's noncommutative version of
Grothendieck's inequality (see e.g. \cite{Psl}) and Grothendieck's
inequality for exact operator spaces in \cite{JP} use
probabilistic techniques. This latter result is inspired by
`Grothendieck's program for operator algebras' which motivates
fundamental research in the theory of operator spaces. In
Pisier/Shlyakhtenko's Grothendieck theorem for operator spaces
(see \cite{PS}) and for the  results in this paper the use of free
probability is crucial.


We follow Grothendieck's ideas and investigate positive
sesquilinear forms on $C^*$-algebras.  Let us first recall the
so-called `little Grothendieck inequality' on $C(K)$-spaces.
Grothendieck \cite{G} showed  that for every bounded linear map
$v:C(K)\to \ell_2$ there is a probability measure $\mu$ on $C(K)$
such that $(0.1)$ holds for all $f\in C(K)$:
 \begin{equation}
  \noo v(f)\rrm^2 \kl \frac{2}{\sqrt{\pi}} \noo v\rrm^2 \int_K
  |f|^2
 d\mu \pl .
 \end{equation}
The constant $\frac{2}{\sqrt{\pi}}$ is optimal for complex $C(K)$
spaces. Bounded positive sesquilinear forms on $C(K)$ (i.e.
possibly degenerate scalar products) are in one-to-one
correspondence with bounded linear maps $v:C(K)\to \ell_2$  via
$B(f,g)\lel (v(f),v(g))$. For a probability measure $\mu$ we may
define the positive sesquilinear form $B_{\mu}$ as follows
 \[ B_{\mu}(f,g) \lel \int_K \bar{f}g d\mu \pl. \]
According to Grothendieck's work, $B_\mu$ is the prototype of an
integral linear form. Integral forms are continuous functionals on
the injective Banach space tensor product
$\overline{C(K)}\ten_{\eps}C(K)\cong C(K\times K)$. Indeed,  we
have
 \begin{align*}
 \bet \intt \summ_{k=1}^n \bar{f}_kg_k d\mu \rag &\le \sup_{t\in
 K} \bet \summ_{k=1}^n \bar{f}_k(t)g_k(t)\rag \kl
  \sup_{t,s\in K} \bet \summ_{k=1}^n \bar{f}_k(t)g_k(s)\rag
 \end{align*}
for all finite sequences $(f_k)$ and $(g_k)$. We say that a
positive sesquilinear form $B$ is majorized by a bilinear form
$\tilde{B}$, in short $B\le \tilde{B}$, if
 \[ B(x,x)\kl |\tilde{B}(x,x)| \pl \]
holds for all $x$. Therefore, the `little Grothendieck inequality'
implies  that every bounded positive, sesquilinear form is
majorized by an integral linear form (even a  positive, integral,
sesquilinear form).

We will now discuss the analogue of this result in the context of
$C^*$-algebras. Let $A_1\subset \B(H)$ and $A_2\subset \B(K)$ be
$C^*$-algebras.  For $C^*$-algebras we shall replace the Banach
space injective tensor norm  by the smallest $C^*$-tensor  norm
$A_1\ten_{\min}A_2$ on $A_1\ten A_2$. This norm is given by the
inclusion $A_1\ten_{\min}A_2\subset B(H\ten K)$. In this context a
bilinear form $B:A_1\times A_1\to \cz$ is called an
\emph{integral} form if there exists a constant $C>0$ such that
 \[ \bet \summ_{k=1}^n B(x_k,y_k)\rag \kl C \noo
 \summ_{k=1}^n x_k\ten y_k\rrm_{A_1\ten_{\min}A_2} \pl \]
holds for all finite sequences $(x_k)\subset A_1,(y_k)\subset
A_2$. The smallest possible constant is the norm of $B$ as a
linear functional on $A_1\ten_{\min}A_2$ and will be denoted by
$\noo B\rrm_I$. From the theory of operator spaces it is clear
that the `appropriate' substitute for `bounded bilinear form' is
the notion of a `jointly completely bounded form'. A bilinear form
$B:A_1\times A_2\to \cz$ is \emph{jcb} (\emph{jointly completely
bounded}) if there exists a constant $C$ such that
 \begin{equation}\label{jcb}  \noo \summ_{k,j} B(a_k,b_j) \pl x_k\ten y_j\rrm_{M_{m}\ten_{\min}M_m}
  \kl C \noo \summ_{k} a_k\ten x_k\rrm_{M_m(A_1)}
  \noo \summ_{j} b_j\ten y_j \rrm_{M_m(A_2)} \pl \end{equation}
holds for all finite sequences $(a_k)\subset A_1,(b_j)\subset A_2$
and $(x_k),(y_j)\subset M_m$. The jcb-norm $\noo B\rrm_{jcb}$ is
given by the infimum over all $C$ satisfying \eqref{jcb}.
Following Grothendieck's categorial approach it is natural to ask
whether every positive, sesquilinear jcb form is majorized by an
integral form. In contrast to the commutative case this fails in
the noncommutative setting:\lz

\noindent{\bf Theorem 1.}\label{main} {\it There exists a
positive, integral, sesquilinear form on $\B(H)$ which can not be
majorized by an integral form. More precisely, for every $\nen$
there exists a positive, integral, sesquilinear form $B$ on
$\B(H)$ of rank $n$ such that
 \begin{equation}\label{lppb}
  (1+\ln n)\pl \noo B\rrm_{jcb} \kl C\| \tilde{B}\|_I
  \end{equation}
holds for all $\tilde{B}$ satisfying $B\le \tilde{B}$.} \lz

The factor $(1+\ln n)$ in \eqref{lppb} is optimal. Indeed,
Pisier/Shlyakhtenko \cite{PS} showed that if $B$ a positive,
integral, sesquilinear form on $\B(H)$ of rank $n$, then there
exists an integral, sesquilinear form $\tilde{B}$ with $B\le
\tilde{B}$ such that
 \begin{equation}\label{uppb}  \| \tilde{B}\|_{I}\kl C\pl (1+\ln n)\noo B\rrm_{jcb} \pl .
 \end{equation}
In fact, Pisier/Shlyakhtenko \cite{PS}  improved  an estimate of
the author (see \cite{Jo1}) of the order $(1+\ln n)^2$ in
\eqref{uppb}. Following Pisier's work (see \cite{Poh}), positive
sesquilinear jcb forms are closely connected to completely bounded
linear maps with values in the operator space OH (for definitions
see below). Our approach to Theorem 1 is probabilistic in nature.
We find an embedding of the operator space OH in a noncommutative
$L_1$ space, imitating the classical embedding of $\ell_2$ via
Gaussian variables. The properties of this embedding of OH then
yield the logarithmic term.

We  recall some operator space notation before giving more
details. An operator space $F$ comes either with a concrete
isometric embedding $\iota:F\to \B(H) $ or with a sequence $(\noo
\pl \rrm_m)$ of matrix norms on $(M_m(F))$  such that
 \[ \noo [x_{ij}] \rrm_{M_m(F)} \lel \noo
 [\iota(x_{ij})]\rrm_{B(H^m)} \pl .\]
Ruan's axioms (see e.g. \cite{ER}) describe axiomatically those
sequences of matrix norms which can occur from an isometric
embedding in $\B(H) $. The morphisms in this category are
completely bounded linear maps $u:E\to F$, i.e. linear maps such
that
 \[ \noo u\rrm_{cb} \lel \sup_m \noo id\ten u:M_m(E)\to M_m(F)\rrm \]
is finite. We denote by $CB(E,F)$ the Banach (operator) space of
completely bounded maps equipped with this norm. Hilbertian
operator spaces are of particular interest. For example the column
and row spaces of matrices
 \[ K^c\lel \B(\cz,K) \pl \subset \B(K) \quad \mbox{and}
 \quad  K^r \lel \B(K,\cz) \subset \B(K) \pl \]
play a fundamental r\^{o}le.  Pisier discovered that the sequence
of norms  on $M_m\ten K$ obtained by the complex interpolation
method
 \begin{equation}\label{ohdef}  M_m(K^{oh}) \lel [M_m(K^c),M_m(K^r)]_{\frac12} \
  \end{equation}
defines a sequence of matrix norms on $M_m\ten K$ satisfying
Ruan's axioms. Thus \eqref{ohdef} defines an operator space
structure on $K$ denoted by $K^{oh}$. In particular, for
$K=\ell_2$ we find a sequence of operators $(T_k)\subset
\B(\ell_2)$ satisfying
 \[ \noo \summ_{k=1}^\8 x_k \ten T_k \rrm_{M_m(\B(\ell_2))}
  \lel \noo \summ_{k=1}^\8 \bar{x}_k\ten x_k \rrm_{\overline{M_m}\ten
  M_m}^{\frac12}
  \lel \sup_{\noo a\rrm_2,\noo b\rrm_2\le 1}
  \kla \summ_{k=1}^\8 tr(ax_k^*bx_k)\mer^{\frac12} \pl
  \]
for all sequences $(x_k)\in M_m$. The {\it operator Hilbert space}
$OH=\ell_2^{oh}$ is the span of the $(T_k)$'s. It is still unclear
how to construct `concrete' operators $(T_k)$ satisfying this
equality. The investigations in this paper may be considered as a
starting point in this direction. We need the concept of the
\emph{standard} dual $E^*$ of an operator $E$. The operator space
structure  on $E^*$ is given by the isometric equality
 \[ M_m(E^*)\lel CB(E,M_m)\pl. \]
Here a matrix $[x^*_{ij}]$ of functionals corresponds to the
linear map $x\mapsto [x^*_{ij}(x)]$. In  particular, duals of
$C^*$-algebras and preduals of von Neumann algebras carry a
natural operator space structure.

With the help of the operator space dual it is easy to explain
why the notion of jcb forms is the natural analogue of bounded
bilinear forms on  Banach spaces. Indeed, a bilinear form
$B:F\times E\to \cz$ is jcb if and only if the corresponding
linear map $T_B:E\to F^*$, $T_B(x)(y)=B(y,x)$ is completely
bounded. We will now explain the connection to the operator space
OH.  Given a positive sesquilinear form $B:\bar{E}\times E\to
\cz$, we may consider $L=\{x\pl|\pl B(x,x)=0\}$. We denote by $K$
the Hilbert space obtained by completion of $E/L$ with respect to
induced scalar product $(x+L,x+L)=B(x,x)$. Then the natural map
$v:E\to K^{oh}$ defined by $v(x)=x+L$ is completely bounded if
and only if $B$ is jcb (see \cite{Poh}). This equivalence yields
a one to one correspondence between positive, sesquilinear jcb
forms and completely bounded maps with values in OH. It   allows
us to derive Theorem 1 from properties of the operator space OH.

A key new ingredient is a formula of Pusz/Woronowicz for the
square root of two sesquilinear forms (see \cite{Wo} and
\eqref{wo}). In section 3, we  show how this formula (and its new
dual version) provides a concrete realization of OH as a subspace
of a quotient of the direct sum $R\oplus C$. Let us be more
specific.  Inspired by \cite{Wo} we consider the probability
measure $d\mu(t)=dt/(\pi \sqrt{t(1-t)})$ on $[0,1]$ and the two
measures $d\nu_1(t)=t^{-1}d\mu(t)$, $d\nu_2(t)=(1-t)^{-1}d\mu(t)$.
Then the direct sum
 \[ \Ha \lel  L_2^c(\nu_1;\ell_2) \oplus L_2^r(\nu_2;\ell_2) \]
is an operator space. On $\Ha$, we define the map $Q:\Ha\to
L_0(\mu;\ell_2)$ by
 \[ Q(x_1,x_2)(t) \lel x_1(t)+x_2(t) \in \ell_2 \pl .\]
It is easily checked that $M_m(\Ha/S) \lel M_m(\Ha)/M_m(S)$
defines a sequence of matrix norms satisfying Ruan's axioms and
therefore
 \[ G \lel \Ha/\krn(Q) \pl \]
is an operator space. Then, we may consider the subspace $F\subset
G$ of equivalence classes $(x_1,x_2)=\krn(Q)$ such that $x_1+x_2$
is a $\mu$-almost everywhere a constant element in $\ell_2$. The
operator space structure of OH is encoded in $\mu$ and the two
densities $1/t$ and $1/(1-t)$:

\hz \noindent {\bf Theorem 2.} {\it  $F$ is $2$-completely
isomorphic to OH.} \hz

The embedding of $G$ into the predual of a von Neumann algebra
uses free probability. We first  extend Voiculescu's inequality
for the norm of the  sums of free independent random variables to
the operator-valued setting. This operator-valued  version of
Voiculescu's inequality provides estimates for the cb-norm of
linear maps. Voiculescu's inequality naturally involves three
terms. Using a  central limit procedure, we may eliminate one of
them. These methods can be used to show that arbitrary quotients
of $R\oplus C$ embed in the predual of a von Neumann algebra. In a
subsequent paper \cite{JO} we will elaborate this fact and
construct an embedding of OH in the predual of a hyperfinite
factor. For our applications, it is important  to know that the
underlying von Neumann algebra is QWEP. Let us recall that a
$C^*$-algebra has the weak expectation property (WEP) if there
exists  a (complete) contraction $P:\B(H) \to A^{**}$ such that
$P|_{A}=id_A$. A $C^*$-algebra is QWEP if it is the quotient
$A=\A/I$ of  some $C^*$-algebra $\A$ with WEP by a two-sided ideal
$I$. It is an open question whether every $C^*$-algebra is QWEP
(see \cite{Ki}).

\hz

\noindent  {\bf Theorem 3.} {\it $G$ is completely isomorphic to a
completely complemented subspace of the predual $N_*$ of a  von
Neumann algebra $N$ with QWEP. In particular, OH embeds into
$N_*$.}

\hz

Recently Pisier \cite{Ps5} showed that no embedding of OH is
possible in the predual of a semifinite von Neumann algebra. Type
III von Neumann algebras are indeed necessary for embedding OH in
noncommutative $L_1$ spaces.

The factor $(1+\ln n)$ is a result of norm calculations in the
predual of the tensor product of two von Neumann algebras. This
approach is again motivated by Grothendieck's work on absolutely
$1$-summing maps. Let us denote by $\pi_1^B$ the absolutely
summing norm for Banach spaces (see \cite{Psl}). Let
$g_1,...,g_n,g_1',...,g_n'$ be independent, normalized, complex
Gaussian variables. Following Grothendieck's work we know that
 \[ \pi_1^B(id_{\ell_2^n}) \lel  \frac{\pi}{4}
  \noo \summ_{i=1}^n g_ig_i'\rrm_{L_1([0,1]^2)} \sim
 \frac{\sqrt{\pi}}{2} \sqrt{n} \pl .\]
The first `equality' remains true in the operator space context if
we replace independent Gaussian random variables by a suitable
tensor product of noncommutative random variables. However,
calculating this tensor norm in the noncommutative context is more
involved. We show that it can be calculated as an element of a
$4$-term quotient of classical Banach spaces (see section 5). The
outcome of these norm calculations provides  the logarithmic
factor:

\hz \noindent  {\bf Theorem 4.} {\it Let $u:G\to N_*$ be the
embedding from Theorem 3 and $(f_k)$ be the unit vectors basis in
$F$ and $\nen$. Then
 \[ \noo \summ_{k=1}^n u(f_k)\ten u(f_k)\rrm_{(N\bar{\ten} N)_*} \sim
 \sqrt{n(1+\ln n)} \pl. \]
} \hz

In section 4, we show that Theorem 4 implies an estimate on the
completely $1$-summing norm of the  identity map on OH$_n$ (see
section 4 for a definition). Using the well-known concept of trace
duality, we obtain  estimates for the  operator space projection
constant
 \[ \la_{cb}(OH_n) \lel\inf_{P:\B(\ell_2)\to OH_n,P|_{OH_n}=id} \noo P\rrm_{cb}
 \pl .
 \]

 \noindent  {\bf Corollary 5.} {\it Let $\nen$.  Then
$\la_{cb}(OH_n) \sim \sqrt{\frac{n}{1+\ln n}}$.}
  \hz

If $P:\B(\ell_2)\to OH_n$ is the optimal projection, then
$B(x,y)=(P(x),P(y))$ provides an example for Theorem 1 where the
logarithmic term is necessary (see again section 4). I learned
from C. le Merdy that the analogue of the `little Grothendieck
inequality' fails for the reduced $C^*$-algebra of the free group
in $n$-generators. Using Haagerup's characterization of
$C^*$-algebras with WEP in terms of selfpolar forms, we can show
in section 2 that WEP is the crucial property.

\hz

\noindent  {\bf Theorem 6.} {\it Let $A$ be a $C^*$-algebra and
$\al>0$. $A$ has WEP if and only if  there exists a constant
$C_{\al}$ such that every positive sesquilinear form $B$ of rank
$n$ on $A$ is majorized by an integral sesquilinear from
$\tilde{B}$ satisfying
 \[ \| \tilde{B}\|_I\kl C_{\al} (1+\ln n)^{\al} \|B\|_{jcb}\pl .\]
Moreover,  if $A$ is not a subalgebra of $C(K,M_m)$ for some $m$,
then this condition is satisfied only for $\al\gl 1$.} \hz

In section 1 we  provide some background and notation. Theorem 6
is proved in section 2. The Pusz/Woronowicz formula and its
application to OH is contained in section 3. Recently, alternative
pairs of measures have been found (see \cite{JX2}) which lead to
nicer representations of OH but are beyond the scope of this
paper. In section 4 we prove Theorem 1 assuming the probabilistic
result Theorem 3 (see section 7) and the norm calculation in
section 5. In section 6, we investigate different notions of
$K$-functionals, in particular $3$-term $K$-functionals. These
$K$-functionals arise naturally in the context of Voiculescu's
inequality in section 7. At the end of the paper we show that the
free product of von Neumann algebras with QWEP is again QWEP (see
Theorem \ref{freeQWEP}). This result might be of independent
interest.

We would like to thank U. Haagerup for the  collaboration on the
proof of Proposition \ref{Voi}. We also thank the anonymous
referees, J. Parcet, Anthony Yew and in particular Q. Xu  for
patient readings of different versions of this paper and many
helpful suggestions.

\section{Preliminaries}
We use standard notation in operator algebras as in
\cite{TAK,KR,Strat}. Let $A_1\subset \B(H)$ and $A_2\subset
\B(K)$ be $C^*$-algebras. The norm on $A\ten B$ induced by the
inclusion $A\ten B\subset \B(H\ten_2 K)$ is called the minimal
$C^*$-norm. We use $A\ten_{\min}B$ for the completion of $A\ten
B$ with  respect to that norm. Here we have used $H\ten_2 K$ for
the unique tensor product making $H\ten K$ a Hilbert space (often
called the Hilbert-Schmidt norm). Let  $N\subset \B(H)$,
$M\subset \B(K)$ be von Neumann algebras. We denote by
$N\bar{\ten}M$ the closure of $N\ten_{\min} M$ in the weak
operator topology. For a $C^*$-algebra $A$, we denote by $A^{op}$
the C$^*$-algebra defined on the same underlying Banach space but
with the reversed multiplication $x\circ y=yx$. By $\bar{A}$, we
denote the $C^*$-algebra obtained from he complex multiplication
$\la.x\lel \bar{\la}x$ on $A$. Thus $A$ and $\bar{A}$ coincide as
real Banach algebras.  We see that the map $j:A^{op}\to \bar{A}$
given by $j(x)=x^*$ is a $C^*$-isomorphism.

The notation $E\ten_{\eps}F$ is used for the completion of $E\ten
F$ with respect to the largest tensor norm in the category of
Banach spaces. Similarly, we will use $E\ten_{\pi} F$ for the
completion of $E\ten F$ with respect to  the biggest tensor norm
in the category of Banach spaces. For Hilbert spaces $H$ and $K$,
the space of trace class operators is denoted by  $H\ten_{\pi}K$
or $S_1(H,K)$. Similarly, we  use the notation
$H\ten_{\eps}K=\K(H,K)$ for the space of compact operators. We
note the trivial inclusions
 \begin{eqnarray}
  H\ten_{\pi} K &\subset& H\ten_2 K \pl \subset \pl
  H\ten_{\eps} K \pl .
  \end{eqnarray}
We assume the reader to be familiar with standard operator space
terminology which can be found in the monographs \cite{ER} or
\cite{P}. We will need some basic facts about the column Hilbert
space $K^c=\B(\cz,K)$ and the row Hilbert space $K^r=\B(K,\cz)$
of  a given Hilbert space $K$. Given an element $x=[x_{ij}]\in
M_m(K^c)$, the norm is given by
 \begin{eqnarray}
  \noo x\rrm_{M_m(K^c)} &=& \noo x\rrm_{\B(\ell_2^m,\ell_2^m(K))}
  \lel  \noo x^*x\rrm_{\B(\ell_2^m)}^{\frac12} \lel
  \noo \bigg[\summ_k (x_{ki}, x_{kj})\bigg]_{ij} \rrm_{M_m}^{\frac12} \pl .
  \end{eqnarray}
Here $(x,y)$ denotes the scalar product of $x$ and $y$. We will
assume that  scalar products are  antilinear in the first
component. Similarly, we have
 \begin{eqnarray}
 \noo x\rrm_{M_m(K^r)} &=&
 \noo \bigg[\summ_k (x_{ik}, x_{jk})\bigg]_{ij} \rrm_{M_m}^{\frac12} \pl .
  \end{eqnarray}
We use the standard notation $C=\ell_2^c$, $C_n=(\ell_2^n)^c$ and
$R=\ell_2^r$, $R_n\lel (\ell_2^n)^r$ for the column and row
Hilbert spaces.  We refer to \cite{Poh,P} for more details on the
operator Hilbert space $K^{oh}=[K^c,K^r]_{\frac12}$ and
interpolation norms. For $K=\ell_2$ and  $K=\ell_2^n$, we will
simply write $OH\lel \ell_2^{oh}$ and $OH_n=(\ell_2^n)^{oh}$. We
denote by $(e_k)$ the natural unit vector basis. Given a sequence
$(x_k)\subset \B(H) $, we define the  associated linear map
$u:OH\to \B(H) $ defined by $u(e_k)=x_k$ and have
 \begin{eqnarray} \label{ohnorm}
 \noo u:OH\to {\rm Im}(u)\rrm_{cb} &=& \noo u\rrm_{cb}
 \lel \noo \summ_k \bar{x}_k\ten
  x_k\rrm_{\overline{\B(H) }\ten_{\min} \B(H) }^{\frac12}
 \pl .
 \end{eqnarray}
The Hilbertian  operator spaces $K^c$, $K^r$ and $K^{oh}$ are
homogeneous, i.e. for $s\in \{ c,r,oh\}$ and every bounded linear
map $u:K^s\to K^s$, we have
 \[ \noo u\rrm_{cb} \lel \noo u\rrm \pl .\]

In terms of general operator space notation, let us recall that a
\emph{complete contraction} $u:E\to F$ is given by a completely
bounded map with $\noo u\rrm_{cb}\le 1$. We say that two operator
spaces $E$ and $F$ are \emph{$\la$-cb isomorphic} if there exists
a linear isomorphism  $u:E\to F$ such that $\noo u\rrm_{cb}\noo
u^{-1}\rrm_{cb}\le \la$. An operator space $E$ is
\emph{$\la$-completely complemented} in an operator space $F$ if
there exist completely bounded maps $u:E\to F$ and $v:F\to E$ such
that $vu=id_{E}$ and $\noo u\rrm_{cb}\noo v\rrm_{cb}\kl \la$.

Quotient operator spaces $G=V/E$ are important in this paper. We
refer to the introduction for the definition of the operator space
structure $M_m(G)=M_m(V)/M_m(E)$. If $T:V\to X$ is a completely
bounded map  which vanishes on $F$, then $T$ induces a unique map
$\hat{T}:V/E\to X$  defined by $\hat{T}(x+F)=T(x)$. We have $\noo
T\rrm_{cb}=\| \hat{T}\|_{cb}$.

The \emph{injective tensor product} $E\ten_{\min}F$ of two
operator spaces $E\subset \B(H)$, $F\subset \B(K)$ is the
completion of $E\ten F$ with respect to the norm induced by the
inclusion map $E\ten_{\min}F\subset \B(H\ten_2 K)$. The norm does
not depend on the underlying completely isometric embedding. The
\emph{projective tensor product} $E\wet F$ is defined such that
 \[ (E\wet F)^*\lel CB(E,F^*)\cong CB(F,E^*) \pl .\]
See \cite{P,ER} for details. Recall  that a bilinear form
$B:E\times F\to \cz$ is jcb if and only if its linear extension
to  $B:E\wet F\to \cz$ is continuous. Moreover, $B$ induces the
linear map $T_B:E\to F^*$, $T_B(e)(f)= B(e,f)$ which satisfies
 \[  \noo B\rrm_{jcb}\lel \| B:E\wet F\to \cz\| \lel \noo T_B:E\to F^*\rrm_{cb}  \pl .\]
Note that the operator space projective tensor product is indeed
projective, i.e. $E/F\wet X\lel E\wet X/F\wet X$. Moreover, if $N$
is an injective von Neumann algebra and  $E_1$ is completely
isometrically embedded in $E_2$, then we have an isometric
inclusion
\begin{eqnarray}\label{inject}
 N_*\wet E_1 &\subset&  N_*\wet E_2 \pl.
\end{eqnarray}

The direct sum $E\oplus_p F$ of given operator spaces $E$ and $F$
is defined for $p=\8$ and a matrix $[x_{kl}]$ with
$x_{kl}=(e_{kl},f_{kl})$ by
 \[ \noo [x_{kl}]\rrm_{M_m[E\oplus_\8F]}
 \lel \max\{ \noo [e_{kl}]\rrm_{M_m(E)},  \noo
 [f_{kl}]\rrm_{M_m(E)}\} \pl .\]
The operator space $E\oplus_1F$ is defined by its canonical
inclusion in $(E^*\oplus_\8 F^*)^*$. For $1\le p\le \8$, the
operator space structure is given by complex interpolation $
E\oplus_p F \lel [E\oplus_\8 F,E\oplus_1
 F]_{\frac1p}$.

We refer to  \cite{Haa1,Haa2}, \cite{Te}  and \cite{Tak2,Tak3} for
the general theory of noncommutative $L_p$-spaces. Let $\phi$ be a
normal, semifinite faithful weight with modular automorphism group
$\si_t^{\phi}$. For $0<p<\infty$ the space $L_p(N)=L_p(N,\phi)$ is
defined as a subset of (unbounded) operators affiliated with
$N\rtimes_{\si_t^{\phi}}\rz$. To be more precise, we denote by
$\theta_s$ the dual action and $\tau$  the unique semifinite,
normal, faithful trace on $N\rtimes_{\si_t^{\phi}}\rz$ such that
$\tau\circ \theta_s=e^{-s}\tau$. The Haagerup $L_p$-space is
defined by
 \[ L_p(N)  \lel L_{p}(N,\phi) \lel \{ d \pl|\pl  d \mbox{ $\tau$-measurable and
 } \theta_s(d)\lel \exp(-\frac{s}{p})d \} \pl .\]
We have  an operator valued weight $T(x)=\int_{\rz} \theta_s(x)ds
$ from $N\rtimes_{\si_t^{\phi}}\rz$ to $N$. For a functional
$\phi\in N_*$, $\phi\circ T$ defines a density $d_{\phi}\in
L_1(N)$ such that $\phi\circ T(x)=\tau(d_{\phi}x)$. The tracial
functional (different from $\tau$) on $L_1(N)$ is defined by
 \[ tr(d_{\phi})\lel \phi(1) \pl.\]
For an element $x$ in $L_p(N)$ the norm is given by $\noo
x\rrm_p=(tr(|x|^p))^{\frac1p}$. Let $M\subset N$ be a von Neumann
subalgebra with a faithful, normal conditional expectation $E:N\to
M$. Then we have natural inclusion mappings $i_p:L_p(M)\to
L_p(N)$. Indeed, let us assume that $\phi$ is a normal, faithful
state. According to \cite{C2} we have
 \begin{equation}\label{cs}
    \si_t^{\phi} \circ E \lel   E \circ \si_t^{\phi\circ E} \pl
 \end{equation}
for all $t$.  Thus we have a natural inclusion
$M\rtimes_{\si_t^{\phi}}\rz \subset N\rtimes_{\si_t^{\phi\circ
E}}\rz$. The restriction of $\theta_s$ is the corresponding dual
automorphism group and similarly  for $\tau$. This yields an
isometric inclusion $L_p(M)\subset L_p(N)$ (see \cite{JX} for
details).

$L_2(N)$ is a Hilbert space with scalar product $(x,y)=tr(x^*y)$.
We will use the notation $L_2^{s}(N)$ instead of $L_2(N)^s$ for
$s\in \{c,r,oh\}$. In the theory of operator spaces it is
customary to use the `$(i,j)-(i,j)$-duality'
 \[ \langle [x^*_{ij}],[x_{ij}]\rangle \lel
 \summ_{i,j=1}^n x^*_{ij}(x_{ij}) \]
between matrices $[x^*_{ij}]\in S_1^n\wet X^*$ and $[x_{ij}]\in
M_n(X)$. Unfortunately, this is not consistent with the natural
trace $tr$ on $n\times n$-matrices, which corresponds to
 \[ <\langle  [x^*_{ij}],[x_{ij}]\rangle>  \lel
 \summ_{i,j=1}^n x^*_{ij}(x_{ji}) \pl . \]
This forces us to define the \emph{operator space structure on}
$L_1(N)$ by its action on $N^{op}$. Since $N$ and $N^{op}$
coincide as Banach spaces, we may consider $\iota:L_1(N)\to
(N^{op})^*$ defined by
 \[ \iota(d)(y) \lel tr(dy) \lel \phi_{d}(y) \pl .\]
Here $\phi_d$ is the linear functional associated with the density
$d$ in $L_1(N)$. This implies that
\begin{equation}\label{opp} \iota(L_1(N))\lel  N^{op}_* \pl .\end{equation}
If $\phi$ is a  semifinite, normal, faithful  weight, then
$\phi_n=tr_n\ten \phi$ is a semifinite, normal, faithful  weight
on $M_n(N)$. Moreover, $tr_n\ten \tau$ is the unique trace
satisfying $(tr_n \ten \tau)\circ \theta_s=e^{-s}(tr_n\ten \tau)$
and $tr_n\ten tr:L_1(M_n(N),tr\ten \phi)\to \cz$ still yields the
evaluation at $1$. Therefore, we get
 \begin{align*}
 \noo [\iota(x_{ij})] \rrm_{S_1^n\wet N^{op}_*}
 &=  \sup_{\noo [y_{ij}]\rrm_{M_n(N^{op})}\le 1}
 \bet \summ_{ij=1}^n \iota(x_{ij})(y_{ij})\rag
  \lel   \sup_{\noo [y_{ij}]\rrm_{M_n(N)}\le 1}
 \bet \summ_{ij=1}^n tr(y_{ji}x_{ij}) \rag \\
  &=  \sup_{\noo [y_{ij}]\rrm_{M_n(N)}\le 1}
 \bet tr_n\ten tr([y_{ij}][x_{ij}]) \rag \lel
 \noo [x_{ij}] \rrm_{L_1(M_n\ten N,tr_n\ten \phi)}
 \pl .
 \end{align*}
The use of $N^{op}$ enables us to `untwist'  the duality bracket
and we have
 \begin{eqnarray}\label{os1}
  S_1^n \wet L_1(N,\phi) &=&  L_1(M_n\ten N,tr_n \ten \phi)
 \end{eqnarray}
Here we distinguish between the predual $N_*^{op}$ and the
concrete realization of $N^{op}_*$ as space of operators
$L_1(N,\phi)$.

Column and row space interchange their roles when combined with
the  projective tensor. The parallel duality send columns to
columns. Therefore the dualities $R_n=C_n^*$, $(L_1(N)\wet
R_n)^*=N\ten_{\min} C_n$ and $(L_1(N)\wet
S_1^n)^*=N\bar{\ten}M_n$ imply that
 \begin{equation} \label{cn}
  \noo (\summ_{k=1}^n x_k^*x_k)^{\frac12}\rrm_1 \lel   \noo \summ_{k=1}^n x_k\ten e_{1,k} \rrm_{L_1(N)\wet S_1^n} \lel
  \noo \summ_{k=1}^n x_k \ten e_k\rrm_{L_1(N)\wet R_n} \pl.
 \end{equation}
We will frequently use this well-known (though surprising) switch
between the $L_1$- and $L_{\infty}$-theory. Given  an arbitrary
Hilbert space $H$, we denote by
 \[ (\pl\p,\pl):(L_1(N)\ten H)\ten (L_1(N)\ten H)\to L_{\frac12}(N)\pl ,\pl  (d_1\ten h_1,d_2\ten
 h_2)=d_1d_2 (h_1,h_2)\]
the vector-valued extension of the scalar product. Then \eqref{cn}
implies that
 \begin{equation} \label{hc}
 \noo (x^*,x)^{\frac12}\rrm_{L_1(N)}\lel \noo x\rrm_{L_1(N)\wet H^r
 } \quad \mbox{and} \quad \noo (x,x^*)^{\frac12}\rrm_{L_1(N)}\lel \noo x\rrm_{L_1(N)\wet H^c}
 \end{equation}
for all $x=\sum_k x_k\ten h_k$ with $x^*=\sum_k x_k^* \ten h_k \in
L_1(N)\ten H$. Instead of the opposite structure  $N^{op}$, we
will often  work with an antilinear duality bracket. More
precisely, the map $\bar{\iota}:N\to \overline{L_1(N,\phi)}^*$
given by
 \begin{equation} \label{os2}
 \bar{\iota}(y)(x) \lel  tr(x^*y) \lel \lb x,y\rb  \pl
 \end{equation}
is a complete isometry.  The symbol $\lb x ,y \rb=tr(x^*y)$ is
reserved for this \emph{antilinear} duality bracket. Let us
illustrate this duality in connection with \eqref{hc}. The dual
space of the first column in $L_1(N\bar{\ten} B(\ell_2))$ with
respect to the antilinear duality bracket is the first column in
$N\bar{\ten} B(\ell_2)$. For such that a column with entries
$(x_k)$, the sum $\sum_k x_k^*x_k$ is converging in the weak
operator topology.  We use the suggestive notation
$B(\ell_2)\bar{\ten}C$ and $B(\ell_2)\bar{\ten}R$. Given an
arbitrary Hilbert spaces $H$, we may still consider columns and
rows  after fixing a unit vector. We deduce from \eqref{hc} that
we have complete isometries
 \begin{equation}\label{adrc}
  \overline{L_1(N)\wet H^r}^*\lel N\bar{\ten}H^c \quad \mbox{and}
  \quad \overline{L_1(N)\wet H^c}^* \lel N\bar{\ten} H^r \pl
  .\end{equation}
In the sequel, we will use
 \begin{equation}\label{ipeq}
  S_2^n[L_2^{oh}(N)]\lel L_2(M_n\ten N) \pl.
 \end{equation}
Equality \eqref{ipeq} follows immediately from \cite{Pip} in the
hyperfinite, semifinite case. In the general case, we may first
assume that $N$ is $\si$-finite. Let $\phi$ be a faithful normal
state with density $d\in L_1(N)$. Then the map $v:N\to L_2(N)$,
$v(x)=d^{\frac14}xd^{\frac14}$ is bounded. Note that $\bar{v}^*v$
is the map $M_{d^{1/2},d^{1/2}}:N\to L_1(N)$ given by
$M_{d^{1/2},d^{1/2}}(x)=d^{\frac12}xd^{\frac12}$. Hence
\cite[Corollary 2.4]{Poh} implies that
 \[ [M_{d^{1/2},d^{1/2}}(N),L_1(N)]_{\frac12} \lel L_2^{oh}(N) \]
completely isometrically. Therefore, we deduce from \cite{Ko} and
\cite[Corollary 1.4]{Pip} that
 \[ L_2(M_n\ten N)\lel [(id\ten M_{d^{1/2},d^{1/2}})(M_n\ten N),L_1(M_n\ten
 N)]_{\frac12}\lel S_2^n[L_2^{oh}(N)] \pl .\]
Since every von Neumann algebra admits a strictly semifinite
normal weight, (i.e. a weight which is an orthogonal sum of
states), \eqref{ipeq} follows by approximation for arbitrary von
Neumann algebras.

Let $N$ be a semifinite von Neumann algebra and $\tau$  a normal,
faithful, semifinite trace. Then the classical $L_p$-space
$L_p(N,\tau)=[N,L_1(N^{op},\tau)]_{\frac1p}$ is (completely)
isometrically isomorphic to  the Haagerup $L_p$-space
$L_p(N,\tau)$ (see  \cite{Te} for the explicit isomorphism).
Moreover, \eqref{os1} and \eqref{os3} below  also hold in the
category of semifinite $L_p$-spaces.

An important result of Effros and Ruan (see e.g. \cite{ER}) shows
that the projective tensor product is compatible with von Neumann
algebras. Indeed, given von Neumann algebras $N$ and $M$ then
 \begin{eqnarray} \label{ER}
 (N_*\wet M_*)^* &=&  N\bar{\ten}M \pl .
 \end{eqnarray}
Indeed, in the $\si$-finite with n.s.f. states $\phi$ and $\psi$,
the modular group of $\phi\ten \psi$ is given by $\si_t^{\phi\ten
\psi}\lel \si_t^{\phi}\ten \si_t^{\psi}$. We may then use
$L_1(N)\cong N^{op}_*$, $L_1(M)\cong M^{op}_*$.  This yields a
completely isometric isomorphism
 \begin{eqnarray}  \label{os3}
  L_1(N\bar{\ten}M,\phi\ten \psi) &\cong & (N\bar{\ten}M)^{op}_*
  \lel N^{op}_*\wet M^{op}_*  \pl \cong \pl  L_1(N,\phi)\wet
  L_1(M,\psi)\pl .
 \end{eqnarray}
The general case follows by approximation from the $\si$-finite
case.


\section{A logarithmic characterization of C$^*$-algebras with
WEP}

The notion of $(2,oh)$-summing maps introduced in \cite{Poh} will
be an important tool in this section. It allows us to find the
smallest integral, sesquilinear form majorizing a given positive,
sesquilinear form. A linear map $u:E\to OH$ is called
\emph{$(2,oh)$-summing} if there exists a constant $C>0$ such that
 \begin{equation} \label{2oh}  \summ_k \noo u(x_k)\rrm^2\kl C^2  \noo \summ_{k} \bar{x}_k\ten
 x_k\rrm_{\bar{E}\ten_{\min}E}
 \end{equation}
holds for all finite sequences $(x_k)\subset E$. The
$(2,oh)$-summing norm is defined by $\pi_2^{oh}(u)\lel \inf C$,
where the infimum is taken over all $C$ satisfying \eqref{2oh}.\lz

\begin{samepage}\begin{lemma}\label{transl} Let $E$ be an operator space and
$B:\bar{E}\times E\to \cz$ be a positive sesquilinear form. Let
$L=\{x\in E: B(x,x)=0\}$ and $K$ the completion of $E/L$ with
respect to the induced scalar product $(x+L,y+L)=B(x,y)$. Then
the linear map $u:E\to K^{oh}$, $u(x)=x+L$  satisfies
 \begin{enumerate}
  \item[i)] $B(x,y)=(u(x),u(y))$ for all $x,y\in E$,
 \item[ii)] $\noo u\rrm_{cb}^2 \lel \noo B\rrm_{jcb}$,
 \item[iii)] $\inf_{B\le \tilde{B}} \noo \tilde{B}\rrm_{I}\lel
 \pi_2^{oh}(u)^2$.
  \end{enumerate}
\end{lemma}\end{samepage}

\begin{proof}[\bf Proof:] Equality i) is obvious from the definition of $u$.
Assertion ii) is proved in \cite[Corollary 2.4]{Poh}. For the
estimate $"\ge"$ in iii), we assume that $\tilde{B}$ is an
integral form such that $B\le \tilde{B}$. Let $z_k$ be scalars
such that $|z_k|=1$ and
$z_k\tilde{B}(x_k,x_k)=|\tilde{B}(x_k,x_k)|$. Then, we deduce
from a version of the Cauchy-Schwarz inequality due to Haagerup
(see \cite[(7.2)]{Ps5}) that
 \begin{align*}
 \summ_k B(x_k,x_k)&\le  \summ_k |\tilde{B}(x_k,x_k)|
 \lel   \bet \summ_k \tilde{B}(x_k,z_kx_k) \rag \kl
 \|\tilde{B}\|_I \noo \summ_k \bar{x}_k\ten
 z_kx_k\rrm_{\min} \\
 &\le \|\tilde{B}\|_I \noo \summ_k \bar{x}_k\ten x_k\rrm_{\min}^{\frac12}
 \noo \summ_k \bar{z}_k\bar{x}_k\ten z_kx_k\rrm_{\min}^{\frac12}
 \lel  \|\tilde{B}\|_I  \noo \summ_{k} \bar{x}_k\ten
 x_k\rrm_{\min} \pl .
 \end{align*}
holds for all $(x_k)\subset E$. Taking the infimum over all
$\tilde{B}\ge B$, we obtain
 \[  \summ_k \noo u(x_k)\rrm^2\lel
 \summ_k B(x_k,x_k) \kl \inf_{B\le \tilde{B}}
 \|\tilde{B}\|_I  \noo \summ_{k} \bar{x}_k\ten
  x_k\rrm_{\min} \pl. \]
This implies that $\pi_2^{oh}(u)^2\le \inf_{B\le \tilde{B}}
 \|\tilde{B}\|_I$. Conversely, we assume \eqref{2oh}. We apply
a variant of the Grothendieck-Pietsch separation argument in the
context of $(2,oh)$-summing maps, see \cite[Theorem 5.7]{Poh}.
For this we shall assume that the operator space $E$ is given by
a concrete representation $E\subset \B(H)$. Then there exists an
index set $I$, an ultrafilter $\U$ and nets $(a_i), (b_i)$ in the
unit sphere of $S_4(H)$ such that
 \begin{equation} \label{domm}
   \noo u(x)\rrm \kl \pi_2^{oh}(u) \pl \lim_{i,\U} \noo a_ixb_i\rrm_2 \pl
 \end{equation}
holds for all $x\in E$. For fixed $i$,   we note that
 \begin{align*}
  \noo a_ixb_i\rrm_2^2 &= tr(a_ixb_ib_i^*x^*a_i^*) \lel
  tr(x^*a_i^*a_ixb_ib_i^*) \lel
  ((b_ib_i^*)^{t},(x^*\ten x)(a_i^*a_i))
   \pl .
  \end{align*}
We refer to \cite{P,ER} for the fact that $B_i(x,y)\lel
tr(x^*a_i^*a_iyb_ib_i^*) \lel ((b_ib_i^*)^{t},(x^*\ten
y)(a_i^*a_i))$ is nuclear with norm $\le 1$. Integral forms are
closed under  pointwise limits with  respect to  bounded nets of
nuclear forms. Therefore
 \[ \tilde{B}(x,y)\lel \pi_2^{oh}(u)^2\pl \lim_i tr(x^*a_i^*a_iyb_ib_i^*) \]
is a positive  integral bilinear form with $ \|\tilde{B}\|_I\le
\pi_2^{oh}(u)^2$. Inequality  \eqref{domm} shows that
$B(x,x)=(u(x),u(x))$ is dominated by $\tilde{B}$. \qd

\begin{lemma}\label{triv}   Let $A$ be a $C^*$-algebra with WEP and $u:A\to
OH_n$. Then
 \[ \pi_2^{oh}(u) \kl C(1+\ln n)^{\frac12} \noo u\rrm_{cb} \pl .\]
\end{lemma}\lz

\begin{proof}[\bf Proof:] We recall that $A$ has WEP, if there exists
a contraction $P:\B(H)\to A^{**}$ such that $A\subset
A^{**}\subset \B(H) $ and that $P|_A=id_A$.  It is well-known that
$P$ is indeed completely contractive (see e.g. \cite[Lemma
2.1.]{Fub}). Then $v=u^{**}P:\B(H) \to OH_n$ satisfies
 \[ \noo v\rrm_{cb}\kl \noo P\rrm_{cb} \pl \noo u^{**}\rrm_{cb}
 \kl \noo u\rrm_{cb} \pl .\]
We can apply \eqref{uppb} for $B(x,y)=(v(x),v(y))$ and deduce the
assertion from Lemma \ref{transl}.\qd

\begin{lemma}\label{cbohne} Let $N$ be a von Neumann algebra.
Let  $a,b\in L_4(N)$ and  $M_{ab}:N\to L_2^{oh}(N)$ defined by
$M_{ab}(x)\lel axb$. Then
 \[ \noo M_{ab}:N\to L_2^{oh}(N)\rrm_{cb} \kl \noo a\rrm_4 \pl
 \noo b\rrm_4 \pl .\]
\end{lemma}\lz

\begin{proof}[\bf Proof:] Let $x\in M_n(N)$ and $a,b\in L_4(N)$.
According to \cite[Lemma 1.7]{Pip}, we have
 \[ \noo (id\ten M_{ab}) (x) \rrm_{M_n(L_2^{oh}(N))} \lel \sup_{\al,\beta} \noo
 (\al\ten a)x(\beta \ten b) \rrm_{S_2^n[L_2^{oh}(N)]} \pl .\]
Here the supremum is taken over all $\al$, $\beta$ in the unit
ball of $S_4^n$. The assertion follows from \eqref{ipeq} and
H\"older's inequality
 \begin{align*}
 \noo (\al\ten a)x(\beta \ten b) \rrm_{S_2^n[L_2^{oh}(N)]}  &=
 \noo (\al\ten a)x(\beta \ten b) \rrm_{L_2(M_n\ten N)} \kl
  \noo \al\ten a\rrm_4
 \noo x\rrm_{M_n(N)} \noo \beta \ten b\rrm_4 \\
 &\le   \noo \al \rrm_4
  \noo a\rrm_4 \pl  \noo x\rrm_{M_n(N)} \pl \noo \beta\rrm_4 \noo b\rrm_4  \pl . \qedhere
 \end{align*}
\qd

\begin{proof}[\bf Proof of Theorem 6:]
We use the one-to-one correspondence between completely bounded
linear  maps $u:A\to OH$ and  positive sesquilinear jcb-forms on
$A$ (see Lemma  \ref{transl}). Theorem 6 means  that $A$ has WEP
if and only if
 \begin{equation}\label{oko}
  \pi_2^{oh}(u)\le C(1+\ln n)^{\beta} \noo u\rrm_{cb}
  \end{equation}
holds for every linear map $u:A\to OH_n$. If $A$ has WEP, then
Lemma \ref{triv} implies that \eqref{oko} holds for
$\beta=\frac12$ and a universal constant $C$. Conversely, we
assume that \eqref{oko} holds for some $\beta>0$ and some constant
$C$. Let us consider the von Neumann algebra $N=A^{**}$. Recall
from \cite{Te} that $N$ is in standard form on $L_2(N)$. This
means in particular that $N$ acts on $L_2(N)$ by left
multiplication $\pi(x)h=xh$ and $J(h)=h^*$ is an antilinear
isometry $J$ such that $N'\lel JNJ$ (see \cite{Te}). Let $h\in
L_2(N)$ be a unit vector. Then, we can find $a,b\in L_4(N)$ of
norm $1$ such that $h=ab$. According to Lemma \ref{cbohne}, the
maps $M_{a^*a}:N\to L_2^{oh}(N)$ and $M_{bb^*}:N\to L_2^{oh}(N)$
are complete contractions. Now, we consider elements $x_1,...,x_n$
in $A$. Let $P$ and $Q$ be orthogonal  projections onto ${\rm
span}\{M_{a^*a}(x_k)| 1\le k\le n\}$ and ${\rm
span}\{M_{bb^*}(x_k)| 1\le k\le n\}$, respectively. Then
$PM_{a^*a}$ and $QM_{bb^*}$ have rank as most $n$. Therefore we
may apply \eqref{oko} and deduce
 \begin{eqnarray}
 \begin{minipage}{15cm}\vspace{-0.3cm}
 \for
 & &(h,\summ_{k=1}^n x_kJx_kJh) \lel      (h,\summ_{k=1}^n x_khx_k^*) \lel
  \summ_{k=1}^n tr(h^*x_khx_k^*) \lel \summ_{k=1}^n tr(a^*x_kabx_k^*b^*) \nonumber \\
 & & \pll \le \kla \summ_{k=1}^n \noo P(a^*x_ka)\rrm^2 \mer^{\frac12}
 \kla \summ_{k=1}^n \noo Q(bx_kb^*)\rrm^2 \mer^{\frac12} \kl
 C^2(1+\ln n)^{2\beta} \noo \summ_{k=1}^n \bar{x}_k\ten
 x_k\rrm_{\min} \pl . \nonumber
 \mel
 \end{minipage}
 \end{eqnarray}
In order to eliminate the log-term, we use Haagerup's trick and
consider the positive  operator \for
 \lefteqn{ [(\summ_{k=1}^n x_kJx_kJ)^*(\summ_{k=1}^n
 x_kJx_kJ)]^m}\\
   & & \lel  \summ_{k_1,...,k_{2m}=1}^n   x_{k_1}^*(Jx_{k_1}^*J)x_{k_2}(Jx_{k_2}J) \cdots
     \cdots     x_{k_{2m-1}}^*(Jx_{k_{2m-1}}^*J)      x_{k_{2m}}(Jx_{k_{2m-1}}J) \\
   & & \lel  \summ_{k_1,...,k_{2m}=1}^n (Jx_{k_1}^*x_{k_2}\cdots
   x_{k_{2m-1}}^*x_{k_{2m}}
   J)   x_{k_1}^*x_{k_2}\cdots
   x_{k_{2m-1}}^*x_{k_{2m}}   \\
      \mel
We apply $(2.4)$ to  the finite family
$x_{k_1,....,k_{2m}}=x_{k_1}^*x_{k_2}\cdots
   x_{k_{2m-1}}^*x_{k_{2m}}$  and deduce that
\begin{align*}
 &(h,[(\summ_{k=1}^n x_kJx_kJ)^*(\summ_{k=1}^n
 x_kJx_kJ)]^mh)\\
 &=
  \summ_{k_1,...,k_{2m}=1}^n (h,x_{k_1}^*x_{k_2}\cdots
   x_{k_{2m-1}}^*x_{k_{2m}}Jx_{k_1}^*x_{k_2}\cdots
   x_{k_{2m-1}}^*x_{k_{2m}}
   Jh) \\
 &\pl \le C^2 (1+\ln n^{2m})^{2\beta} \pl \noo \summ_{k_1,....,k_{2m}=1}^n
 \bar{x}_{k_1,....,k_{2m}} \ten x_{k_1,....,k_{2m}}
 \rrm_{\min} \\
 &\pl = \! C^2 (1+\ln n^{2m})^{2\beta} \p \noo [(\summ_{k=1}^n \bar{x}_k\ten
 x_k)^* (\summ_{k=1}^n \bar{x}_k\ten x_k)]^m \rrm_{\min} \!\! \lel
 C^2 (1+2m \ln n)^{2\beta} \noo \summ_{k=1}^n \bar{x}_k\ten
 x_k\rrm^{2m}_{\min} \! .
 \end{align*}
Taking the supremum over $\noo h\rrm\le 1$, we deduce from
positivity that
 \begin{align*}
  \noo \summ_{k=1}^n x_kJx_kJ \rrm
  &=  \noo [(\summ_{k=1}^n x_kJx_kJ)^*(\summ_{k=1}^n
 x_kJx_kJ)]^m\rrm^{\frac{1}{2m}} \nonumber \\
 &\le  C^{\frac2m}
 (1+2m)^{\frac{2\beta}{2m}} (1+ \ln n)^{\frac{2\beta}{m}} \pl
  \noo \summ_{k=1}^n \bar{x}_k\ten
 x_k\rrm_{\min} \pl . \end{align*}
Taking the limit for $m\to \8$, we obtain  (in the language of
\cite{P3})
 \begin{align}
 \noo \summ_{k=1}^n L_{x_k}R_{x_k^*}\rrm &=
  \noo \summ_{k=1}^n x_kJx_kJ \rrm \kl \noo \summ_{k=1}^n \bar{x}_k\ten
 x_k\rrm_{\bar{A}\ten_{\min}A} \lel
  \noo \summ_{k=1}^n \bar{x}_k\ten
  x_k\rrm_{\overline{\B(H) }\ten_{\min}\B(H)  } \pl . \label{previous}
 \end{align}
Let us recall an equality proved by  Pisier \cite[Theorem 2.1]{P3}
(and \cite[Theorem 2.9]{Haa} in the non-semifinite case)
 \begin{eqnarray}
 \noo \summ_{k=1}^n \bar{x}_k\ten x_k\rrm_{\bar{A}\ten_{\max} A}
 &=&
 \noo  (x_1,...,x_n)\rrm_{[R_n(A),C_n(A)]_{\frac12}}^2 \lel
   \noo \summ_{k=1}^n L_{x_k}R_{x_k^*}\rrm_{L_2(A^{**})}  \pl .
 \end{eqnarray}
Combining this with  \eqref{previous},  we deduce that
 \begin{align*}
 \noo \summ_{k=1}^n \bar{x}_k \ten x_k\rrm_{\bar{A}\ten_{\max} A}
 &\le     \noo \summ_{k=1}^n \bar{x}_k\ten
 x_k\rrm_{\overline{\B(H) }\ten_{\min} \B(H) }
  \kl \noo \summ_{k=1}^n \bar{x}_k\ten
  x_k\rrm_{\overline{\B(H) }\ten_{\max} \B(H) }
  \pl .
 \end{align*}
According to \cite[Theorem 3.7]{Haa} $A$ has WEP. In section 4 we
will show  that  $\al=2\beta=1$ is the best possible exponent for
non sub-homogeneous $C^*$-algebras. \qd

\begin{rem} {\rm Let us consider the special case where $A=N$ is
a von Neumann algebra. We have shown above that if $N$ satisfies
the logarithmic Grothendieck inequality, then
  \[ \noo (x_1,...,x_n)\rrm_{[R_n(N),C_n(N)]_{\frac12}} \kl \noo
   \summ_{k=1}^n \bar{x}_k \ten x_k
   \rrm_{\bar{N}\ten_{\min}N}^{\frac12}\pl .
   \]
According to Pisier's characterization \cite{P3} (see
\cite[Theorem 2.9]{Poh} for a simple proof in the semifinite
case), we deduce that $N$ is injective. Thus for von Neumann
algebras, we don't have to use Haagerup's deep (and unfortunately
unpublished) results.}
\end{rem}

\section{The Pusz/Woronowicz formula and the operator space OH}

In this section we will show the connection between the
Pusz/Woronowicz formula for square roots of sesquilinear forms and
subspaces of quotients of Hilbert spaces. For our understanding of
the problem this concrete formula (and its dual version) plays an
important role. Following \cite{Wo} we consider two  positive
\emph{commuting} operators $A,B$ on a Hilbert space. According to
\cite{Wo}
 \begin{align}\label{wo}
    ((AB)^{\frac12}x,x) &=
   \inf_{x=a(t)+b(t)} \intt_0^1 \frac{(Aa(t),a(t))}{t}+
  \frac{(Bb(t),b(t))}{1-t} \pl \frac{dt}{\pi \sqrt{t(1-t)}}
  \pl .
 \end{align}
Here the infimum is taken over piecewise constant functions in
$H$; see \cite[Appendix]{Wo} for a proof. Let us denote by
$H_{\sqrt{AB}}$ the Hilbert space $H$ equipped with the scalar
product $\sqrt{AB}(x,y)=(\sqrt{AB}x,y)$. Motivated by \eqref{wo},
we define the probability  $\mu$ and the measures $\nu_1$,
$\nu_2$ on $[0,1]$ as follows:
 \begin{equation}\label{mumeas}  d\mu(t)\lel \frac{dt}{\pi  \sqrt{t(1-t)}} \quad ,\quad
 \nu_1(t)\lel t^{-1}d\mu(t)\quad \mbox{and} \quad
 d\nu_2(t) \lel (1-t)^{-1} d\mu(t) \pl.
 \end{equation}
We denote by $H_{A}$ and  $H_{B}$ the space $H$ equipped with the
Hilbertian norm $\noo x\rrm_{H_A}=(Ax,x)^{\frac12}$ and $\noo
x\rrm_{H_B}=(Bx,x)^{\frac12}$, respectively. If $A$ and $B$ are
invertible, the canonical inclusion maps $H_{A}\subset H$ and
$H_{B}\subset H$ are continuous. Then we may  define the linear
map $Q:L_2(\nu_1,H_A)\oplus_2 L_2(\nu_2,H_B)\to L_0(\mu;H)$ by
$Q(f,g)(t)\lel f(t)+g(t)\in H$. We denote by
 \[ K \lel L_2(\nu_1,H_A) \oplus_2 L_2(\nu_2,H_B)/\krn(Q) \]
the quotient space and
 \[ E\lel \{(f,g)+\krn(Q) \pl|\pl Q(f,g) \mbox{ constant a.e.}
 \}  \pl . \]
In this terminology the Pusz/Woronowicz formula reads as follows:
\lz

\begin{lemma} \label{W1}
If $A$ and $B$ are invertible, then  $H_{\sqrt{AB}}$ is
isometrically isomorphic to the subspace $E\subset K$.
\end{lemma}\lz

The dual version of Lemma \ref{W1} is based on a characterization
of linear functionals on $E$.\lz

\begin{lemma}\label{lf}
Let $A$ and $B$ be bounded and  invertible. Linear functionals
$\phi:E\to \cz$ are in one-to-one correspondence with pairs
$(f,g)\in L_2(\nu_1,H_A)\oplus L_2(\nu_2,H_B)$ such that
 \begin{equation} \label{caracf}  \frac{Af(t)}{t} \lel \frac{Bg(t)}{1-t} \quad \mu \mbox{ a.e.}
 \quad,
 \quad
  \phi(x) \lel \intt_0^1 (\frac{Af(t)}{t},x) \pl
 d\mu(t) \pl .\end{equation}
Moreover,
 \[ \noo \phi\rrm \lel \inf (\noo
f\rrm_{L_2(\nu_1,H_A)}^2 +\noo
 g\rrm_{L_2(\nu_2,H_B)}^2)^{\frac12} \pl .\]
Here the infimum is taken over all pairs satisfying
\eqref{caracf}.
\end{lemma}\lz

\begin{proof}[\bf Proof:] By the Hahn-Banach theorem, the norm
one functionals on $E$ are in one to one correspondence with the
restrictions $\hat{\phi}|_E$ of norm one functionals
$\hat{\phi}:L_2(\nu_1,H_A)\oplus_2 L_2(\nu_2,H_B)\to \cz$.  A
norm one functional $\hat{\phi}$ is given by a couple $(f,g)$
such that $\| \hat{\phi}\|^2=\noo f\rrm_{L_2(\nu_1,H_A)}^2+\noo
g\rrm_{L_2(\nu_2,H_B)}^2$ and
 \begin{align*}
  \hat{\phi}(h_1,h_2)&=
  \intt_0^1 (Af(t),h_1(t)) \frac{d\mu(t)}{t} +
  \intt_0^1 (Bg(t),h_2(t)) \frac{d\mu(t)}{1-t} \pl .
  \end{align*}
Thus the functional $\hat{\phi}$ vanishes for all $(h,-h)$ if and
only if $\frac{Af(t)}{t}=\frac{Bg(t)}{1-t}$ $\mu$-almost
everywhere. Finally, given $x\in E$, we see that $\sqrt{t}x\in
L_2(\nu_1,H_A)$ and  $(1-\sqrt{t})x\in L_2(\nu_1,H_A)$. Since
$\hat{\phi}$ is an extension of $\phi$, we deduce that
 \begin{align*}
  &\phi(x)\lel \hat{\phi}Q(\sqrt{t}x,\sqrt{1-t}x) \\
  &\pl \lel  \intt_0^1 (\frac{Af(t)}{t},\sqrt{t}x)+
  (\frac{Bg(t)}{1-t},(1-\sqrt{t})x) d\mu(t) \lel
  \intt_0^1
  (\frac{Af(t)}{t},x) d\mu(t)\pl  . \qedhere
 \end{align*}
\qd

\begin{samepage}\begin{lemma} \label{W2} Let $A,A^{-1}$, $B,B^{-1}$ be
bounded and $y\in H$. Then
 \[ ((AB)^{\frac12}y,y) \lel \inf \intt_0^1 \frac{(Af(t),f(t))}{t} d\mu(t)  +
  \intt_0^1 \frac{(Bg(t),g(t))}{1-t} d\mu(t) \pl ,\]
where the infimum is taken over all tuples $(f,g)$ of $H$-valued
measurable functions satisfying
$\frac{Af(t)}{t}=\frac{Bg(t)}{1-t}$ $\mu$ a.e. and
 \[ B^{\frac12}y \lel \intt_0^1 \frac{A^{\frac12}f(t)}{t} \pl d\mu(t) \pl .\]
\end{lemma}\lz\end{samepage}

\begin{proof}[\bf Proof:] Since $\sqrt{AB}$ defines a scalar product, we have
 \[ \noo y \rrm_{H_{\sqrt{AB}}}  \lel ((AB)^{\frac12}y,y)^{\frac12}
 \lel \sup_{\sqrt{AB}(x,x)\le 1} |((AB)^{\frac12}y,x)|    \pl .\]
Thus the norm of $y$ in $H_{\sqrt{AB}}$ coincides with the norm
of the linear functional
 \[ \phi_{y}(x)\lel ((AB)^{\frac12}y,x) \pl \]
on $H_{\sqrt{AB}}$. According to Lemma \ref{lf}, we  can find
$(f,g)$ such that \eqref{caracf} is satisfied and
 \for
   ((AB)^{\frac12} y,y) &=&  \noo \phi_{y} \rrm^2
    \lel
  \intt_0^1  (Af(t),f(t)) \frac{d\mu(t)}{t} +
   \intt_0^1  (Bg(t),g(t)) \frac{d\mu(t)}{1-t} \pl .
   \mel
From \eqref{caracf} and the definition of $\phi_y$ we deduce that
 \[ B^{\frac12}y \lel A^{-\frac12} (AB)^{\frac12}y
 \lel A^{-\frac12} \intt_0^1 \frac{Af(t)}{t} \pl
  d\mu(t) \lel
  \int_0^1 \frac{A^{\frac12}f(t)}{t} \pl d\mu(t) \pl .\]
Conversely, any pair satisfying these conditions induces the same
functional $\phi_{y}$ and thus provides an upper estimate for the
norm of $y$. \qd

As pointed out in the introduction, the operator space OH$_n$
will be obtained by amplifications of these densities. To be more
specific, we consider the map $Q:L_2(\nu_1;\ell_2^n)\oplus_1
L_2(\nu_2;\ell_2^n)\to L_0(\mu;\ell_2^n)$ defined by
 \[ Q(f,g)(t) \lel f(t)+g(t)\in \ell_2^n  \pl .\]
We denote by $G_n$ the quotient space
 \[ G_n \lel  L_2^c(\nu_1;\ell_2^n) \oplus_1
 L_2^r(\nu_2;\ell_2^n)/\krn(Q)\pl  .\]
In the limiting case $n=\infty$ we will simply write $G$. The
interesting subspaces $F_n$, $F$ of $G_n$  are given by those
equivalence classes $(f,g)+\krn Q$  such that $Q(f,g)$ is a
constant function with values in $\ell_2^n$, $\ell_2$. We denote
by $f_1,..,f_n$ the canonical unit vector basis in $F_n$ given by
$f_k=(\sqrt{t}e_k,(1-\sqrt{t})e_k)+\krn(Q)$. The following lemma
is proved using the polar decomposition and the density of
invertible matrices.\lz

\begin{lemma} \label{t1}  Let $x_1,...,x_n\in M_m$. Then
\[ \noo \summ_{k=1}^n \bar{x}_k \ten x_k \rrm_{\overline{M_m}\ten_{\min} M_m}^{\frac12}
 \kl  \sup_{\noo a\rrm_4\le 1, \noo b\rrm_4\le 1, a> 0, b>0} \kla \summ_{k=1}^n \noo bx_ka\rrm_2^2 \mer^{\frac12} \pl. \]
Here $a>0$ means that $a\ge 0$ and $a$ is invertible.
\end{lemma} \lz

\begin{lemma} \label{1s} Let $(e_k)$ be the natural unit vector basis of OH$_n$.
The identity map $id:F_n\to $OH$_n$ given by $(f_k)=e_k$ has
$cb$-norm less than $\sqrt{2}$.
\end{lemma} \lz

\begin{proof}[\bf Proof:] Let $x_1,...,x_n\in M_m$ and assume
$ \noo \sum_{k=1}^n x_k\ten f_k\rrm_{M_m(F_n)}\pl <\pl 1$. By the
definition of the matrix norms for quotient operator spaces, we
find elements $f\in M_m(L_2^c(\nu_1;\ell_2^n))$ and $g\in
M_m(L_2^r(\nu_2;\ell_2^n))$ such that
 \[ x_k \lel f_k(t) + g_k(t) \quad \mbox{$\mu$-a.e. }\]
and
\[
   \max\left\{\noo \intt_0^1 \summ_{k=1}^n f_k(t)^*f_k(t) \frac{d\mu(t)}{t}
 \rrm_{M_m}^{\frac12},
   \noo \intt_0^1 \summ_{k=1}^n g_k(t)g_k(t)^* \frac{d\mu(t)}{1-t} \rrm_{M_m}^{\frac12}\right\}  \kl 1 \pl .\]
Let $a,b$ be  positive, invertible, norm one elements in $S_4^m$.
On the Hilbert space $H=\ell_2^n(S_2^m)$ with the scalar product
 \[ ((x_k),(y_k)) \lel \summ_{k=1}^n tr(x_k^*y_k) \]
we  define $A(x_k)=(x_ka^4)$ and $B(x_k)=(b^4x_k)$. Clearly, these
operators commute and we deduce from $(3.1)$ that

 \begin{align*}
 \summ_{k=1}^n \noo b x_ka\rrm_2^2
 &= \summ_{k=1}^n tr(a^*x_k^*b^*bx_ka)
 \lel  \summ_{k=1}^n tr(a^2x_k^*b^2x_k) \lel
  ( (AB)^{\frac12}(x_k),(x_k))\\
 &\le \intt_0^1 (A(f_k(t)),(f_k(t))) \frac{d\mu(t)}{t} +
 \intt_0^1 (B(g_k(t)),(g_k(t))) \frac{d\mu(t)}{1-t} \\
 &=  \intt_0^1 \summ_{k=1}^n tr(a^4f_k^*(t)f_k(t)) \frac{d\mu(t)}{t} +
 \intt_0^1 \summ_{k=1}^n tr(g_k(t)^*b^4g_k(t)) \frac{d\mu(t)}{1-t} \\
 &= tr\Bigg(a^4 \intt_0^1 \summ_{k=1}^n f_k^*(t)f_k(t) \frac{d\mu(t)}{t}\Bigg) +
 tr\Bigg(b^4\intt_0^1 \summ_{k=1}^n g_k(t)g_k(t)^*
 \frac{d\mu(t)}{1-t}\Bigg) \\
 &\le   \noo a^4\rrm_1 \noo
 \intt_0^1 \summ_{k=1}^n f_k^*(t)f_k(t) \frac{d\mu(t)}{t}\rrm_{M_m}
 +
 \noo b^4\rrm_1  \noo \intt_0^1 \summ_{k=1}^n g_k(t)g_k(t)^*\frac{d\mu(t)}{1-t}\rrm_{M_m}
  \!\! \le   \! 2 \p .
 \end{align*}
From Lemma \ref{t1}, we deduce that $\noo id: F_n\to
OH_n\rrm_{cb} \le\sqrt{2}$. \qd

\begin{lemma} \label{2s}
Let $(f_k^*)_{k=1}^n$ be the dual basis of $F_n^*$ satisfying
$f_k^*(f_j)=\delta_{kj}$. The identity map  $id: F_n^*\to OH_n$
given by $id(f_k^*)=e_k$ has $cb$-norm less than $\sqrt{2}$.
\end{lemma} \lz

\begin{proof}[\bf Proof:] We have to consider a norm one element  $z\in
M_m(F_n^*)\cong CB(F_n,M_m)$. Then, we may write $z=\sum_{k=1}^n
z_k\ten f_k^*$. Let us denote the corresponding complete
contraction which satisfies $z_k=u_z(f_k)\in M_m$ as $u_z:F_n\to
M_m$. According to Wittstock's theorem there exists a complete
contraction $v:G_n\to M_m$. Since $G_n$ is a quotient space, we
see that $vQ:L_2^c(\nu_1;\ell_2^n)\oplus_1
L_2^r(\nu_2;\ell_2^n)\to M_m$ is a complete contraction. Thus
there are $x\in M_m(L_2^r(\nu_1;\ell_2^n))$ and $y\in
M_m(L_2^c(\nu_2;\ell_2^n))$ of norm less than $1$ such that
 \[ vQ(h^1,h^2) \lel \intt_0^1 \summ_{k=1}^n x_k(t)h^1_k(t) \frac{d\mu(t)}{t}
 + \intt_0^1 \summ_{k=1}^n y_k(t)h^2_k(t) \frac{d\mu(t)}{1-t}  \pl .\]
Again, we can use the fact that $v$ vanishes on $\krn(Q)$ and get
 \[   \frac{x_k(t)}{t} \lel \frac{y_k(t)}{1-t} \quad \mbox{$\mu$-a.e.}\]
for all $k=1,...,n$. In order to identify our original map $u_z$,
we compute
 \begin{align*}
  u_z(f_k) = vQ(e_k\sqrt{t},e_k(1-\sqrt{t})) =
  \intt_0^1 \sqrt{t} x_k(t) \pl \frac{d\mu(t)}{t} +
  \intt_0^1 (1-\sqrt{t}) y_k(t) \pl \frac{d\mu(t)}{1-t} =
   \intt_0^1 x_k(t) \pl \frac{d\mu(t)}{t} \pl .
  \end{align*}
Now, let us consider invertible positive elements $a,b\in S_4^m$.
As above, we define the operators
$A(z_k)_{k=1}^n=(z_ka^4)_{k=1}^n$ and
$B(z_k)_{k=1}^n=(b^4z_k)_{k=1}^n$. We consider
$\tilde{x}_k(t)=b^2x_k(t)a^{-2}$ and get
 \for
  B^{\frac12}((z_k)_{k=1}^n) &=& (b^2z_k)_{k=1}^n \lel \intt_0^1 (b^2x_k(t))_{k=1}^n \pl  \frac{d\mu(t)}{t}
  \lel \intt_0^1 A^{\frac12}((\tilde{x}_k(t))_{k=1}^n)  \pl
\frac{d\mu(t)}{t}  \pl .
  \mel
On the other hand  for $\tilde{x}=(\tilde{x}_k)_{k=1}^n$, we
deduce that
 \begin{align*}
 \intt_0^1 (A \tilde{x}(t),\tilde{x}(t)) \frac{d\mu(t)}{t}
 &=  \intt_0^1 (\tilde{x}(t),A \tilde{x}(t))  \frac{d\mu(t)}{t}
 \lel
  \intt_0^1 \summ_{k=1}^n tr(\tilde{x}_k^*(t) \tilde{x}_k(t)a^4) \frac{d\mu(t)}{t} \\
 &=  \intt_0^1 \summ_{k=1}^n tr(a^{-2}x_k^*(t)b^2b^2 x_k(t)a^{-2}a^4) \frac{d\mu(t)}{t}\\
 &=  tr\bigg(b^4\intt_0^1 \summ_{k=1}^n x_k(t)x_k^*(t) \pl \frac{d\mu(t)}{t}\bigg) \\
 &\le \noo b^4\rrm_1 \pl \noo \intt_0^1 \summ_{k=1}^n x_k(t)x_k^*(t) \pl
 \frac{d\mu(t)}{t}\rrm_{M_m}\\
 &\le  \noo x\rrm_{M_m(L_2^r(\nu_1;\ell_2^n))}^2 \kl 1 \pl .
 \end{align*}
Similarly, we define $\tilde{y}_k(t)=b^{-2}y_k(t)a^2$,
$\tilde{y}(t)=(\tilde{y}_k(t))_{k=1}^n$ and get
 \for
  \frac{A((\tilde{x}_k(t))_{k=1}^n)}{t} &=& \frac{(b^2x_k(t)a^2)_{k=1}^n}{t}
  \lel \frac{(b^2y_k(t)a^2)_{k=1}^n}{1-t}
  \lel \frac{B((\tilde{y}_k(t))_{k=1}^n)}{1-t}  \quad \mbox{$\mu$-a.e.} \pl .
  \mel
The same calculation as above yields
 \[  \intt_0^1 (B \tilde{y}(t),\tilde{y}(t)) \frac{d\mu(t)}{1-t}
 \kl \noo a^4\rrm_1  \noo y\rrm_{M_m(L_2^c(\nu_1;\ell_2^n))}^2 \kl 1 \pl . \]
Therefore Lemma \ref{W2} implies that
  \for
 \summ_{k=1}^n \noo bz_ka\rrm_2^2 &=& \noo (\sqrt{AB}(x_k)_{k=1}^n)\rrm_2^2
 \kl  \intt_0^1 (A\tilde{x}(t),\tilde{x}(t)) \frac{d\mu(t)}{t} +
 \intt_0^1 (B\tilde{y}(t),\tilde{y}(t)) \frac{d\mu(t)}{1-t} \kl 2 \pl .
 \mel
Since, $a,b$ are  arbitrary, we deduce from Lemma \ref{t1} that
 \begin{align*}
  \noo \summ_{k=1}^n \bar{z}_k\ten z_k\rrm^{\frac12} &\le  \sqrt{2} \pl \noo \summ_{k=1}^n z_k\ten f_k^*\rrm_{M_m(F_n^*)} \qedhere \pl .
  \end{align*}
\qd

\begin{cor}\label{cbisom} $F_n$ is $2$-completely isomorphic to OH$_n$.
\end{cor} \lz

\begin{proof}[\bf Proof:] Since OH$_n$ is selfdual, see \cite{Poh}, it
suffices to apply  Lemma \ref{1s} and Lemma \ref{2s}:
 \begin{align*}
  \noo id:F_n\to OH_n\rrm_{cb} \noo id:OH_n\to F_n\rrm_{cb}
 &=  \noo id:F_n\to OH_n\rrm_{cb} \noo id:F_n^*\to OH_n\rrm_{cb}
 \kl  2  \pl . \qedhere
 \end{align*}\qd

\begin{rem}{\rm A similar result holds in the
context of $L_p$ spaces. We use the standard notation
$H^{c_p}=[H^c,H^r]_{\frac1p}$ and $H^{r_p}=[H^r,H^c]_{\frac1p}$.
Let $2\le p\le \8$ and $p'\le q\le p$. We consider $2\le r \le
\8$ such that  $\frac{1}{r}+\frac{1}{p}=\frac12$. As above, we
consider the subspace $F_n(p)$ of `constant' functions in the
quotient space $G_n=L_2^{c_p}(\nu_1;\ell_2^n)\oplus_q
L_2^{r_p}(\nu_2;\ell_2^n)/\krn(Q)$. The same proofs  as in Lemma
\ref{1s}  and in Lemma \ref{2s} applies with the help of the
formula
 \[ \noo \summ_{k=1}^n x_k\ten e_k\rrm_{S_p^m[OH_n]}
  \lel \sup_{\noo a\rrm_{2r}\le 1, \noo b\rrm_{2r}\le 1}
  \kla \summ_{k=1}^n \noo bx_ka\rrm_2^2 \mer^{\frac12} \pl. \]
Therefore, we have
 \[ \noo id:F_n(p)\to OH_n\rrm_{cb} \kl
 2^{\frac12-\frac1p} \quad \mbox{and} \quad \noo id:F_n(p)^*\to OH_n\rrm_{cb} \kl 2^{\frac12-\frac1p} \pl . \]
These inequalities imply that
 \[ d_{cb}(F_n(p),OH_n) \kl 2^{1-\frac{2}{p}} \pl .\]
We obtain a concrete embedding of OH as a subspace of a quotient
of $S_p=L_p(\B(\ell_2),tr)$ (with $q=p$) and of $S_{p'}$ (with
$q=p'$ and using $H^{c_p}=H^{r_{p'}}$). For an independent,
alternative approach we refer to \cite{X1}.}\end{rem} \lz

\begin{rem} {\rm A slight modification of this approach
yields the space $C_p=[C,R]_{\frac1p}$. Indeed, let
$\al=\frac{1}{p}$. Using the substitution $u=t^{-1}-1$, we find
 \[ \intt_0^1 \frac{1}{(1-t)+tB} \frac{dt}{t^{\al}(1-t)^{1-\al}}  \lel \intt_0^\infty
 \frac{1}{u+B} \pl \frac{du}{u^{1-\al}}\lel
 \frac{c(\al)}{B^{1-\al}} \pl .\]
Following \cite{Wo}, we deduce  for arbitrary commuting operators
$A$, $B$ that

 \begin{samepage} \begin{align*}
 &(x,A^{1-\al}B^{\al}x) \lel \\
 &\quad \inf_{x=f(t)+g(t)} \intt_0^1
 \frac{(f(t),Af(t))}{t} \frac{dt}{c(\al)t^\al(1-t)^{1-\al}}
 + \intt_0^1
 \frac{(g(t),Bg(t))}{1-t} \frac{dt}{c(\al)t^\al(1-t)^{1-\al}} \pl
 .\end{align*}\end{samepage}
Similar as in Lemma \ref{1s}, it then easily follows that  the
space $F_{\al}$ of `constants'  in the quotient
$L_2^c(t^{-1}\mu_{\al})\oplus_1
L_2^r((1-t)^{-1}\mu_{\al})/\krn(Q)$ satisfies $\noo id:F_{\al}\to
[C,R]_{\frac1p}\rrm_{cb}\le \sqrt{2}$. The analogue of Lemma
\ref{2s} is considerably more involved and again based on the dual
Pusz/Woronowicz formula. In the forthcoming publication \cite{JX2}
we will develop better tools for finding the `right pair of
densities'.}
\end{rem}\lz

Assuming Theorem 3, we immediately  obtain an embedding of OH in
a noncommutative $L_1$ space.\lz

\begin{theorem}\label{oh}  OH embeds into the predual of a von Neumann algebra with QWEP.
\end{theorem}\lz

\begin{proof}[\bf Proof:] Note that the isomorphism $id:F_n\to
OH_n$ from Corollary \ref{cbisom} satisfies $id(F_{n-1})\subset
OH_{n-1}$. By the density of $\bigcup F_n$ in
 \[ F\subset L_2^c(\nu_1;\ell_2) \oplus_1   L_2^r(\nu_2;\ell_2)/\krn(Q)\lel G \]
we obtain an isomorphism $id:F\to OH$ with $\noo id\rrm_{cb}\le
\sqrt{2}$. Similarly, the  density  of $\bigcup_n OH_n$ in $OH$
implies that $\noo id^{-1}\rrm_{cb}\le \sqrt{2}$. Hence, OH is
$2$-cb-isomorphic to $F$ and Theorem 3 implies the assertion.\qd


\begin{cor} \label{ohf} There exists a constant $C>0$ such that for all $\nen$, there is an integer $m$ and an injective  linear map $u:OH_n\to S_1^m$ such that
\[ \noo u\rrm_{cb}\noo u^{-1}:u(OH_n)\to OH_n\rrm_{cb} \kl C \pl .\]
\end{cor}

\begin{proof}[\bf Proof:] This is an immediate consequence of the strong
principle of local reflexivity in \cite{EJR} and the fact that $N$
is QWEP.\qd

\section{The projection constant of the operator space OH$_n$}

In this section, we will provide the proof of Theorem 1, assuming
the probabilistic result Theorem 3. The main  tool is a
characterization of the completely summing norms for linear maps
between certain operator spaces, see Lemma 4.4 and Lemma 4.5. Some
notation is required. For an operator space $E$ we denote by
$\iota_E:E\to E^{**}$ the canonical completely isometric embedding
of $E$ into its bidual $E^{**}$. Let $E$ and $F$ be operator
spaces. The $\Gamma_{\infty}$-norm of a  linear map $v:F\to E$ is
defined by
 \[ \Gamma_\8(v) \lel \inf \noo \al \rrm_{cb}  \noo \beta
 \rrm_{cb}  \pl .\]
Here the infimum is taken over all $\al:F\to \B(H) $, $\beta:\B(H)
\to E^{**}$ such that $\iota_Ev=\al\beta$. For a finite rank map
$v:F\to E$, we also define
 \[ \gamma_\8(v) \lel \inf \noo \al \rrm_{cb}  \noo \beta
 \rrm_{cb}  \pl .\]
Here the infimum is taken over all $m\in \nz$, $\beta:F\to  M_m$
and $\al:M_m\to E$ such that $v=\al\beta$. Following Effros/Ruan
(see e.g.\cite{ER})  a linear map $v:E\to F$ is said to be
\emph{{\rm (}completely{\rm )} 1-summing} if
 \[ \pi_1^o(v)\lel \noo id_{S_1}\ten v:S_1\ten_{\min}E\to S_1\wet F\rrm
 \]
is finite. The $\gamma_\8$-norm is related to the $1$-summing norm
via trace duality.  \lz

\begin{lemma}\label{trace} Let $u:E\to F$ be a linear map.  Then
 \[ \pi_1^o(u) \lel \sup\{ | tr(vu)| \pl |\pl \gamma_\8(v)\le 1\}
 \pl .\]
\end{lemma}

\begin{proof}[\bf Proof:] Indeed, we have
 \begin{align*}
  \pi_1^o(u) &= \sup\left\{  |\langle (id\ten u)(x),y\rangle| \pl|\pl
  \noo x\rrm_{S_1^m\ten_{\min} E}\le 1 \pl ,\pl \noo
  y\rrm_{M_m(F^*)}\le 1\right\}\\
  &= \sup \left\{ |tr(T_y^*uT_x)| \pl|\pl \noo T_x:M_m\to E\rrm_{cb}\le
  1\pl,\pl
  \noo T_y:M_m^*\to
  F^*\rrm_{cb} \le 1 \right\} \\
  &= \sup \left \{ |tr(vu)| \pl|\pl \gamma_\8(v)\le 1 \right\} \pl
  . \qedhere
 \end{align*}\qd
We will need the following result from \cite{EJR}.\lz

\begin{lemma}[EJR]\label{webfc} Let $E$ and $F$ be finite dimensional operator
spaces and $v:F\to E$. Then
 \[ \gamma_\8(v)\lel \Gamma_\8(v) \pl .\]
\end{lemma}\lz

The connection between the 1-summing norm and the operator space
projective tensor norm works only  for subspaces of
noncommutative $L_1$ spaces. Therefore, we follow \cite{Pip} and
define
 \for
 d_{SL_1}(E) &=&
 \inf_{w:E\to E_1\subset S_1^m}  \noo w\rrm_{cb}\noo
 w^{-1}\rrm_{cb} \pl.
 \mel
We are interested in estimates for finite dimensional spaces.
This leads to the following definition for infinite dimensional
operator spaces
 \for
 d_{SL_1}(Y) &=& \sup_{E\subset Y}  d_{SL_1}(E) \pl
 .
 \mel
Here the supremum is taken over all finite dimensional subspaces
$E\subset Y$. The following fact follows immediately from the
definition. \lz

\begin{lemma}\label{ptfd}
Let $x\in X\ten Y$ and $\eps>0$.  Then  there  are finite
dimensional subspaces $E\subset X$ and $F\subset Y$ such that
 \[ \noo x\rrm_{E\wet F} \kl (1+\eps) \noo
 x\rrm_{X\wet Y} \pl .\]
\end{lemma}\lz

\begin{lemma} \label{linf1}
Let $X$ and $Y$ be operator spaces, $E\subset X$ and $F\subset Y$
be finite dimensional subspaces. Let  $x\in E\ten F$ be a tensor
with associated linear map $T_x:E^*\to Y$, $T_x(e^*)=(e^*\ten
id)(x)$. Then
 \[ \pi_1^{o}(T_x)\kl d_{SL_1}(Y)\noo x\rrm_{X\wet
 Y} \pl. \]
\end{lemma}\lz

\begin{proof}[\bf Proof:] According to Lemma \ref{ptfd}, we can find $\tilde{F}\supset F$ such that
 \[ \noo x\rrm_{X\wet \tilde{F}}\le (1+\eps)\noo x\rrm_{X\wet Y} \pl .\]
Let  $w:\tilde{F}\to F_1\subset S_1^n$ be a linear isomorphism
with completely contractive inverse $w^{-1}:F_1\to \tilde{F}$. In
order to estimate the 1-summing norm, we have to consider complete
contractions $\beta:Y\to M_m$ and $\al:M_m\to E^*$. By Wittstock's
extension theorem, there is a complete contraction
$\hat{\beta}:S_1^n\to M_m$ such that $\hat{\beta}|_{F_1}=\beta
w^{-1}$. Then $\al \hat{\beta}$ corresponds to an element $z\in
E^*\ten_{\min} M_n$ of norm less than one. The injectivity of the
projective tensor product on $S_1^n$ (see \eqref{inject}) yields
 \for
  |tr(T_x\al\beta)| &=& |tr(\al \beta T_x)|
  \lel |tr(\al\hat{\beta}wT_x)| \lel |\langle z, (id_E\ten w)(x) \rangle|\\
  &\le& \noo z\rrm_{E^*\ten_{\min} M_n} \noo (id_E\ten w)(x)  \rrm_{E\wet S_1^n }
  \lel \noo z\rrm_{E^*\ten_{\min} M_n} \noo (id_X\ten w)(x)
  \rrm_{X\wet S_1^n} \\
  &\le& \noo (id_X\ten w)(x)
  \rrm_{X\wet F_1} \kl \noo w\rrm_{cb} \noo x\rrm_{X\wet \tilde{F}}
  \kl (1+\eps) \noo w\rrm_{cb} \noo x\rrm_{X\wet Y}  \pl .
 \mel
Taking the infimum over all $w$, we may replace $\noo w\rrm_{cb}$
by $d_{SL_1}(\tilde{F})$ and then by $d_{SL_1}(Y)$.  \qd

We have  a  partial converse to this inequality. \lz

\begin{lemma} \label{linf2}
Let $M$ and $N$ be von Neumann algebras such that $N$ is QWEP. Let
$Y$ be an  operator space and $w:Y\to M_*$ be a complete
contraction. Let $X$ be an operator space, $E\subset X$ be a
finite dimensional subspace and $u:E\to N_*$ be a complete
contraction. If $x\in E\ten Y$  and $T_x:E^*\to Y$ is the
associated map, then
 \[ \pi_1^o(T_x:E^*\to Y) \gl \noo (u\ten
 w)(x)\rrm_{N_*\wet M_*} \pl .\]
\end{lemma}{\lz}

\begin{proof}[\bf Proof:] Let us assume that $N=A/I$ where $A$ has
WEP and let $\pi:A\to N$ be the quotient homomorphism. Then the
map  $\pi^*:N^*\to A^*$ has a  completely contractive left
inverse. Indeed, let $z$ be a central projection in $A^{**}$ such
that $I^{**}=zA^{**}$. Then we have a  canonical isomorphism
$(A/I)^*\cong N^*$. Therefore the mapping $v:A^*\to (A/I)^*\cong
N^*$ given by $v(x)=zx$ satisfies $id_{N^*}\lel  v\pi^*$. Since
$N_*$ is completely complemented in $N^*$, we also have a
completely contractive left inverse for the restriction
$\pi^*|_{N_*}$. This implies that
 \[ \pi^*\ten id_{M_*}: N_*\wet M_*\to  A^*\wet M_* \]
is completely isometric. We consider the tensor $(\pi^*u\ten
w)(x)\in A^*\wet M_*$.  By the Hahn-Banach theorem we may find
 \[ y\in (A^*\wet M_*)^* \lel CB(M_*,A^{**}) \cong
 M^{**}\bar{\ten}A^{**} \]
of norm less than one such that
 \[ |\langle y, (\pi^*u\ten w)(x)\rangle| \lel \noo  (\pi^*u\ten
 w)(x)\rrm_{A^*\wet M_*} \pl .\]
Let $\eps>0$. According to Kaplansky's density theorem (see
\cite{EJR,ER} for details), we may  find  a finite rank element
$y_{\eps} \in M\ten A$ of norm $\le 1$ such that
 \[ \noo  (\pi^*u\ten  w)(x)\rrm_{N_*\wet M_*}\kl (1+\eps) |\langle
 y_{\eps},(\pi^*u\ten w)(x)\rangle| \pl .\]
The  tensor $y_{\eps}$ corresponds to a complete contraction
$S_{y_{\eps}}:M_*\to A$ (see \cite[Theorem 7.3.2]{ER}). We need
the finite dimensional subspace  $F=S_{y_\eps}wT_x(E^*)$. We
consider $\tilde{T}=S_{y_{\eps}}wT_x:E^* \to F$ as a map with
values in its range $F$. Using rank one tensors it is easy to
check that
 \[ \langle  y_{\eps},(\pi^* u\ten w)(x)\rangle \lel tr(u^*\pi S_{y_\eps}wT_x)  \lel tr(u^*\pi|_F\tilde{T})  \pl .\]
We apply Lemma \ref{trace}, Lemma \ref{webfc} and basic
properties of the 1-summing norm:
 \begin{align*}
 |tr(u^*\pi S_{y_{\eps}}wT_x)| &= |tr(u^*\pi|_F\tilde{T})|
  \le \gamma_{\infty}(u^*\pi|_F) \pi_1^o(\tilde{T}) \kl
 \Gamma_{\infty}(u^*\pi|_F) \pi_1^o(S_{y_{\eps}}wT_x)\\
  &\le
  \Gamma_{\infty}(u^*\pi|_F)  \noo S_{y_\eps}\rrm_{cb} \noo w\rrm_{cb}
  \pi_1^o(T_x) \kl
    \Gamma_{\infty}(u^*\pi |_F) \pi_1^o(T_x) \pl .
  \end{align*}
Since $A$ is WEP, there is a complete contraction $P:\B(H) \to
A^{**}$ such that $P|_A\lel id_A$ where $A\subset A^{**}\subset
\B(H)$. This implies that
 \[ \Gamma_{\infty}(u^*\pi|_F)\lel
 \Gamma_{\infty}(u^{***} \pi|_F) \lel
  \Gamma_{\infty}(u^{***}P\pi|_F) \kl
 \noo u^{***}P\rrm_{cb} \lel \noo u\rrm_{cb} \pl .\]
Thus, we get $\noo  (u\ten  w)(x)\rrm_{N_*\wet M_*}\kl (1+\eps)
\pi_1^o(T_x)$. Letting $\eps\to 0$, the assertion follows.\qd

In the following we use the notation  $G$ for the spaces
introduced in section 3. The following proposition will be proved
in section 7. \lz

\begin{prop}\label{kick00}
There exists a von Neumann algebra $N$ with QWEP and a completely
contractive injective map  $u:G\to N_*$ such that
 \[ \noo (u\ten u)(x) \rrm_{N_*\wet  N_*} \gl \frac{1}{9} \noo
 x\rrm_{G\wet G} \pl \]
for all $x\in G\ten G$. Moreover, $\noo u^{-1}:u(G)\to
G\rrm_{cb}\le 3$ and $u(G)$ is completely complemented in $N_*$.
\end{prop}\lz

\begin{cor}\label{pp} Let $F_n\subset G$ be
a finite dimensional subspace and $x\in F_n\ten G$. Then
 \[ \frac{1}{9}\noo x\rrm_{G\wet G}\kl \pi_1^o(T_x:F_n^*\to G) \kl 3 \noo x\rrm_{G\wet G}  \]
\end{cor}\lz

\begin{proof}[\bf Proof:] Since $G$ is $3$-cb isomorphic to a subspace of
$N_*$ and $N$ is QWEP, the strong principle of local reflexivity
in \cite{EJR} implies that $d_{SL_1}(G)\le 3$. Thus by Lemma
\ref{linf1}
 \[  \pi_1^o(T_x:F_n^*\to G) \kl 3 \noo x\rrm_{G\wet G}
 \pl .
 \]
Now, we prove the converse inequality. By Proposition
\ref{kick00}, we may apply  Lemma \ref{linf2} to get
 \begin{align*}
  \noo x\rrm_{G\wet G} &\le  9 \noo (u\ten u)(x)\rrm_{N_*\wet N_*}
 \kl 9  \pl  \pi_1^o(T_x:F_n^*\to G)
 \pl . \qedhere
 \end{align*} \qd

\begin{cor}\label{th4}
Let $(f_i)$ be the canonical unit vector basis in the space $F_n$
constructed in section 3. Let $u:OH_n\to OH_n$ be a linear map
represented by a matrix $[a_{ij}]$. Then
 \[  \frac{1}{6} \pl \pi_1^{o}(u) \kl
 \noo \summ_{i,j=1}^n a_{ij} f_i\ten
  f_j \rrm_{G_n\wet G_n} \kl 18 \pl \pi_1^{o}(u) \pl .\]
\end{cor}

\begin{proof}[\bf Proof:]  Let $u:OH_n\to OH_n$ be a linear map  represented by
a  matrix $(a_{ij})$. We deduce from Lemma \ref{1s}, Lemma
\ref{2s} and Corollary \ref{pp} that
 \begin{align*}
  \pi_1^o(u:OH_n\to OH_n) &\le 2 \pl
  \pi_1^o(u:F_n^*\to
  F_n) \kl 6 \noo \summ_{i,j=1}^n a_{ij} f_i\ten
  f_j \rrm_{G\wet G} \pl .
  \end{align*}
Conversely, we have
 \begin{align*}
 \noo \summ_{i,j=1}^n a_{ij} f_i\ten
  f_j \rrm_{G\wet G}
  &\le  9 \pl \pi_1^o(u:F_n^*\to F_n) \kl 18 \pl
  \pi_1^o(u:OH_n\to OH_n) \pl . \qedhere
  \end{align*}
\qd

The following norm calculations in $G_n\wet G_n$ will be postponed
to the next section.\lz

\begin{prop}\label{ncalc} Let $[a_{ij}]$ be an $n\times n$ matrix. Then
 \[ \noo \summ_{i,j=1}^n a_{ij} f_i\ten f_j \rrm_{G_n\wet G_n} \kl
 18
 \pl
 \sqrt{1+\ln n} \kla \summ_{i,j=1}^n |a_{ij}|^2 \mer^{\frac12} \pl
 .\]
Moreover, for $n\ge 7$
 \[ \noo \summ_{i=1}^n  f_i\ten f_i \rrm_{G_n\wet G_n} \gl
 (16\sqrt{2}\pi)^{-1} \pl \sqrt{n(1+\ln n)} \pl .\]
\end{prop}\lz

As an application, we derive an independent proof of \eqref{uppb}.
\lz

\begin{samepage}
\begin{cor}\label{psc}  Let $u:\B(H) \to OH$ be a completely bounded map. Then
 \begin{equation}\label{lg1}  \kla \summ_{k=1}^n \noo u(x_k)\rrm_{OH}^2 \mer^{\frac12} \kl 108
 \pl
 \sqrt{1+\ln n} \pl \noo  u\rrm_{cb} \pl \noo \summ_{k=1}^n \bar{x}_k
 \ten x_k\rrm_{\overline{\B(H) }\ten_{\min} \B(H) }^{\frac12} \pl
 \end{equation}
for all $\nen$ and $x_1,..,x_n\in \B(H) $. Moreover, let
$B:\overline{\B(H) }\times \B(H) \to \cz$ be a positive
sesquilinear form of rank $n$. Then there exists a positive
integral sesquilinear form $\tilde{B}$ such that $B\le \tilde{B}$
and
 \[ \|  \tilde{B}\|_I\kl 2\times 108^2 (1+\ln n) \noo B\rrm_{jcb} \pl .\]
\end{cor} \end{samepage}

 \begin{proof}[\bf Proof:] Let $u:\B(H)\to OH$ be a
completely bounded map and $x_1,....,x_n\in \B(H)$. We can find
an orthogonal projection $P:OH\to {\rm span}\{u(x_k)|1\le k\le
n\}$ of rank at most $n$. Since OH is homogeneous, we may assume
that $P(OH)=OH_n$. A glance at  \eqref{ohnorm} shows that
$w:OH_n\to \B(H)$ given by $w(e_k)=x_k$ satisfies
 \[ \noo w\rrm_{cb} \lel \noo \summ_{k=1}^n \bar{x}_k\ten
 x_k\rrm_{\overline{\B(H) }\ten_{\min} \B(H) }^{\frac12} \pl .\]
This implies that the composition map $v=Puw$ satisfies
 \for
 \gamma_\8(v)&=&  \Gamma_\8(uw)\kl \noo u\rrm_{cb} \noo
 w\rrm_{cb} \lel
 \noo u\rrm_{cb} \pl \noo \summ_{k=1}^n \bar{x}_k
 \ten x_k\rrm_{\overline{\B(H) }\ten_{\min} \B(H) }^{\frac12} \pl .
 \mel
Now, we apply trace duality. Let  $(a_{ij})\in \ell_2^n(OH_n)\lel
\ell_2^{n^2}$ be of norm one such that
 \[ \kla \summ_{k=1}^n \noo Pu(x_k)\rrm_{OH_n}^2 \mer^{\frac12}
 \lel \summ_{k,l} a_{kl}(e_l,v(e_k)) \lel tr(av) \pl .
 \]
We deduce from Lemma \ref{trace}, Corollary \ref{th4} and
Proposition \ref{ncalc}  that
 \begin{align*}
 \kla \summ_{k=1}^n \noo u(x_k)\rrm_{OH}^2 \mer^{\frac12}  &=
 |tr(av)| \kl \pi_1^o(a)  \gamma_\8(v)\\
 & \le 6\times  18 \pl \sqrt{1+\ln n} \pl \noo a\rrm_2 \pl
 \noo u\rrm_{cb} \pl \noo \summ_{k=1}^n \bar{x}_k
 \ten x_k\rrm_{\overline{\B(H) }\ten_{\min} \B(H) }^{\frac12} \\
 & \le 108 \pl  \sqrt{1+\ln n} \pl \noo u\rrm_{cb} \pl \noo \summ_{k=1}^n \bar{x}_k
 \ten x_k\rrm_{\overline{\B(H) }\ten_{\min} \B(H) }^{\frac12} \pl
 .
 \end{align*}
This completes the proof of the first assertion. Following
\cite[(9.3)]{Poh} (and ultimately \cite{TJ}) this implies that
 \[
 \kla \summ_{k}  \noo u(x_k)\rrm_{OH}^2 \mer^{\frac12}
 \kl \sqrt{2} 108 \pl  \sqrt{1+\ln n} \pl \noo u\rrm_{cb} \pl \noo \summ_{k} \bar{x}_k
 \ten x_k\rrm_{\overline{\B(H) }\ten_{\min} \B(H) }^{\frac12} \pl
 \]
for all $u$ of rank at most $n$ and arbitrary sequence $(x_k)$.
An appeal to  Lemma \ref{transl} concludes the proof.\qd

Now, we apply a typical trace duality argument.\lz

\begin{cor} Let $\nen$. Then
 \[ \gamma_\8(id_{OH_n}) \kl 288 \sqrt{2}\pi  \sqrt{\frac{n}{1+\ln n}} \pl .\]
\end{cor}\lz

\begin{proof}[\bf Proof:] Using a well-known averaging trick, we
have
 \begin{equation}\label{trdua} \gamma_{\infty}(id_{OH_n})\pi_1^o(id_{OH_n})\lel
 n \pl.
 \end{equation}
(Indeed, according to Lemma \ref{trace} there exists $v$ with
$\gamma_{\infty}(v)=1$ and $\pi_1^{o}(id)\lel |tr(v)|$. Let $\si$
be the normalized Haar measure on the unitary group $U_n$. Then
$\tilde{v}=\int_{U_n} uvu^{-1}d\si(u)$ satisfies
$\gamma_{\infty}(\tilde{v})\le 1$ and $n\tilde{v}=tr(v)id$. This
implies that $\gamma_{\infty}(id)\pi_1^o(id)\lel
\gamma_{\infty}(id) |tr(v)|\le n$. The converse inequality is
obvious.) For $n\ge 7$ we deduce from Corollary \ref{th4} and
Proposition \ref{ncalc} that
 \begin{equation}\label{alln}  \pi_1^o(id_{OH_n})\gl \frac{1}{18}\noo \summ_{i=1}^n
 f_i\ten
 f_i\rrm_{G_n\wet G_n} \gl \frac{1}{288\sqrt{2}\pi} \sqrt{n(1+\ln n)} \pl .
  \end{equation}
For $n\le 7$ we use the well-known Banach space estimate
$\pi_1^o(id_{OH_n})\gl \frac{2}{\sqrt{\pi}}\sqrt{n}$. Hence
\eqref{alln} is valid for all $\nen$ and the assertion follows
from \eqref{trdua}.\qd

We are ready for the proof of Corollary 5.\lz

\begin{cor}\label{proj} Let $OH_n\subset \B(\ell_2)$, then there exists
a projection $P:\B(\ell_2)\to OH_n$ such that
 \[ \noo P\rrm_{cb} \kl 288\sqrt{2} \pi \pl  \sqrt{\frac{n}{1+\ln n}} \pl .\]
In particular,
 \[  \frac{1}{108} \sqrt{\frac{n}{1+\ln n}}\kl  \la_{cb}(OH_n)\kl 288\sqrt{2} \pi \pl  \sqrt{\frac{n}{1+\ln n}} \pl .\]
\end{cor}\lz

\begin{proof}[\bf Proof:] We write
$id_{OH_n}=vw$, with $w:OH_n\to \B(H) $ and  $v:\B(H) \to OH_n$,
and
 \[ \noo v\rrm_{cb} \noo w\rrm_{cb} \kl 288 \sqrt{2} \pi \pl  \sqrt{\frac{n}{1+\ln n}}
 \pl .\]
According to Wittstock's extension theorem, we can find an
extension $\hat{w}:\B(\ell_2)\to \B(H) $ with the same cb-norm as
$w$. Then $P=v\hat{w}$ is the corresponding projection. The lower
estimate follows easily from Corollary \ref{psc} (see also
\cite{PS}). \qd

\begin{cor} The order $(1+\ln n)$  in
\eqref{uppb} is  best possible. Moreover, there is a sesquilinear
jcb form which can not be majorized by an integral form.
\end{cor} \lz

\begin{proof}[\bf Proof:] Let $\iota:OH_n\to \B(\ell_2)$ be a completely
isometric embedding and $x_k=\iota(e_k)$. According to Corollary
\ref{proj}, we can find a projection $P:\B(\ell_2)\to OH_n$ of
cb-norm
 \[ \noo P\rrm_{cb} \kl 288 \sqrt{2} \pi  \sqrt{\frac{n}{1+\ln n}} \pl .\]
We define $B_n(x,y)\lel (P(x),P(y))$ and assume that $B_n\le
\tilde{B}$. Then, we deduce from Lemma \ref{transl} that
 \begin{align*}
 n &= \summ_{k=1}^n \noo e_k\rrm^2 \lel  \summ_{k=1}^n B_n(x_k,x_k)
 \kl \| \tilde{B}\|_I \noo \summ_{k=1}^n \bar{x}_k\ten
 x_k\rrm_{\overline{\B(H) }\ten_{\min}\B(H) } \lel
 \| \tilde{B}\|_I  \pl .
 \end{align*}
This implies that
 \[ \noo B_n\rrm_{jcb}
   \lel \noo P\rrm_{cb}^2 \kl (288 \sqrt{2}\pi)^2 \frac{n}{1+\ln n} \quad
 \mbox{and} \quad
  n\kl \| \tilde{B}\|_I \pl .\]
Combining these estimates, we deduce that
 \begin{align*} (1+\ln n) \noo B_n\rrm_{jcb}&\le
 (288 \sqrt{2}\pi)^2
  \| \tilde{B}\|_I \pl .
  \end{align*}
For the second assertion, we define $B\lel \sum_{k\in \nz} k^{-2}
\frac{k^4}{2^{k^4}}B_{2^{k^4}}$. The triangle inequality shows
that $B$ is jcb. However every if $B\le \tilde{B}$, then
$k^{-2}\frac{k^4}{2^{k^4}}B_{2^{k^4}}\le  \tilde{B}$ implies that
$\| \tilde{B}\|_I\ge c k^{-2}k^4$ for all $k\in \nz$.
 \qd

\begin{rem} {\rm This argument also shows that $\sqrt{1+\ln n}$ in Lemma
\ref{triv} is best possible. Moreover, we see that $\al=1$ in
Theorem 5  is best possible for non sub-homogeneous
$C^*$-algebras. Indeed, since $\gamma_{\infty}(id_{OH_n})\le 288
\pi \sqrt{2n/(1+\ln n)}$, we may find a complete contraction
$w:OH_n\to M_m$ and a linear map $u:M_m\to OH_n$ such that
$uw=id_{OH_n}$ and $\noo u\rrm_{cb} \le 288 \pi \sqrt{2n/(1+\ln
n)}$. Moreover, it is well-known (see e.g. \cite{JNRX}) that
$M_m$ is $(1+\eps)$ completely complemented in a non
subhomogeneous $C^*$-algebra $A$, i.e. there is a complete
contraction $\al:M_m\to A$ and a map $\beta:A\to M_m$ of cb-norm
$\le (1+\eps)$ such that $\beta \al=id$. Then
$B(x,y)=(u\beta(x),\beta u(y))$ provides the `counterexample' on
$A$.}
\end{rem}

\section{Norm calculations in a quotient space}

Although the calculations in this section are of technical nature,
the idea is very simple: The influence of the singularities at the
corners $(0,1)$ and $(1,0)$ is minimized by rectangular
decompositions. The next lemma justifies the use of these
decompositions. \lz

\begin{lemma}\label{start}
$G_n\stackrel{\wedge}{\ten} G_n$ is  isometrically  isomorphic to
the quotient space  of
\[ L_2(\nu_1\ten \nu_1;\ell_2^{n^2}) \oplus_1
  L_2(\nu_1;\ell_2^n)\ten_{\pi}  L_2(\nu_2;\ell_2^n)
  \oplus_1
  L_2(\nu_2;\ell_2^n)\ten_{\pi}  L_2(\nu_1;\ell_2^n)
  \oplus_1
  L_2(\nu_2\ten \nu_2;\ell_2^{n^2}) \]
with respect to
 \[ S \lel \{(f,g,h,k) \pl |\pl f(t,s)
 +g(t,s)+h(t,s)+k(t,s)\lel 0 \pll \mu\ten \mu \mbox{ -
 a.e.}\} \pl .\]
\end{lemma} \lz

\begin{proof}[\bf Proof:] By the properties of the projective operator
space tensor product, we have
 \for
 G_n\wet G_n &=& (L_2^c(\nu_1;\ell_2^n)\oplus_1
 L_2^r(\nu_2;\ell_2^n))\wet (L_2^c(\nu_1;\ell_2^n)\oplus_1
 L_2^r(\nu_2;\ell_2^n))/\krn(Q\ten Q) \pl .
 \mel
We note that $H^c \stackrel{\wedge}{\ten}K^c= H\ten_2 K=
 H^r \stackrel{\wedge}{\ten}K^r$ and $H^c
 \stackrel{\wedge}{\ten}K^r=H\ten_{\pi} K\lel H^r
 \stackrel{\wedge}{\ten}K^c$.
Therefore the properties  of $\oplus_1$ and $\wet$ imply that
\begin{align*}
 & \big(L_2^c(\nu_1;\ell_2^n)\oplus_1
 L_2^r(\nu_2;\ell_2^n)\big)\wet \big(L_2^c(\nu_1;\ell_2^n)\oplus_1
 L_2^r(\nu_2;\ell_2^n)\big)\\
  &\pl  \lel  \big(L_2^c(\nu_1;\ell_2^n)\wet L_2^c(\nu_1;\ell_2^n)\big) \oplus_1
       \big(L_2^c(\nu_1;\ell_2^n)\wet L_2^r(\nu_2;\ell_2^n)\big) \\
  &  \pll \pll \pll     \oplus_1 \big(L_2^r(\nu_2;\ell_2^n)\wet L_2^c(\nu_1;\ell_2^n)\big) \oplus_1
    \big(L_2^r(\nu_2;\ell_2^n)\wet L_2^r(\nu_2;\ell_2^n)\big)\\
  &\pl  \lel L_2(\nu_1\ten \nu_1;\ell_2^{n^2})\oplus_1
  \big (L_2(\nu_1;\ell_2^n)\ten_{\pi} L_2 (\nu_2;\ell_2^n)\big) \\
  &  \pll \pll \pll \oplus_1  \big(L_2(\nu_2;\ell_2^n) \ten_{\pi}  L_2(\nu_1;\ell_2^n)\big) \oplus_1
  L_2^r(\nu_2\ten \nu_2;\ell_2^{n^2}) \pl .
 \end{align*}
Using $(1.1)$, we observe  that all four components can be
represented by $\mu\ten \mu$ measurable functions. Applying $Q\ten
Q$ yields the assertion. \qd

\begin{cor}\label{ffunct} Let  $a=[a_{ij}]$ be an $n\times n$-matrix. Then
 \for
  \noo \summ_{i,j=1}^n  a_{ij}f_i\ten f_j \rrm_{G_n\stackrel{\wedge}{\ten} G_n}
  &\ge& \frac{|tr(a)|}{\sqrt{n}}  \pl \sup \bet \intt f(t,s) \frac{d\mu(t)}{t}
  \frac{d\mu(s)}{s} \rag \pl ,
  \mel
where the supremum  is taken over all measurable functions
$(f,g,h,k)$ such that
 \[ \frac{f(t,s)}{ts} \lel \frac{g(t,s)}{(1-t)(1-s)}
  \lel \frac{h(t,s)}{t(1-s)} \lel
  \frac{k(t,s)}{(1-t)s} \pll  \mu\ten \mu \mbox{  a.e. }\]
and
 \begin{align}
 \max\{  \| f\|_{L_2(\nu_1\ten \nu_1)},
 \| g\|_{L_2(\nu_2\ten \nu_2)}\} &\le 1 \pl ,\label{sqn1}\\
  \max\{ \noo h\rrm_{L_2(\nu_1)\ten_{\eps}L_2(\nu_2)},
  \| k\|_{L_2(\nu_2)\ten_{\eps}L_2(\nu_1)}\} &\le \sqrt{n} \pl .
  \label{sqn}
 \end{align}
\end{cor}\lz

\begin{proof}[\bf Proof:] Let $(f,g,h,k)$ be given as above. Consider a
decomposition of $a$ in matrix valued functions
$a=a^1(t,s)+a^2(t,s)+a^3(t,s)+a^4(t,s)$ such that
 \for
 \lefteqn{ \noo a^1\rrm_{L_2(\nu_1\ten\nu_1;\ell_2^{n^2})}
  +  \noo
  a^2\rrm_{L_2(\nu_2\ten\nu_2;\ell_2^{n^2})}}\\
  & &  + \noo a^3\rrm_{L_2(\nu_1;\ell_2^n)\ten_{\pi}
 L_2(\nu_2;\ell_2^n)} +
 \noo a^4\rrm_{L_2(\nu_2;\ell_2^n)\ten_{\pi} L_2(\nu_1;\ell_2^n)} \kl
 (1+\eps) \pl \noo \summ_{i,j=1}^n a_{ij} f_i\ten f_j \rrm_{G_n\stackrel{\wedge}{\ten}
 G_n}\pl .
 \mel
Then the Cauchy Schwarz inequality implies that
 \for
 \lefteqn{ \bet  \intt \summ_{i=1}^{n} a^1_{ii}(t,s) f(t,s)
 \frac{ d\mu(t)}{t} \frac{d\mu(s)}{s} \rag }\\
 & & \kl \kla \intt \summ_{i=1}^{n} |a^1_{ii}(t,s)|^2
 \frac{ d\mu(t)}{t} \frac{d\mu(s)}{s} \mer^{\frac12} \pl
 \kla  \intt \summ_{i=1}^{n} |f(t,s)|^2
 \frac{ d\mu(t)}{t} \frac{d\mu(s)}{s}
  \mer^{\frac12}\\
 & & \lel  \noo a^1\rrm_{L_2(\nu_1\ten
    \nu_1;\ell_2^{n^2})} \pl\sqrt{n} \pl  \noo f\rrm_{L_2(\nu_1\ten
    \nu_1)} \kl \sqrt{n} \pl \noo
    a^1\rrm_{L_2(\nu_1\ten
    \nu_1;\ell_2^{n^2})}\pl .
    \mel
Similarly,
 \[
  \bet  \intt \summ_{i=1}^{n} a^2_{ii}(t,s) g(t,s)
    \frac{ d\mu(t)}{1-t} \frac{d\mu(s)}{1-s} \rag
    \kl \sqrt{n}\noo a^2\rrm_{L_2(\nu_2\ten
    \nu_2;\ell_2^{n^2})}   \pl .\]
For every operator $h:L_2\to L_2$, we recall that $\noo h\ten
id_{\ell_2^n}\rrm=\noo h\rrm$. Hence, we deduce from trace duality
and \eqref{sqn} that
 \begin{align*}
   \bet  \intt \summ_{i=1}^{n} a^3_{ii}(t,s) h(t,s)
     \frac{ d\mu(t)}{t} \frac{d\mu(s)}{1-s} \rag
    &\le   \noo h\ten id_{\ell_2^n}\rrm_{
    L_2(\nu_1;\ell_2^n)\ten_{\eps}
    L_2(\nu_2;\ell_2^n)} \pl
      \noo
    a^3\rrm_{L_2(\nu_1;\ell_2^n)\ten_{\pi}
    L_2(\nu_2;\ell_2^n)} \\
    &\le   \noo
    h\rrm_{L_2(\nu_1)\ten_{\eps}L_2(\nu_2)}  \pl \noo a^3\rrm_{L_2(\nu_1;\ell_2^n)\ten_{\pi}
    L_2(\nu_2;\ell_2^n)} \\
    &\le
     \sqrt{n} \noo a^3\rrm_{L_2(\nu_1;\ell_2^n)\ten_{\pi}
    L_2(\nu_2;\ell_2^n)} \pl .
    \end{align*}
Similarly,
 \begin{align*}    \bet  \intt \summ_{i=1}^{n} a^4_{ii}(t,s) k(t,s)
    \frac{ d\mu(t)}{1-t} \frac{d\mu(s)}{s} \rag &\le  \noo k\rrm_{L_2(\nu_2)\ten_{\eps}L_2(\nu_1)} \pl
    \noo a^4\rrm_{L_2(\nu_2;\ell_2^n)\ten_{\pi}
    L_2(\nu_1;\ell_2^n)}  \\
    &\le  \sqrt{n}\pl
    \noo a^4\rrm_{L_2(\nu_2;\ell_2^n)\ten_{\pi}
    L_2(\nu_1;\ell_2^n)}  \pl .
    \end{align*}
Therefore, we get
 \begin{align*}
   \bet \intt  \summ_{i=1}^n a_{ii} f(t,s)
 \frac{d\mu(t)}{t}\frac{d\mu(s)}{s}\rag
 &\le \bet \intt \summ_{i=1}^n a^1_{ii}f(t,s)
 \frac{d\mu(t)}{t}\frac{d\mu(s)}{s} \rag +
 \bet \intt \summ_{i=1}^n a^2_{ii}g(t,s)
 \frac{d\mu(t)}{1-t}\frac{d\mu(s)}{1-s} \rag \\
 & \pll +\bet \intt \summ_{i=1}^n a^3_{ii}h(t,s)
 \frac{d\mu(t)}{t}\frac{d\mu(s)}{1-s} \rag +
 \bet \intt \summ_{i=1}^n a^4_{ii}k(t,s)
 \frac{d\mu(t)}{1-t}\frac{d\mu(s)}{s} \rag \\
 &  \le  \sqrt{n}
 (1+\eps)
 \noo \summ_{i,j=1}^n a_{ij} f_i\ten f_j \rrm_{G_n\stackrel{\wedge}{\ten} G_n} \pl
 . \qedhere
 \end{align*}
\qd

We will now prove the lower estimate in Proposition \ref{ncalc}.

\begin{lemma} \label{lest} Let $n\ge 7$. Then
 \[ \noo \summ_{i=1}^n f_i\ten f_i\rrm_{G_n\stackrel{\wedge}{\ten} G_n} \gl
 \frac{1}{16\sqrt{2}\pi}
 \sqrt{n(1+\ln n)} \pl. \]
\end{lemma}

\begin{proof}[\bf Proof:]
Let $0<\delta<\frac12$, to be determined later. We consider the
rectangle $I=[\delta,\frac12]\times [\frac12,1-\delta]$ and the
function
 \[ v(t,s)\lel  \frac{1}{ts+(1-t)(1-s)}  \pl 1_I \pl .\]
Following Corollary \ref{ffunct} we define  $f(t,s)\lel ts v(t,s)$
and $g(t,s)=(1-t)(1-s)v(t,s)$ and observe that
 \begin{align*}
  \intt_I f(t,s)^2 \frac{d\mu(t)}{t}\frac{d\mu(s)}{s}
 &+ \intt_I g(t,s)^2
 \frac{d\mu(t)}{1-t}\frac{d\mu(s)}{1-s}
  \lel    \intt_I v(t,s)^2[ts+(1-t)(1-s)] d\mu(t)d\mu(s)\\
 &=   \intt_I \frac{1}{ts+(1-t)(1-s)}  d\mu(t)d\mu(s) \\
 & \le   4\pi^{-2} \intt_{\delta}^{\frac12}
 \intt_{\frac12}^{1-\delta} \min(t^{-1},(1-s)^{-1})
 \frac{dt}{\sqrt{t}} \frac{ds}{\sqrt{1-s}} \\
 & =   4\pi^{-2} \intt_{\delta}^{\frac12}
 \intt_{\delta}^{\frac12} \min(t^{-1},s^{-1})
 \frac{dt}{\sqrt{t}} \frac{ds}{\sqrt{s}}
  \lel  8\pi^{-2} \intt_{\delta}^{\frac12}
 \intt_{\delta}^t \frac{ds}{\sqrt{s}}  \pl
 \frac{dt}{t\sqrt{t}} \\
 & \le  16 \pi^{-2} \intt_{\delta}^{\frac12} \sqrt{t}
 \frac{dt}{t\sqrt{t}} \kl 16 \pi^{-2} (-\ln \delta)
 \pl .
 \end{align*}
In order to  estimate the norm for $h(t,s)=t(1-s)v(t,s)$, we use
$(1.1)$. Hence it suffices to estimate the $L_2$-norm:
 \begin{align*}
  \noo h\rrm_{L_2(\nu_1\ten \nu_2)}^2 &\le  4 \intt_I
   \min(t^{-2},(1-s)^{-2}) \pl \pl t^2(1-s)^2 \pl
   d\nu_1(t) d\nu_2(s)  \\
  &\le 8\pi^{-2} \intt_{\delta}^{\frac12} \intt_{\delta}^{\frac12}
  \min(t^{-2},s^{-2}) ts \frac{dt}{\sqrt{t}}
  \frac{ds}{\sqrt{s}} \lel
  16\pi^{-2} \intt_\delta^{\frac12}
  \intt_{\delta}^t  \sqrt{s} ds \p
  \frac{dt}{t\sqrt{t}} \\
  &\le
 \frac{32}{3\pi^2} \intt_{\delta}^{\frac12}
  t^{\frac32} \frac{dt}{t\sqrt{t}}
  \kl \frac{16}{3\pi^2} \pl .
  \end{align*}
Finally, we need an  $L_2$-norm estimate of $k(t,s)=(1-t)sv(t,s)$:
 \begin{align*}
  &\noo k\rrm^2_{L_2(\nu_2\ten \nu_1)}
  \kl  4 \intt_I  \min(t^{-2},(1-s)^{-2}) \pl \pl  (1-t)^2 s^2 \pl
  \frac{d\mu(t)}{(1-t)} \frac{d\mu(s)}{s}\\
  &\pl \le  8 \pi^{-2} \pl \intt_{\delta}^{\frac12} \intt_{\delta}^{\frac12}
   \min(t^{-2},s^{-2}) \frac{dt}{\sqrt{t}}
   \frac{ds}{\sqrt{s}} \! \kl\!
 16  \pi^{-2} \, \intt_{\delta}^{\frac12} \intt_{\delta}^{t}
  \frac{ds}{\sqrt{s}}
    \frac{dt}{t^2\sqrt{t}} \kl
       32 \pi^{-2} \pl \intt_{\delta}^{\frac12}
  \sqrt{t}    \frac{dt}{t^2\sqrt{t}}
  \kl  32 \pi^{-2} \delta^{-1}  .
  \end{align*}
We note that
 \for
 \lefteqn{ \intt_I f(t,s) \frac{d\mu(t)}{t}\frac{d\mu(s)}{s}
 \lel
  \intt_{\delta}^{\frac12}\intt_{\frac12}^{1-\delta} \frac{1}{ts+(1-t)(1-s)} d\mu(t)d\mu(s)
  }\\
 & & \gl  \frac{1}{2\pi^{2}}  \pl  \intt_{\delta}^{\frac12}\intt_{\delta}^{\frac12}
 \min(t^{-1},s^{-1})  \frac{dt}{\sqrt{t}} \frac{ds}{\sqrt{s}} \lel
  \frac{2}{\pi^{2}}  \pl
 \intt_{\delta}^{\frac12} (\sqrt{t}-\sqrt{\delta})
    \frac{dt}{t\sqrt{t}} \\
  & &   \gl  \frac{1}{\pi^{2}}  \pl  \intt_{4\delta}^{\frac12}
 \frac{dt}{t} \lel \frac{(-\ln 8\delta)}{\pi^2} \pl
 .
 \mel
We define $\delta=\frac{1}{ne}$,
$C=\sqrt{\max\{\frac{32}{\pi^2},\frac{16}{\pi^2},\frac{16}{3\pi^2}\}}=\frac{4\sqrt{2}}{\pi}$
and $\tilde{f}=\frac{f}{C\sqrt{-\ln\delta}}$. For $n\ge 6$ we have
$\ln ne\gl e$ and hence  \eqref{sqn1} and \eqref{sqn} are
satisfied for the corresponding quadruple
$(\tilde{f},\tilde{g},\tilde{k},\tilde{h})$. Note that
$-\ln8\delta=\ln \frac{ne}{8}\gl \frac14 \ln ne$ for $n\gl 7$. The
assertion follows from
 \begin{align*}
 \frac{-\ln 8\delta}{\pi^2 C\sqrt{\ln ne}}
 \gl \frac{\ln ne}{4\pi^2C\sqrt{\ln ne}} \lel \frac{\sqrt{1+\ln n}}{16\pi \sqrt{2}}  \pl .
 \end{align*}
For very large $n$ we can choose $C=4/\pi$ and asymptotically get
$\ge (1-\eps_n) \frac{\sqrt{1+\ln n}}{4\pi}$.
 \qd

The rest of this section is devoted to the upper estimate. For a
measure $\nu$ and positive measurable densities $g,h$, we use the
$L_p$-sum
 \[ L_2(g\nu)+_p L_2(h\nu) \lel L_2(g\nu)\oplus_p L_2(h\nu)/ker Q\]
where $Q(f_1,f_2)=f_1+f_2$. Given a measurable function $k$, we
define the norm of $k$ in $L_2(g\nu)+_p L_2(h\nu)$ as the norm of
the equivalence class $[(k,0)]$.  For $p=2$, this is  again a
Hilbert space and we have an explicit formula (see \cite[Theorem
5.2.2 and Theorem 5.4.4]{BL}).\lz

\begin{lemma} \label{L2q}  Let $\nu$ be a measure and $g,h$
strictly  positive measurable functions. For a measurable function
$k$, the norm of $k$  in $L_2(g\nu)+_2L_2(h\nu)$ is given by
 \[  \noo  k\rrm_{L_2(g\nu)+_2L_2(h\nu)}
  \lel \kla \intt \frac{|k|^2}{(g^{-1}+h^{-1})} d\nu
  \mer^{\frac12} \pl .\]
\end{lemma}\lz

Using this formula, the  following estimates are established in a
very similar way to the estimates in Lemma \ref{lest}. We leave
them to the interested reader.\lz

\begin{cor} \label{calc} Let $0<\delta<\frac12$. Then
 \begin{eqnarray}
 \noo 1_{[\delta,\frac12]}\ten
  1_{[\frac12,1-\delta]}\rrm_{L_2(\nu_1\ten
  \nu_1)+_1L_2(\nu_2\ten \nu_2)} &\le&  4\sqrt{2} \pi^{-1} (-\ln \delta)^{\frac12} \\
 \noo 1_{[\frac12,1-\delta]}\ten
  1_{[\delta,\frac12]}\rrm_{L_2(\nu_1\ten
  \nu_1)+_1L_2(\nu_2\ten \nu_2)} &\le&
  4\sqrt{2} \pi^{-1} (-\ln \delta)^{\frac12} \\
  \noo
  1_{[0,\frac12]}\ten
  1_{[0,\frac12]}\rrm_{L_2(\nu_1\ten
  \nu_1)+_1L_2(\nu_2\ten \nu_2)} &\le& 2\sqrt{2}\\
 \noo
  1_{[\frac12,1]}\ten
  1_{[\frac12,1]}\rrm_{L_2(\nu_1\ten
  \nu_1)+_1L_2(\nu_2\ten \nu_2)} &\le& 2\sqrt{2} \pl
  .
 \end{eqnarray}
\end{cor}

The next inequality yields the upper estimate in the logarithmic
`little Grothendieck inequality'.\lz

\begin{lemma}\label{logest} Let $a$ be an $n\times n$
matrix, then
 \[ \noo \summ_{ij=1}^n a_{ij} f_i\ten f_j \rrm_{
 G_n\stackrel{\wedge}{\ten} G_n} \kl 18  \pl \sqrt{1+\ln n}
  \pl
  \noo a\rrm_2\pl .\]
\end{lemma}

\begin{proof}[\bf Proof:] Given $a\in \ell_2^{n^2}$ and
$0<\delta<\frac12$, we decompose $a \lel a^1(t,s) + a^2(t,s)$
where
 \begin{align*}
   a^1(t,s) \lel
 a\ten \bigg(1_{[0,\frac12]}(t)
  1_{[0,\frac12]}(s)+ 1_{[\delta,\frac12]}(t)
  1_{[\frac12,1-\delta]}(s) +  1_{[\frac12,1-\delta]}(t)
 1_{[\delta,\frac12]}(s)+
 1_{[\frac12,1]}(t)
  1_{[\frac12,1]}(s)\bigg)
  \end{align*}
and
 \[ a^2(t,s) \lel a\ten 1-a^1(t,s)  \pl .\]
According to Corollary \ref{calc}, we get
 \for
  \noo a^1\rrm_{L_2(\nu_1\ten \nu_1;\ell_2^{n^2})+_1
  L_2(\nu_2\ten \nu_2;\ell_2^{n^2})}
  &\le& (4\sqrt{2} + 8\sqrt{2} \pi^{-1} (-\ln
  \delta)^{\frac12}) \noo a\rrm_2 \pl .
  \mel
In order to estimate $a^2$, we note that
 \for
  \| 1_{[0,\delta]}\ten 1_{[\frac12,1]}
  \|_{L_2(\nu_2)\ten_{\pi} L_2(\nu_1)} &=&
   \| 1_{[0,\delta]} \|_{L_2(\nu_2)} \|
   1_{[\frac12,1]}\|_{L_2(\nu_1)}  \kl
   \frac{2^{\frac54}\delta^{\frac14}}{\sqrt{\pi}}  \frac{2}{\sqrt{\pi}}
   \lel \frac{2^{\frac{13}{4}}}{\pi} \delta^{\frac14}   \pl
  \pl .
   \mel
Similarly, we get
 \begin{align*}
 &\max\{
  \|1_{[\delta,\frac12]}\ten 1_{[1-\delta,1]}
  \|_{L_2(\nu_2)\ten_{\pi} L_2(\nu_1)},
  \| 1_{[\frac12,1]}\ten 1_{[0,\delta]}
   \|_{L_2(\nu_1)\ten_{\pi} L_2(\nu_2)},
   \| 1_{[1-\delta,1]}\ten 1_{[\delta,\frac12]}
   \|_{L_2(\nu_1)\ten_{\pi} L_2(\nu_2)}\} \\
 &\pll  \kl \frac{2^{\frac{13}{4}}}{\pi} \delta^{\frac14}
 \pl .
 \end{align*}
Using $L_2(\ell_2^n)\ten_\pi L_2(\ell_2^n)= (L_2^c\wet L_2^r)\wet
((\ell_2^n)^c\wet (\ell_2^n)^r)$, we deduce that
 \for
  \noo a^2\rrm_{L_2(\nu_1;\ell_2^n)\ten_{\pi}
 L_2(\nu_2;\ell_2^n)+_1
 L_2(\nu_2;\ell_2^n)\ten_{\pi}
 L_2(\nu_1;\ell_2^n)}
 &\le& \frac{4\pl 2^{\frac{13}{4}}}{\pi} \pl \delta^{\frac14}
 \noo a\rrm_1 \kl
 \frac{4\pl 2^{\frac{13}{4}}}{\pi} \pl \delta^{\frac14}
 \pl \sqrt{n} \noo a\rrm_2 \pl .
 \mel
We can choose $\delta=\frac{1}{e^2n^2}$ and the assertion follows
from  Lemma \ref{start}.\qd

\begin{proof}[\bf Proof of Proposition \ref{ncalc}:] Combine Lemma
\ref{logest} and Lemma \ref{lest}.\qd

\section{K-functionals}
In section 7 we  will see that Voiculescu's inequality leads to
three terms. By duality, we have to consider three term
$K$-functionals. However, the quotient structure discussed before
only involves two terms. In this section we justify the abstract
central limit procedure relating two and three term
$K$-functionals. In the following, $N$ will be a semifinite von
Neumann algebra with a normal faithful trace $\tau$. We fix a
positive $\tau$-measurable operator $d\in L_0(N,\tau)$ with full
support. The corresponding strictly semifinite weight is given by
 \[ \varphi(x)\lel \tau(dx) \pl. \]
We will use the standard notation
 \[  n_{\varphi}\lel \{x\in N\pl|\pl \varphi(x^*x)<\infty\} \pl . \]
Similarly, we will use the notation $n_{\psi}=\{x\in N\p |\p
\psi(x^*x)<\infty\}$ for every operator-valued weight $\psi$
defined on $N$. The arguments in this section generalize to
arbitrary strictly semifinite weights, but for our applications it
suffices to consider $N=L_\infty(\tilde{\mu};M_2)$ (more precisely
$L_{\infty}(\nz\times [0,1],\tilde{\mu};M_2)$, where
$\tilde{\mu}=m\ten \mu$ is given as the tensor product of the
counting measure $m$ and the measure $\mu$ defined by
\eqref{mumeas} in section 3). The three term $K$-functional $K_t$
is defined on (a subspace of) $L_0(N,\tau)$ as follows:
 \[ \noo x\rrm_{\kz_t}\lel \inf_{x=x_1+x_2d^{\frac12}+d^{\frac12}x_3}  t^{\frac12}\noo x_1\rrm_1 +
 \noo x_2\rrm_2+ \noo x_3\rrm_2 \pl. \]
We will use the quotient map  $q_t:L_1(N)\oplus L_2(N)\oplus
L_2(N)\to L_0(N,\tau)$ given by
 \[ q_t(x_1,x_2,x_3)\lel t^{-\frac12}x_1+x_2d^{\frac12}+d^{\frac12}x_3 \pl .\]
The operator space structure of $\kz_t$ is then defined as the
quotient space
 \[ \kz_t\lel \kz_t(N,d) \lel L_1(N)\oplus_1 L_2^r(N)\oplus_1
 L_2^c(N)/\krn(q_t) \pl. \]
Let us start with some elementary properties. (Note that the
intersection in the following lemma depends on $d$.) \lz

\begin{lemma}\label{d-dual}
\begin{enumerate}
\item[i)]  The dual of $\overline{\kz_t}$ with respect to the
antilinear duality bracket is $n_{\varphi}\cap n_{\varphi}^*$
equipped with the operator space structure of $N\cap L_2^c(N)\cap
L_2^r(N)$.
\item[ii)] Let $M$ be  another von Neumann algebra and
 \[ \psi(m\ten x)\lel m\varphi(x) \]
the induced operator valued weight on $M\bar{\ten}N$.  The dual
space of $\overline{L_1(M)\wet \kz_t}$ is $n_{\psi}\cap
n_{\psi}^*$ equipped with the operator space structure of
$M\bar{\ten}N\cap M\bar{\ten}L_2^c(N)\cap M\bar{\ten}L_2^r(N)$.
\item[iii)]  Let $e_n\lel 1_{[\frac{1}{n},n]}(d)$ denote the
spectral projections of $d$. The maps  $P_n(x)=e_nxe_n$ extend to
complete contractions on $\kz_t$ such that $\bigcup_n (id\ten
P_n)(L_1(M)\wet \kz_t)$ is norm dense in $L_1(M)\wet \kz_t$.
\end{enumerate}
\end{lemma}

\begin{proof}[\bf Proof:] We follow the well-known principle in
interpolation theory that the dual (unit ball) of a sum is the
intersection (of the dual unit balls). More precisely, let
$\l:\overline{\kz}_t\to \cz$ be a linear functional. Since
$\overline{L_1(N)}^*=N$, we find an element $y_1\in N$ such that
$\l(q_t(x_1,0,0))\lel \tau(x_1^*y_1)$ for all $x\in L_1(N)$.
Similarly, we find $y_2\in L_2(N)$, $y_3(N)$ such that
 \[ \l(0,x_2,0)\lel \tau(x_2^*y_2) \quad \mbox{and} \quad
 \l(0,0,x_3)\lel \tau(x_3^*y_3) \]
for all $x_2,x_3\in L_2(N)$.  Since $q_t(xe_nd^{\frac
12},-t^{-\frac1 2}xe_n,0)=0$, we deduce that
 \[ \tau(d^{\frac1 2}e_nx^*y_1)\lel t^{\frac 12} \tau(e_nx^*y_2) \pl \]
holds for all $x$ and $n$. This implies that
$y_1d^{\frac12}=t^{\frac12}y_2$  and hence $\tau(dy_1^*y_1)=t\noo
y_2\rrm_2^2$ is finite. Similarly, we find that
$d^{\frac12}y_1\lel t^{\frac12}y_3$, and hence $y_1$ and $y_1^*$
are in $n_{\varphi}$. This yields an isometric embedding
 \[ \overline{\kz}_t^*\lel \{(y,t^{-\frac12}yd^{\frac12},
 t^{-\frac12}d^{\frac12}y)\pl| \pl y\in n_{\varphi}\cap
 n_{\varphi}^*\} \subset N\oplus_{\infty}L_2(N)\oplus_{\infty} L_2(N)
 \pl .\]
Repeating the same argument for $\overline{L_1(M)\wet \kz_t}$, we
apply \eqref{adrc} and obtain an isometric embedding
 \[ \overline{L_1(M)\wet \kz_t}^{\pl *} \subset
 M\bar{\ten}N\oplus_{\infty}
 M\bar{\ten}L_2^r(N)\oplus_{\infty}M\bar{\ten}L_2^c(N) \pl .\]
Moreover, for finite rank tensors $z=\sum_i m_i\ten x_i$ we have
 \begin{equation}\label{conn}
  \noo \psi(z^*z)\rrm_M^2\lel \noo \sum_{ik} \phi(x_k^*x_i)
  m_k^*m_i\rrm_M\lel \noo \summ_{i} m_i\ten x_id^{\frac12}
  \rrm_{M\bar{\ten}L_2^c(N)}^2 \pl .\end{equation}
By weak$^*$-density of the finite rank tensor in $n_{\psi}$ we
obtain ii). Assertion i) follows immediately by applying ii) for
the matrix algebras $M=M_n$, $\nen$. For the proof of iii), we
observe that $P_n(x)=e_nxe_n$ is a complete contraction on
$L_1(N)$, $L_2^r(N)$ and $L_2^c(N)$  for all $\nen$. Moreover,
$(P_n,P_n,P_n)(\krn(q_t))\subset \krn(q_t)$. Since $(e_n)$
converges strongly to $1$, we have point-norm convergence of
$(P_n)$ to the identity in $L_1(N)$, $L_2(N)$, respectively. This
implies that $\bigcup_n id\ten P_n(L_1(M)\wet \kz_t)$ is norm
dense. \qd


The two term $K$-functional is defined as follows:
 \[ K\lel K(N,d)\lel L_2^r(N)\oplus_1 L_2^c(N)/\krn(q)\quad ,\quad  q(x_2,x_3)\lel
 x_2d^{\frac12}+d^{\frac12}x_3 \pl .\]
The identity map $I_t:K\to \kz_t$ is completely contractive. (If
$x=x_2d^{\frac12}+d^{\frac12}x_3$, then  we may choose $x_1=0$ in
the definition of $\kz_t$.)\lz

\begin{prop}\label{complement} $K(N,d)$ is a direct limit of the
$K(e_nNe_n,e_nd)$'s.  Let $(t_k)$ be a sequence with $\lim_k
t_k=\infty$ and let $\U$ be a free ultrafilter. Then $K$ is
completely contractively complemented in $\prodd_{k,\U}
\kz_{t_k}$.
\end{prop}\lz

\begin{proof}[\bf Proof:]
We will not repeat the argument for the first assertion, which is
very similar to the proof of Lemma \ref{d-dual} iii). Similarly as
in Lemma \ref{d-dual}, we can show the dual $\overline{K}^*$ is
the subspace of $L_2^c(N)\oplus L_2^r(N)$ consisting of pairs
$(y_2,y_3)$ such that $d^{\frac12}y_2=y_3d^{\frac12}$. The linear
mapping $I=(I_{t_k}):K\to \prod_{k,\U} \kz_{t_k}$ is clearly a
contraction. The only difficulty is to construct a right inverse.
Using the $P_n$'s, we may assume that $d$ and $d^{-1}$ are
bounded. For $y\in N$ we may define
 \[ \l_{t_k}(y) (q_{t_k}(x_1,x_2,x_3))\lel \tau((t_k^{-\frac 12}x_1^* +
 d^{\frac12}x_2^*+x_3^*d^{\frac12})y) \pl .\]
Using the duality between sums and intersections, we deduce from
Lemma \ref{d-dual} that
 \begin{align*}
 &\lim_k \noo (id\ten \l_{t_k})(y)\rrm_{M_m(\overline{\kz}_{t_k}^*)}\\
 &\pl = \lim_k  \max\{t_k^{-\frac12}\noo y\rrm_{M_m(N)},\| y(1\ten
 d^{\frac12})\|_{M_m(L_2^c(N))}, \| (1\ten
 d^{\frac12})y\|_{M_m(L_2^r(N))}\} \\
 &= \max\{\|y(1\ten
 d^{\frac12})\|_{M_m(L_2^c(N))}, \| (1\ten
 d^{\frac12})y\|_{M_m(L_2^r(N))}\}
 \pl
 \end{align*}
holds for  all $y\in M_m(N)$.  This shows that the map
$\l=(\l_{t_k}):L\to \prod_{k} \overline{\kz}_{t_k}^*$, defined on
the subspace
 \[ L \lel \{(yd^{\frac12},d^{\frac12}y)\pl|\pl y\in N\}\pl \subset \pl \overline{K}^* \pl, \]
is a complete contraction. Since $Nd^{\frac12}$ is dense in
$L_2(N)$ and the mapping $T(x)=d^{\frac12}xd^{-\frac12}$ is
bounded, we deduce that $L$ is norm dense in $\bar{K}^*$. (Here we
use  that  $d$ and $d^{-1}$ are bounded; in general we have only
weak$^*$-density.) We will use the obvious inclusion $\prod_{t_k}
\overline{\kz}_{t_k,\U}^*\subset \overline{\prod_{t_k,\U}
\kz_{t_k}}^*$. By continuity we may extend $\l$ to  a complete
contraction $l:\overline{K}^*\to \overline{\prod_{t_k,\U}
\kz_{t_k}}^*$ such that
 \begin{align}\label{rrr}
 &l\big((yd^{\frac12},d^{\frac12}y)\big)\big(I((x_2,x_3)+\krn(q)\big)
  \lel  \lim_{k,\U} \l_{t_k}(y) (q_{t_k}(x_1,x_2,x_3)) \lel
 \tau((d^{\frac12}x_2^*+x_3^*d^{\frac12})y)
 \end{align}
holds for all $x_2,x_3\in L_2(N)$ and $y\in N$. For $\eta\in K^*$
we use the notation $\bar{\eta}(x)\lel \overline{\eta(x)}$. Then,
we may define the adjoint $l':\prod_{t_k,\U} \kz_{t_k} \to
K^{**}$ by
 \[  \eta(l'(\xi)) \lel \overline{l(\bar{\eta})(\xi)} \pl .\]
It is easily checked that $l'$ is also completely contractive.
Using the reflexivity of $K$ and \eqref{rrr}, we deduce that
$id_K=l'I$. \qd

In the preceding sections, we had to work with two densities. We
refer to \cite{JX2} for a more systematic  explanation of why two
densities are necessary for homogeneous operator spaces in
$QS(R\oplus C)$. Technically, two densities are easily obtained
from the classical $2\times 2$-matrix trick. Let
$(\Om,\Si,\tilde{\mu})$ be a measure space and
$N=L_\infty(\Om,\Si,\tilde{\mu};M_2)$. The natural trace is given
by $\tau(x)=\int tr_2(x(\om)) d\tilde{\mu}(\om)$. Let $d_1,d_2$ be
two non-singular densities in $L_0(N,\tau)$ of $\tau$-measurable
operators. Then the diagonal
 \[ d\lel \kla \begin{array}{cc} d_1 & 0 \\
  0& d_2 \end{array}\mer \]
belongs to the space $L_0(N,\tau)$ and $\varphi(x)=\tau(dx)$ is
faithful.\lz

\begin{lemma}\label{corner} Let $N$ be as above.
The $(1,2)$-corner of $K$ is completely complemented in $K$ and
completely isometrically isomorphic to
 \[ K(d_1,d_2) \lel L_2^r(\tilde{\mu})\oplus_1 L_2^c(\tilde{\mu})/\{(f,g)\pl|\pl fd_2^{\frac12}+d_1^{\frac12}g=0\} \pl .\]
\end{lemma}

\begin{proof}[\bf Proof:] The orthogonal projection
$\Pa(x)=\kla \begin{array}{cc} 0& x_{12}\\0&0\end{array}\mer$ is
completely contractive on $L_2^r(N)$ and $L_2^c(N)$. Moreover, we
have $(\Pa,\Pa)(\krn(q))\subset \krn(q)$ and thus $\Pa$ extends to
a complete contraction $\hat{\Pa}$ on the quotient space
$L_2^r(N)\oplus L_2^c(N)/\krn(q)$. The range of $\hat{\Pa}$ is
given by pairs  $(x,y)+\krn(q)$ such that
$x=\kla\begin{array}{cc}0 & x_{12}\\0&0\end{array}\mer$ and
$y=\kla\begin{array}{cc}0 & y_{12}\\0&0\end{array}\mer$. Then, we
observe that
 \[ xd^{\frac12}+d^{\frac12}y \lel \kla\begin{array}{cc}0 & x_{12}d_2^{\frac12}
 +d_1^{\frac12}y_{12} \\0&0\end{array}\mer
 \pl .\]
Thus $\hat{\Pa}\big(L_2^r(N)\oplus L_2^c(N)/\krn(q)\big)$ is
completely isometrically isomorphic to $K(d_1,d_2)$.\qd

\begin{lemma}\label{fin-trans}
Let $(\Om,\Si,\tilde{\mu})$ be a measure space and
$\tilde{\nu}_1=d_1^{-1}\tilde{\mu}$,
$\tilde{\nu}_2=d_2^{-1}\tilde{\mu}$. Then
 \[ L_2^c(\tilde{\nu}_1)\oplus_1 L_2^r(\tilde{\nu}_2)/\{(f,g)\pl|\pl f+g=0 \mbox{ a.e.}\} \]
is completely isometrically isomorphic to $K(d_1,d_2)$.
\end{lemma}

\begin{proof}[\bf Proof:] Let $S=\{(f,g)\pl|\pl
fd_2^{\frac12}+d_1^{\frac12}g=0\}$.  The complete isometry is
induced by  the map
  \[ \iota:L_2^c(\tilde{\nu}_1)\oplus L_2^r(\tilde{\nu}_2)\to K(d_1,d_2)\pl ,\pl  \iota(f,g)\lel
 (fd_2^{-\frac12},d_1^{-\frac12}g)+S \pl \]
because $\krn(\iota)=\{(f,g)\pl|\pl f+g=0 \mbox{ a.e.}\}$. \qd




\section{Sums of free mean zero variables}

In this section we will consider free products in the sense of
\cite{VDN} to prove the  probabilistic estimates we require. The
estimates for cb-norms are obtained by considering  free products
with amalgamation (over matrix algebras). Duality will then
provide the complementation for the $3$-term $K$-functional. Let
us recall the notion of operator-valued free probability needed in
this context. We assume that a von Neumann algebra $B$ is given
along with a family $(A_j)$ of von Neumann algebras $A_j$ which
all contain $B$. In addition, we assume that there are normal,
faithful conditional expectations $E_i:A_i\to B$. Let $M$ be a
$C^*$-algebra containing $B$ with a conditional expectation
$E:M\to B$. We also assume that $\pi_i:A_i\to M$ are
$^*$-homomorphisms such that $E\circ \pi_i=E_i$ and $\pi_i|_B=id$.
Now, the image algebras $B_i=\pi_i(A_i)$ are \emph{free over B} if
 \[ E(b_1\cdots b_n) \lel 0 \]
holds for all $n$, $b_1\in B_{i_1},...,b_n\in B_{i_n}$ with
$i_1\neq i_2\neq \cdots \neq i_n$ and $E(b_1)=\cdots =E(b_n)=0$.
The scalar case corresponds to $B=\cz$ and a state $E:M\to \cz$.
Then the $C^*$-algebra generated by the $B_i$'s is isomorphic to
the free product $\ast_{i\in I} (A_i,E_i)$ (see \cite{DykB} for
details). We are more interested in the von Neumann algebra free
product which will be described later.  We refer the reader to
\cite{Vo1,Dyk2,DykB} for a detailed description of the free
product with  amalgamation and the Fock space construction, an
essential tool for our estimates. Let us use the standard notation
 \[ {\rm \AA_i}\lel (1-E_i)(A_i) \]
for the $B$-bimodule of mean $0$ elements. Following \cite{JD} we
use the notation $L_{\infty}^c(A_i,E_i)$ for the completion of
$A_i$ with respect to the norm $\noo
x\rrm_{L_{\infty}^c(A_i,E_i)}=\|E_i(x^*x)\|^{\frac12}$. In
\cite{DykB} the space  $L_{\infty}^c(A_i,E_i)$ is denoted by $
\L_2(A_i,E_i)$. We consider the $B$ bimodules
 \[ \Hh_i  \lel L_\infty^c(A_i,E_i)\ominus B  \pl. \]
The Fock space is the $B$-bimodule
 \begin{equation}\label{hb}
  \Hh \lel B \oplus \summ_{n\ge 1,i_1\neq \cdots \neq i_n}
 \Hh_{i_1}\ten_B \cdots
 \ten_B \Hh_{i_n}
  \pl .\end{equation}
Let us denote by $Q_{i_1,...,i_n}$ the orthogonal projection onto
the submodule  $\Hh_{i_1}\ten_B \cdots \ten_B \Hh_{i_n}$.  We
denote by $Q_{\emptyset}$ the projection onto $B$. As in \cite{La}
we use the notation $\L(\Hh)$ for the $C^*$-algebra of adjointable
right module maps. Indeed, a right module map $T:\Hh\to \Hh$ is
called adjointable if there exists $S:\Hh\to \Hh$ such that
 \[ \langle S(x),y\rangle \lel \langle x,T(y)\rangle \]
for all $x,y\in \Hh$. (Here $\langle \p ,\rangle$ is the
$B$-valued sesquiliner form.) Note that the $Q_{i_1,...,i_n}$ are
adjointable $B$-bimodule maps,  and that $E(T)\lel
Q_{\emptyset}TQ_{\emptyset}$ defines a conditional expectation
from $\L(\Hh)$ onto $B$. The free product with amalgamation may be
constructed by defining the $^*$-homomorphism $\pi_i:A_i\to
\L(\Hh)$ as follows: If $a\in B$, the $\pi_i(a)$ acts by left
multiplication on $\Hh$. For $a$ with $E_i(a)=0$ and $i_1\neq i$
we have
 \[ \pi_i(a)(h_{i_1}\ten \cdots  \ten h_{i_n})\lel a
 \ten h_{i_1}\ten \cdots  \ten h_{i_n} \pl .\]
For $i_1=i$ we have
 \[ \pi_i(a)(h_{i_1}\ten \cdots  \ten
 h_{i_n})\lel (ah_{i_1}-E_i(ah_{i_1}))\ten  h_{i_2} \ten \cdots  \ten
 h_{i_n} + E_i(ah_{i_1})h_{i_2}\ten \cdots \ten h_{i_n}
  \pl .\]
Then $\ast_{i\in I}A_i$ is  defined as the $C^*$-algebra generated
by $\pi_i(A_i)$. It turns out that then the image algebras
$B_i\lel \pi_i(A_i)$ are free over $E$. The conditional version of
Voiculescu's inequality (\cite{Voir}) reads as follows:

\begin{prop} \label{Voi} Let {\rm $a_i\in {\rm \AA}_{i}$},
such that only finitely many $a_i$'s are different from $0$. Then
\[ \noo \summ_{i} \pi_i(a_{i})\rrm \kl \sup_{i} \noo a_{i}\rrm +
\noo \summ_{i}E_i(a_{i}^*a_{i}) \rrm^{\frac 12}+ \noo \summ_{i}
E_i(a_{i}a_{i}^*)\rrm^{\frac 12} \pl.
\]
\end{prop}

\noindent For our estimates, we  follow Voiculescu \cite{Voir} and
define the projections
 \[ P_i\lel \sum_{i=i_1\neq \cdots \neq i_n}Q_{i_1\cdots i_n}\pl .\]

\begin{lemma} \label{fre1} Let {\rm $a\in {\rm \AA}_i$}. Then
$(1-P_i)\pi_i(a)(1-P_i)\lel 0$.
\end{lemma}

\begin{proof}[\bf Proof:] Given $h_{i_1}\ten \cdots \ten h_{i_n}\in
 \Hh_{i_1}\ten_B \cdots  \ten_B \Hh_{i_n}
 $ and $i_1\neq i$, we observe that
 \[ h\lel \pi_i(a)(h_{i_1}\ten \cdots \ten h_{i_n})\lel a \ten h_{i_1}\cdots h_{i_n}\]
is an element of $\Hh_i\ten_B
 \Hh_{i_1}\ten_B \cdots  \ten_B \Hh_{i_n}$.
Thus $P_i(h)=h$ and $(1-P_i)(h)=0$. By linearity this yields the
assertion.\qd

\begin{cor} \label{fre2} Let  $a\in A_i$. Then $(1-P_i)\pi_i(a) (1-P_i)\lel E_i(a)(1-P_i)$.
\end{cor}

\begin{proof}[\bf Proof:] This is obvious from Lemma \ref{fre1} by writing
 \begin{align*}
 \pi_i(a)&=\pi_i(a-E_i(a))+\pi_i(E_i(a)) \lel
\pi_i(a-E_i(a))+E_i(a)\pl . \qedhere
\end{align*}\qd

We will now give the easy proof of Voiculescu's inequality:\lz

\begin{proof}[\bf Proof of \ref{Voi}:] We deduce from Lemma
\ref{fre1} that
 \for
 \summ_{i} \pi_i(a_{i}) &=& \summ_{i} P_i\pi_i(a_{i})P_i + \summ_{i}
 (1-P_i)\pi_i(a_{i})P_i + \summ_{i}
 P_i\pi_i(a_i)(1-P_i) \\
 & & \pll \pll + \summ_{i}
 (1-P_i)\pi_i(a_i)(1-P_i) \\
 &=& \summ_{i} P_i\pi_i(a_i)P_i + \summ_{i}
 (1-P_i)\pi_i(a_i)P_i + \summ_{i}
 P_i\pi_i(a_i)(1-P_i)  \pl .
 \mel
Since the $P_i$'s  are mutually orthogonal, we have
 \[ \noo \summ_{i} P_i\pi_i(a_i)P_i \rrm \lel \sup_i \noo
 P_i\pi_i(a_i)P_i\rrm \kl \sup_i \noo a_i\rrm \pl .\]
Now, we  consider the second term. By orthogonality, positivity
and the module property, we deduce from  Corollary \ref{fre2} that
 \begin{align*}
 & \noo \summ_{i}
 (1-P_i)\pi_i(a_i)P_i \rrm^2 \lel  \noo \summ_{i,l}
 (1-P_i)\pi_i(a_i)P_iP_{l}\pi_l(a_{l})^*(1-P_{l}) \rrm \\
 &  \lel  \noo \summ_{i} (1-P_i)\pi_i(a_i)P_i\pi_i(a_i)^*(1-P_i) \rrm
 \kl  \noo \summ_{i}(1-P_i)\pi_i(a_ia_i^*)(1-P_i) \rrm\\
 &  \lel  \noo \summ_{i} (1-P_i)E_i(a_ia_i^*)(1-P_i) \rrm\\
 &  \lel \noo \summ_{i} E_i(a_ia_i^*)^{\frac12}(1-P_i)E_i(a_ia_i^*)^{\frac12}
 \rrm\\
 &  \kl  \noo \summ_{i}
 E_i(a_ia_i^*)^{\frac12}E_i(a_ia_i^*)^{\frac12}  \rrm
 \lel  \noo \summ_{i} E_i(a_ia_i^*)  \rrm \pl .
 \end{align*}
The calculation for the third term is the same.\qd

The  converse of Voiculescu's inequality will be formulated in the
predual of the von Neumann algebraic free  product. This
definition follows a general scheme for $C^*$-modules over von
Neumann algebras (see \cite{Pa}). For simplicity we will assume in
the following that $\varphi$ is a normal faithful state on $B$.
Let $K$ be a $C^*$-module over a von Neumann algebra $B$. Then the
Hilbert space $K\ten_B L_2(B)$ is defined by the completion of $K$
with respect to the norm
 \[ \noo h\rrm_2 \lel \varphi(<h,h>)^{\frac12} \pl .\]
Let $T\in \L(K)$. We have
 \[ \noo Th\rrm^2 \lel \varphi(<Th,Th>) \kl \noo T\rrm^2
 \varphi(<h,h>) \pl.
 \]
Since $\varphi$ is faithful, we find  a faithful
$^*$-representation $\pi:\L(K)\to \B(K\ten_B L_2(B))$.  We apply
this to $\Hh$ from \eqref{hb} above  and  obtain the von Neumann
algebra $\bar{\ast}_{i\in I} (A_i,E_i)$ as the closure of
$\pi(\ast_{i\in I} (A_i,E_i))$ in the weak operator topology on
$\Hh\ten_B L_2(B)$. The following fact is certainly well-known to
experts. We give a proof for lack  of a reference. \lz

\begin{lemma}\label{faith} $E_\emptyset(T)\lel Q_{\emptyset}TQ_{\emptyset}$ is
a normal faithful conditional expectation from $\bar{\ast}_{i\in
I}(A_i,E_i)$ onto $B$.
\end{lemma}

\begin{proof}[\bf Proof:] The idea of the proof is very simple.
For  $a\in {\rm \AA}_i$ we define the modified right action
 \[ \pi^r_i(a) (h_{i_1}\ten \cdots \ten h_{i_m}) \lel
 \begin{cases}
 h_{i_1}\ten \cdots \ten h_{i_m} \ten a & \mbox{ if }   i_m \neq
 i\\
 0 &   \mbox{ if }   i_m
 =i
 \end{cases} \pl .\]
It is elementary to check that $\pi_i^r(a)$ commutes with
$\pi_k(a_k)$ for all $a\in {\rm \AA}_i$ and $a_k\in A_k$.
Therefore the right action commutes with the von Neumann algebra
$\bar{\ast}_{i\in I}(A_i,E_i)$. Let $C$ be the algebra generated
by the $\pi^r_i(A_i)$'s and the identity. Then $CQ_{\emptyset}$ is
dense in $\Hh$. Let $x$ be a positive element in the von Neumann
algebra generated $\ast_{i\in I} A_i$ such that
$E_\emptyset(x)=0$. Then $x^{\frac12}Q_{\emptyset}=0$. Since
$x^{\frac12}$ commutes with elements in $C$ and $CQ_{\emptyset}$
is dense in $\Hh$, we obtain $x^{\frac12}h=0$ for all $h\in \Hh$.
Thus $x=0$. The careful reader will have observed that the only
shortcoming of this argument is that the right actions
$\pi_i^r(a)$ are not necessarily continuous and therefore the
passage from elements in $\ast_{i\in I}(A_i,E_i)$ to elements in
the von Neumann algebra is not justified. In order to avoid this
difficulty, we consider the right action on $\Hh\ten_B L_2(B)$.
By assumption the states $\varphi_i\lel \varphi\circ E_i$ are
faithful. Let $a\in A_i$ be an element in the domain of
$\Delta_i^{-1/2}$, where $\Delta_i$ is the generator of the
modular group $\si_t^{\varphi_i}$. Then there exists a constant
$C$ such that
 \[ \varphi(a^*y^*ya)\kl C \varphi(y^*y) \]
holds for all $y\in A_i$.  We define the right action $R_a:K\ten_B
L_2(B)\to K\ten_B A_i\ten_B L_2(B)$ by $R(a)x\lel x\ten a$. We
observe that
 \[ \noo R(a)x\rrm_2^2 \lel \varphi(a^*<x^*,x>a)\kl C
 \varphi(<x^*x>)\lel C\noo  x\rrm_{K\ten_BL_2(B)}^2 \pl .\]
In particular, $R(a)$ is continuous on $K\ten_B L_2(B)$. Let
$D_i^{\frac12} $ be  the cyclic and separating vector for
$\varphi_i$ in $L_2(A_i)$. If $(a_s)$ is a bounded net converging
strongly to $a$, then
  \[ \lim_s \noo h\ten (a-a_s)\rrm_2^2
   \lel \lim_s \noo <h,h>^{\frac12}(a-a_s)D_i^{\frac12}\rrm_2 \lel
   0 \pl .\]
It is well-known  (see \cite{KR}) that the algebra of analytic
elements  is strongly dense in $A_i$. Moreover, since
$\si_t^{\varphi}\circ E\lel E\circ \si_t^{\varphi_i}$, we see that
for analytic $x$ the expectation $E(x)$ is also analytic.
Therefore elements in ${\rm \AA}_i$ can be approximated by
analytic elements in ${\rm \AA}_i$.  We  replace the algebra $C$
from above  by the algebra generated by right actions
$\pi^r_i(a_i)$ where the $a_i$'s are analytic elements in ${\rm
\AA}_i$, $i\in I$.  Then $CQ_{\emptyset}$ is dense in $\Hh\ten_B
L_2(B)$ and the argument at the beginning of this proof is
justified.\qd

\newcommand{\ems}{\emptyset}
In the following, we use $\M_{free}\lel \bar{\ast}_{i\in
I}(A_i,E_i)$ and denote by  $\phi \lel \varphi \circ E_{\ems}$ the
normal faithful state on $\M_{free}$. We denote by $D_{\phi}$
($D_{\varphi}$, $D_{\varphi_i}$) the density of $\phi$ in
$L_1(\M_{free})$ (the density of $\varphi$ in $L_1(B)$, the
density of $\varphi_i=\varphi\circ E_i$ in $L_1(A_i)$,
respectively). These notation suggest that all these densities and
states are compatible. Indeed, we deduce from \eqref{cs} that
  \[  \si_t^{\varphi}\circ E\lel E \circ \si_t^{\varphi_i} \pl .\]
The same argument holds also for the inclusion $\pi_i(A_i)\subset
\bar{\ast}_{i\in I} (A_i,E_i)$. Moreover, following \cite[Lemma
1.1]{DykB}, we have a conditional expectation $\E_i:\M_{free}\to
\pi_i(A_i)$ given by
 \[ \E_i(x) \lel (Q_i+Q_{\ems})x(Q_i+Q_{\ems})\in  \L(B\oplus
 {\rm \AA}_i) \lel \L(L_{\infty}^c(A_i,E_i)) \pl \]
such that
 \[ \E_i(\pi_i(a_1))a_2\lel a_1a_2\]
holds for all $a_1,a_2\in A_i$. Moreover, if $\pi_L:A_i\to
\L(A_i)$ denotes the left action of $A_i$ on
$\L_{\infty}^c(A_i,E_i)$, then $\E_i(\M_{free})\subset
\pi_L(A_i)$. Obviously we have $E_i\circ \E_i \lel E_{\ems}$  and
hence $\phi \lel \varphi_i \circ \E_i$. Applying \eqref{cs} once
more, we deduce that
 \begin{equation} \label{modgr}
   \si_t^{\varphi_i}\circ \E_i\lel \E_i\circ \si_t^{\phi}
 \end{equation}
holds for all $i\in I$. Now, the converse of Voiculescu's
inequality is easily verified:\lz

\begin{lemma} \label{conv} Let $(a_i)$ be a sequence of elements
{\rm $a_i\in {\rm \AA}_i$} such that only finitely many elements
are non zero. Then
 \begin{eqnarray}
  \noo \summ_{i} \pi_i(a_i)D_{\phi}\rrm_{L_1(\M_{free})} &\le& \summ_{i} \noo
  a_iD_{\varphi_{i}} \rrm_{L_1(A_{i})}\pl , \label{s1}
   \\
  \noo \summ_{i} \pi_i(a_i)D_{\phi}\rrm_{L_1(\M_{free})} &\le& \noo
  (\summ_{i}  D_{\varphi} E_i(a_i^*a_i)D_{\varphi})^{\frac12}
    \rrm_{L_1(B)}\pl , \label{s2}
   \\
     \noo \summ_{i} D_{\phi}\pi_i(a_i)\rrm_{L_1(\M_{free})} &\le& \noo
  (\summ_{i} D_{\varphi} E_i(a_ia_i^*)D_{\varphi})^{\frac12}
    \rrm_{L_1(B)}\pl . \label{s3}
  \end{eqnarray}
\end{lemma}

\begin{proof}[\bf Proof:] Using  $\varphi\circ E_{\ems}=\phi$, we have a natural
family of embeddings $\iota_{p}:L_p(B)\to L_p(\M_{free})$
satisfying
 \begin{equation}\label{theta} \iota_p(D_{\varphi}^{\frac{1-\theta}{p}}bD_{\varphi}^{\frac{\theta}{p}})
 \lel
 D_{\phi}^{\frac{1-\theta}{p}}bD_{\phi}^{\frac{\theta}{p}}
 \end{equation}
for all $b\in B$ (see \cite{JX}). Since $\varphi_i\circ \E_i=\phi$
we also have a family of isometric embeddings $\iota_p:L_p(A_i)\to
L_p(\M_{free})$ such that
 \begin{equation}\label{theta2} \iota_p(D_{\varphi_i}^{\frac{1-\theta}{p}}aD_{\varphi_i}^{\frac{\theta}{p}})
 \lel
 D_{\phi}^{\frac{1-\theta}{p}}a D_{\phi}^{\frac{\theta}{p}}
 \end{equation}
for all $a\in A_i$. The triangle inequality implies \eqref{s1}.
For the proof of \eqref{s2}, we consider $z=\sum_{i} \pi_i(a_i)$
with $a_i\in {\rm \AA}_{i}$. Note that by freeness
 \[ E_\ems(z^*z) \lel \summ_{j,i} E_{\ems}(\pi_j(a_j)^*\pi_i(a_i))
 \lel \summ_{i} E_i(a_i^*a_i) \pl .\]
Since $E_\ems$ is a conditional expectation we deduce from $\noo
\E_{\ems}(x^*x)\rrm_{\infty}^{1/2}\le \noo x\rrm_{\infty}$ and by
duality (see \cite[Corollary 2.12]{JD}) that
 \begin{equation}\label{jdd}
 \noo D_{\phi}^{\frac12}E_\ems(z^*z)
 D_{\phi}^{\frac12}\rrm_{\frac12} \kl \noo
  D_{\phi}^{\frac12}z^*zD_{\phi}^{\frac12} \rrm_{\frac12} \pl.
  \end{equation}
We deduce with \eqref{theta} that
 \begin{align*}
 \noo \summ_{i} \pi_i(a_i) D_{\phi}\rrm_1 &\le  \noo \summ_{i} D_{\phi}
 E_i(a_i^*a_i)D_{\phi} \rrm_{\frac12}^{\frac12} \lel
  \noo \summ_{i}  D_{\varphi}
 E_i(a_i^*a_i)D_{\varphi} \rrm_{\frac12}^{\frac12}
  \pl .
  \end{align*}
Assertion \eqref{s3}  follows by taking adjoints.\qd

We now reformulate these inequalities  in terms of a
complementation result. Let $\N$ be a von Neumann algebra and
$\E:\N\to B$ be a normal faithful conditional expectation onto a
von Neumann subalgebra $B$. Let  $\varphi$ be a normal faithful
state on $B$. We denote by $D_{\varphi\circ \E}$ the density of
$\varphi\circ \E$. As for the three term K-functional, we define a
new norm on $L_1(\N)$ by
 \[ \noo x\rrm_{\kz_n(\N,\E)} \lel \inf_{x=x_1+x_2+x_3} n\noo
 x_1\rrm_1+ \sqrt{n} \noo x_2\rrm_{L_1^c(\N,\E)} + \sqrt{n}\noo
 x_3\rrm_{L_1^r(\N,\E)} \pl .\]
Following \cite{JD}, the space $L_1^c(\N,\E)$ is defined as the
closure of $\N D_{\varphi \circ \E}$ with respect to the norm
 \[ \noo zD_{\varphi \circ\E} \rrm_{L_1^c(\N,\E)} \lel \noo D_{\varphi \circ\E}
 \E(z^*z)D_{\varphi\circ \E} \rrm_{\frac12}^{\frac12} \pl .\]
The space $L_1^r(\N,\E)$ is the space of adjoints of elements in
$L_1^c(\N,\E)$ defined by the norm
 \[ \noo x\rrm_{L_1^r(\N,\E)}\lel \noo x^*\rrm_{L_1^c(\N,\E)} \pl .\]
For both spaces we have contractive injective inclusions
$L_1^c(\N,\E)\subset L_1(\N)$ and $L_1^r(\N,\E)\subset L_1(\N)$
(see \eqref{jdd}). Therefore $\kz_n(\N,\E)$ is well-defined. The
following observations are immediate consequences of
\cite[Corollary 2.12]{JD}. There we used the \emph{antilinear}
duality bracket. However, the adjoint map $J(x)=x^*$ is an
isometry on $\kz_n(\N,\E)$ and thus we may work with the usual
trace duality. \lz

\begin{lemma}\label{dualex} The dual of  $\kz_n(\N,\E)$ with respect
to the duality bracket
 \[ (x,y)_n\lel n tr(xy) \]
is $\N$ equipped with the norm
 \[ \noo y\rrm_{\kz_n(\N,\E)}\lel \max\{\noo y\rrm_\N,n^{-\frac12} \noo
 \E(y^*y)\rrm^{\frac12},n^{-\frac12} \noo \E(yy^*)\rrm^{\frac12}\}
 \pl .\]
\end{lemma}\lz

In Voiculescu's inequality it is very important to work with mean
$0$ elements. This will be achieved by a standard symmetrization
process. We define $A_i=\ell_{\infty}^2(\N)\lel \N\oplus \N$. All
the conditional expectations $E_i$ coincide with $E$ defined by
 \[ E(x,y)\lel \frac{\E(x)+\E(y)}{2} \pl .\]

\begin{prop} \label{comp1}
The space $\kz_{n}(\N,\E)$ is $3$-complemented in
$L_1(\bar{\ast}_{1\le i\le n}(\N\oplus \N,E))$.
\end{prop}

\begin{proof}[\bf Proof:] Let $D_{\phi}$, $D_{\varphi}$,
$D_{\varphi \circ \E}$ and  $D_{\varphi\circ E}$  be the density
of $\phi$, $\varphi$, $\varphi\circ \E$ and $\varphi\circ E$ on
$\M_{free}$, $B$, $\N$ and $A$, respectively. We may identify the
space $L_1(A)$ with $L_1(\N)\oplus_1 L_1(\N)$. Then the state
$\varphi\circ E$ is given  by
 \[ \varphi \circ E(x,y)\lel \frac{tr(xD_{\varphi \circ \E})+
 tr(yD_{\varphi \circ \E})}{2} \pl .\]
Therefore, we may assume that  $D_{\varphi\circ E}\lel (\frac12
D_{\varphi \circ \E},\frac12 D_{\varphi \circ \E})\in
L_1(\N)\oplus_1 L_1(\N)$. For $x\in \N$ we define $\eps x\lel
(x,-x)\in A$. The embedding  will be realized by the map
$T:L_1(\N)\to L_1(\M_{free})$  given by
 \[ T(xD_{\varphi\circ \E}) \lel \summ_{i=1}^n \pi_i(\eps x)D_{\phi}
 \pl .\]
According to \eqref{s1}, we have
 \[ \noo T(xD_{\varphi})\rrm_1 \kl  n \noo (\eps x) D_{\varphi \circ E} \rrm_{L_1(A)}
 \lel  \frac n2 (\noo xD_{\varphi}\rrm_{L_1(\N)}+ \noo
 -xD_{\varphi}\rrm_{L_1(\N)}) \lel n \noo xD_{\varphi}
 \rrm_{L_1(\N)}
  \pl .\]
Similarly, we deduce from the fact that $\eps x$ has mean $0$ and
\eqref{s2}  that
 \for
 \noo T(xD_{\varphi})\rrm_1 &\le  \noo D_{\varphi} \summ_{i=1}^n
 E_i(\pi_i(1\ten x^*x)) D_{\varphi} \rrm_{L_\frac12(B)}^{\frac12}
 \lel      \sqrt{n}  \noo (D_{\varphi} \E(x^*x)D_{\varphi})^{\frac12}
  \rrm_{L_1(B)}
  \pl .
   \mel
For an analytic element $x\in \N_a$, we deduce from \eqref{modgr}
that
 \begin{eqnarray}\label{rtrans}
 T(D_{\varphi \circ \E} x) &=& T(\si_{-i}^{\varphi}(x)D_{\varphi}^\N)
 \lel \summ_{k=1}^n
 \si_{-i}^{\phi}(\pi_i(\eps x))D_{\phi} \lel
  \summ_{k=1}^n D_{\phi} \pi_i(\eps x) \pl .
 \end{eqnarray}
By continuity this equality extends to all elements $x\in \N$.
Hence, \eqref{s3} implies the missing inequality and we deduce
that
 \[ \noo T:\kz_{n}(\N,\E) \to L_1(\M)\rrm \kl 1 \pl .\]
We define the map $S:\N\to \M$ by
 \[ S(y) \lel  \summ_{i=1}^n \pi_i(\eps y) \pl .\]
According to Proposition \ref{Voi}, we have
 \[ \noo S:\kz_{n}(\N,\E)^*\to \M_{free}\rrm\le 1 \pl .\]
Moreover, using trace duality
$(x,y)_{\M_{free}}=tr_{\M_{free}}(xy)$, we get
  \for
  \lefteqn{  (T(xD_{\varphi}),S(y))_{\M_{free}}
  \lel    tr_{\M_{free}}(T(xD_{\varphi})S(y)) }\\
   & &  \lel \summ_{i=1}^n tr_{\M_{free}}(\pi_i(\eps x)D_{\phi}\pi_i(\eps y))
   \lel \summ_{i=1}^n \varphi(yx) \lel
   (xD_{\varphi},y)_n \pl .
   \mel
We denote the restriction of $S^*$ to $L_1(\M_{free})$ by
$S'=S^*|_{L_1(\M_{free}}$. Then we obtain a map
$S':L_1(\M_{free})\to \kz_n(\N,\E)^{**}$ such that $S'T$
coincides with the natural inclusion map of $\kz_n(\N,\E)$ in its
bidual. We want to show that $S'T=id_{\kz_n(\N,\E)}$. It suffices
to show that $S'(L_1(\M_{free}))\subset \kz_n(\N,E)$. For the
proof of this inclusion, we observe that $\ast_{i\in I_n}(A,E_i)$
is strongly dense in $\M_{free}$ and hence $\ast_{i\in
I_n}(A,E_i)D_{\phi}$ is norm dense in $L_1(\M_{free})$. Thus it
is enough  to consider elements of the form $x=zD_{\phi}$, where
$z= \pi_{j_1}(a_1)\cdots \pi_{j_m}(a_m)$, $a_k\in {\rm
\AA}_{j_k}$. Then, we have
  \[ (S'(x),y) \lel tr(xS(y)) \lel
  \summ_{i=1}^n tr(zD_{\phi}\pi_i((\eps y)))
 \lel \summ_{i=1}^n \phi(\pi_i(\eps y)\pi_{j_1}(a_1)\cdots
  \pi_{j_m}(a_m)D_{\phi})
  \pl .\]
Thus for $S'(x)\neq 0$ to hold we must have $m=1$ and we may
assume that $x=\pi_k((a_1,a_2))D_{\phi}$ for some $1\le k\le n$.
In this case we obtain
 \begin{equation}\label{explicit}
  (S'(x),y)  \lel \phi(\pi_k(\eps y) \pi_k(a_1,a_2))
 \lel \frac12 (\varphi(ya_1)-\varphi(ya_2)) \pl .
 \end{equation}
This implies that $S'(\pi_k((a_1,a_2))D_{\phi})\lel \frac12
(a_1-a_2)D_{\varphi}$. By density $S'(L_1(\M_{free}))\subset
L_1(N)$. \qd

A natural example of freeness with amalgamation is given by tensor
products.
 \lz

\begin{exam}\label{matmod} Let $(C_i)_{i\in I}$ be von Neumann
algebras and $(\phi_i)$ a family of normal faithful states.  Let
$B$ be a von Neumann algebra. Let $\pi_i:C_i\to \bar{\ast}_{i\in
I}(C_i,\phi_i)$ be the embedding and $\phi=\ast_{i\in I} \phi_i$
be the free product state and
$E_{free}:B\bar{\ten}\bar{\ast}_{i\in I}(C_i,\phi_i)\to B$ be
given by $E_{free}(x\ten y)\lel x\phi(y)$. Then $(id \ten
\pi_i)(B\bar{\ten}C_i)$ are free over $E_{free}$.
\end{exam}

\begin{proof}[\bf Proof:] Let $a=\sum_{k=1}^m x_k\ten y_k \in B\ten
C_i$. We observe that
 \[ \pi_i(a)-E(\pi_i(a)) \lel \sum_{k=1}^m x_k \ten (\pi_i(y_k)-\phi_i(y_k)1) \pl .\]
This shows that $(id-1\ten \phi_i)(B\bar{\ten} C)\cap B\ten C\lel
B\ten \circled{C}_i$. Thus,  given $a_1,...,a_n$ such that $a_j\in
(id-1\ten \phi_{i_j})(B\bar{\ten} C)\cap B\ten C$, we may write
$a_j=\sum_{k=1}^{m} x_{jk}\ten \pi_{i_j}(y_{jk})$ with $y_{jk}\in
\circled{C}_{i_j}$. (The same $m$ is achieved by adding $0$'s.) If
$i_1\neq \cdots \neq i_n$ the freeness of the $\pi_i(C_i)$'s
implies that
 \begin{align*}
 E(a_1\cdots a_n) &= \summ_{k_1,...,k_n=1}^m
 x_{1k_1}\cdots x_{nk_n}
 \phi(y_{1k_1}\cdots y_{nk_n}) \lel 0 \pl .
 \end{align*}
By Kaplansky's density theorem, the unit ball of
$B\ten_{\min}C_{i_j}$ is strongly dense in $B\bar{\ten}C_{i_j}$.
Therefore  mean $0$ elements $a_j$ in $B\bar{\ten}C_{i_j}$ may be
approximated  in the strong operator topology by bounded nets
$a_j(\al)$ of elements in $B \ten_{\min} \circled{C}_{i_j}$. By
continuity this implies that
 \begin{align*}
  E(a_1\cdots a_n) &=  \lim_{\al} E(a_1(\al)\cdots
 a_n(\al))\lel 0 \pl.  \qedhere \end{align*}
\qd

According to \cite{Vo1} it is well-known that free products of von
Neumann algebras with states are  in general not semifinite.
Therefore, we have to work with Haagerup $L_p$-spaces in this
context. We will now describe the natural isomorphism between the
$3$-term $K$-functional in the setting of Haagerup $L_p$-spaces
and in the setting of semifinite $L_p$-spaces. \lz

\begin{lemma}\label{ident} Let $N$ be a semifinite von Neumann
algebra with trace $\tau_N$. Let $\psi$ be a faithful, normal
state on $N$ with density $d_{\psi}\in L_1(N,\tau_N)$. Let $B$ be
$\si$-finite and semifinite. Then there is a natural isomorphism
 \[ \kz_n(B\bar{\ten} N,1\ten \psi) \cong L_1(B)\wet
 \kz_n(N,d_{\psi}) \pl .\]
\end{lemma}

\begin{proof}[\bf Proof:] The trace on $B$ is denoted
by $\tau_B$.  Let $d_{\varphi}\in L_1(B,\tau_B)$ be the density of
a faithful normal state $\varphi$ on $B$. The density
$D_{\varphi\ten \psi}\in L_1(B\bar{\ten} N)$ of $\varphi\ten \psi$
has full support because $\varphi\ten \psi$ is faithful. The
complete isometry $I: L_p(B\bar{\ten}N)\to L_p(B\bar{\ten}
N,\tau_B\ten\tau_N)$ between the Haagerup $L_p$-space and the
semifinite $L_p$-space is given by
 \[ I(zD_{\varphi\ten \psi}^{\frac1p}) \lel
 z(d_{\varphi}^{\frac{1}{p}}\ten d_{\psi}^{\frac{1}{p}})
 \pl
 .\]
Since $I$ commutes with the modular group (see also \cite{JX}) we
have
 \[ I(D_{\varphi\ten \psi}^{\frac{1-\theta}{p}}zD_{\varphi\ten
 \psi}^{\frac{\theta}{p}}) \lel
 (d_{\varphi}^{\frac{1-\theta}{p}}\ten
 d_{\psi}^{\frac{1-\theta}{p}})
 z (d_{\varphi}^{\frac{\theta}{p}}\ten
 d_{\psi}^{\frac{\theta}{p}}) \pl .\]
We use $\E\lel 1\ten \psi$ for the conditional expectation from
$B\bar{\ten}N$ onto  $B\ten 1\subset N$. In section 1 (see
(1.10)) we discussed the $L_{\frac12}(B)$-valued extension of the
scalar product on the Hilbert space $H=L_2(N)$. Given
$b_{1},b_2\in B$ and $n_1,n_2\in N$, we find
 \begin{align*}
  (b_1d_{\varphi} \ten n_1d_{\psi}^{\frac12},b_2d_{\varphi} \ten n_2d_{\psi}^{\frac12})
 &= d_{\varphi}  b_1^*b_2d_{\varphi} \pl
 (n_1d_{\psi}^{\frac12},n_2d_{\psi}^{\frac12})\lel
   d_{\varphi}  b_1^*b_2d_{\varphi} \pl \psi(n_1^*n_2) \\
   &=  d_{\varphi}\E((b_1\ten n_1)^*(b_2\ten
 n_2))d_{\varphi} \pl .
 \end{align*}
Then $\tau(d_{\psi})=1$ and \eqref{hc} imply that
 \begin{align*}
  \noo zD_{\varphi\ten \psi}\rrm_{L_1^c(B\bar{\ten}N,\E)}
   &= \noo D_{\varphi\ten \psi}\E(z^*z)D_{\varphi\ten
   \psi}\rrm_{\frac12}^{\frac12} \lel
   \noo (d_{\varphi}\ten d_{\psi})
  \E(z^*z) (d_{\varphi}\ten d_{\psi})
  \rrm_{\frac12}^{\frac12} \\
   &= \noo d_{\varphi} \E(z^*z)d_{\varphi}
   \rrm_{\frac12}^{\frac12}\lel
   \noo (z(d_{\varphi}\ten d_{\psi}^{\frac12}),
      z(d_{\varphi}\ten d_{\psi}^{\frac12}))
      \rrm_{L_{\frac12}(B)}^{\frac12}
      \\
   &= \noo z(d_{\varphi}\ten d_{\psi}^{\frac12})\rrm_{L_1(B)\wet
   L_2^r(N)} \pl .
 \end{align*}
This shows that $\iota_r(x) \lel I(x(1\ten d_{\psi}^{\frac12}))$
is an isometry between $L_1(B)\wet L_2^r(N)$ and
$L_1^c(B\bar{\ten}N,\E)$. Using adjoints, we see that
$\iota_c(x)\lel I((1\ten d_{\psi}^{\frac12})x)$ also is an
isometry. Since these isometries are compatible with the map $q_n$
defined before Proposition \ref{d-dual}, we deduce that the map
$\iota_n(x)\lel \sqrt{n}I(x)$ yields  an isometric isomorphism
between $\kz_n(B\bar{\ten}N)$ and
 $L_1(B)\wet \kz_n(N,d_{\psi})$.\qd

\begin{cor}\label{comp2} Let $N$ be a semifinite von Neumann algebra and let $\psi$ be a
faithful  normal state with density $d_{\psi}$. Let
$\tilde{\psi}(x,y)\lel \frac12(\psi(x)+\psi(y))$ be the
corresponding state on $N\oplus N$. Then $\kz_{n}(N,d_{\psi})$ is
$3$-completely complemented in $L_1(\bar{\ast}_{1\le i\le
n}(N\oplus N,\tilde{\psi}))$ for all $\nen$.
\end{cor}\lz

\begin{proof}[\bf Proof:] According to Lemma \ref{ident} with $B=\cz$ it suffices to show the assertion for
$\kz_n(N,\psi)$. Here $\psi$ is considered as a conditional
expectation onto $\cz 1$ and the operator space structure is the
one given by the isomorphism in Lemma \ref{ident}. We define
$C=N\oplus N$ and denote by $\pi_i:C\to \bar{\ast}_{1\le i\le n}
(C,\tilde{\psi})$ the natural embeddings. We shall write
$\phi=\ast_{1\le i\le n} \tilde{\psi}$ for the free product state.
The map $T:\kz_n(N,\psi)\to L_1(\bar{\ast}_{1\le i\le
n}(C_i,\psi))$ is given by Proposition \ref{comp1}:
 \[ T(xD_{\psi})\lel \summ_{i=1}^n \pi_i(\eps x)D_{\phi}
 \pl .\]
In order to show that the maps $T$ and $S'$ obtained in
Proposition \ref{comp1} are indeed completely bounded, we use
operator-valued free products with respect to $B=M_m$. Then we
shall use $\N_m=M_m(N)$ with  the conditional expectation
$\E_m(x\ten y)\lel \psi(y)(x\ten 1)$. Note that  $\N_m\oplus
\N_m=M_m(C)$ and the conditional expectations $E_i:M_m(C)\to M_m$
given by $E_i(x\ten y)\lel \tilde{\psi}(y)x$ are compatible with
the definitions before Proposition \ref{comp1}. According to
Example \ref{matmod} we know that the algebras $(id\ten
\pi_i)(M_m(C))$ are free over $1\ten \phi$. $M_m(\bar{\ast}_{1\le
i \le n}(C,\tilde{\psi}))$ is generated as a von Neumann algebra
by the algebras $(id\ten \pi_i)(M_m(C))$. Hence
$M_m(\bar{\ast}_{1\le i \le n}(C,\tilde{\psi}))$ and
$\bar{\ast}_{1\le i\le n}(M_m(C),E_i)$ are isomorphic. Therefore
we are in a position to apply Proposition \ref{comp1}. We use the
normalized trace $\varphi_m(x) \lel \frac{1}{m} tr(y)$ on
$B=M_m$. We observe that $\varphi_m\circ E$ restricted to $1\ten
N$ induces $\psi$ and hence $D_{\varphi_m \circ \E}\lel
D_{\varphi_m\ten \psi}$. Moreover, in the Haagerup $L_1$-spaces
$L_1(M_m(N),\varphi_m\ten \psi)$, $L_1(M_m(\bar{\ast}_{1\le i\le
n}(C_i,\tilde{\psi})),\varphi_m\ten\phi)$, we have
$D_{\varphi_m\ten \psi}\lel 1 \ten  D_{\psi}$, $D_{\varphi_m\ten
\phi}\lel 1 \ten D_{\phi}$, respectively. We denote by $T^{(m)}$,
$S^{(m)}$ the maps constructed in Proposition \ref{comp1} for
$\kz_n(\N_m,\E_m)$. We find that
 \[ (id_{L_1(M_m)}\ten T)(x(1\ten D_{\psi})) \lel
 \summ_{i=1}^n (id\ten \pi_i)(\eps x) (1\ten
 D_{\phi})\lel
 T^{(m)}(xD_{\varphi_m\ten \psi}) \pl .\]
Hence $\noo id_{L_1(M_m)}\ten T\rrm=\noo T^{(m)}\rrm\le 1$ and $T$
is a complete contraction. Using the concrete form of $S^{(m)'}$
constructed in Proposition \ref{comp1}, we find that
$S^{(m)'}=id_{L_1(M_m)}\ten S^{'}$ and hence $\noo S'\rrm_{cb}\le
3$. \qd

\begin{theorem}\label{Ko} Let
$N$ be a semifinite von Neumann algebra with trace $\tau$ and $d$
a positive density in $L_0(N,\tau)$ with full support. Then
$K(N,d)$ is $3$-completely complemented in the predual of a von
Neumann algebra $\M$. If $N$ has QWEP, then there is such an $\M$
with QWEP.
\end{theorem}\lz

\begin{proof}[\bf Proof:] According to Proposition  \ref{complement} we
have  $K(N,d)=\lim_m K(N,e_md)$. Here  $e_m=1_{[\frac1m,m]}(d)$
are the spectral projections of $d$. Using an ultraproduct, Lemma
\ref{ultraqwep}  and the reflexivity of $K$, it therefore suffices
to show the assertion for a density with $\tau(d)<\infty$. By
normalization we may assume $\tau(d)=1$. We define the state
$\psi(x)=\tau(dx)$. According to Lemma \ref{complement}, $K(N,d)$
is $1$-completely contractively complemented in
 \[ \prodd_{n,\U} \kz_{n}(N,d)  \pl .\]
By Corollary  \ref{comp2}, $\prod_{n,\U}\kz_n(N,d)$ is
$3$-completely contractively complemented in
 \begin{eqnarray}\label{frr}
  \M_* &=& \prodd_{n,\U} L_1(\bar{\ast}_{1\le i\le n}(N\oplus
  N,\tilde{\psi}))
  \end{eqnarray}
for every free ultrafilter $\U$ on the integers. This completes
the proof in the general case. If we assume in addition that $N$
is QWEP, then $e_mNe_m$ is QWEP for every $m\in \nz$. According to
Theorem \ref{freeQWEP} (below), the von Neumann algebras
$\bar{\ast}_{1\le i\le n}(N\oplus N,\tilde{\psi})$ are  QWEP.
Since $A$ QWEP implies that $A^{op}$ QWEP, this implies with
\eqref{opp} that $\M$ from \eqref{frr} is QWEP. Using Lemma
\ref{ultraqwep} (below) the  QWEP property is stable under
ultraproducts and the assertion follows.\qd

By standard  properties of the projective tensor product, we
obtain the following application of Theorem \ref{Ko}.   \lz

\begin{cor} \label{kick1} Let $N$  be a semifinite von
Neumann algebra with trace $\tau$ such that $N$ is QWEP. Let $d$
be a density in $L_0(N,\tau)$ with  full support. Then
 \[ K(N,d)\wet K(N,d) \]
is $9$-completely complemented in
 \[ \M_*\wet \M_* \]
for some von Neumann algebra $\M$ with QWEP.
\end{cor}\lz

\begin{samepage} \begin{proof}[\bf Proof of Proposition \ref{kick00} and Theorem 3:] We
consider $\Om\lel \nz\times [0,1]$ and $\tilde{\mu}=m\ten \mu$,
with  $m$ the counting measure on $\nz$ and
$d\mu(t)=dt/{(\pi\sqrt{t(1-t)})}$.  By Lemma \ref{corner} and
Lemma \ref{fin-trans}, we see that $G$ is completely contractively
complemented in $K(M_2(L_{\8}(\nz\times[0,1],\tilde{\mu})),1\ten
d)$ where $d(t)=\kla
\begin{array}{cc} t &0 \\ 0 & 1-t \end{array}\mer$.
Hence the assumptions of  Theorem \ref{Ko} and Corollary
\ref{kick1} are satisfied. This completes the proof of Theorem 3.
The lower estimate in Proposition \ref{kick00} follows immediately
by complementation. \qd \end{samepage} \lz

\begin{rem} {\rm Using Speicher's central limit theorem (see \cite{Sp}) it is not too
difficult to identify the underlying von Neumann algebra of  the
embedding of OH (and indeed an arbitrary quotient of $R\oplus C$)
as a free quasi-free state factor of Shlyakhtenko \cite{S1}. In
case of OH this turns out to be a free quasi-free factor of type
III$_1$. After an early draft of this paper circulated, Pisier
\cite{Ps6} found a more direct approach to Theorem \ref{Ko}
without using the three term $K$-functional. This yields an easier
way to identify OH in the predual of a type III factor. For
related results in the $L_p$-setting and more information  on  the
possible types of these factors, we refer to \cite{X1}. The
approach via the three term $K$-functional is used in \cite{JO}
for a `concrete' embedding of OH in the predual of the hyperfinite
III$_1$ factor.}
\end{rem}\lz

At the end of this section, we will provide the results on the
QWEP property needed above. \lz

\begin{lemma}\label{ultraqwep}  Let $(M^s)$ be  a family of QWEP von
Neumann algebra. Then the von Neumann algebra $\M=(\prod_{s,\U}
M^s_*)^*$ also has QWEP.
\end{lemma}

\begin{proof}[\bf Proof:] According to Kirchberg \cite{Ki}, we know that
$\prod M^s$ is  QWEP and thus $(\prod M^s)^{**}$ is QWEP.
Following Groh \cite{Gr}, we observe that  the space of
functionals $\prod_{s,\U} M^s_*$ on $\prod M^s$ is left and right
invariant under the action of $\prod M^s$. Hence there is a
central projection $z_\U$ such that $\M\cong z_\U(\prod
M^s)^{**}$. Thus $\M$ is QWEP. \qd

\begin{theorem}\label{freeQWEP} Let $N$ and $M$ be von Neumann algebras with QWEP
and let  $\phi$, $\psi$ be normal faithful states on $N$, $M$,
respectively. Then the von Neumann algebra $(N,\phi)\bar{\ast}
(M,\psi)$ is QWEP.
\end{theorem}\lz

We need some preparation. The following result can alternatively
be proved using Dykema's deep analysis of free product of matrix
algebras (see \cite{Dyk}). We prefer a more direct approach using
results from Shlyakhtenko \cite{S1}. \lz

\begin{lemma} \label{fidim} Let $A_1$, $A_2$ be matrix algebras with
normal faithful states $\phi_1$ and $\phi_2$. Then the free
product  $(A_1,\phi_1)\bar{\ast}(A_2,\phi_2)$ is QWEP.
\end{lemma}

\begin{proof}[\bf Proof:] First we observe that we  may assume
$A_1=A_2$ and $\phi_1=\phi_2$. Indeed, we consider
$(A,\phi)=(A_1\ten A_2,\phi_1\ten \phi_2)$. We denote by $\pi_1$
(and  $\pi_2$) the embedding of $A$ in the first (respectively
second) component of the free product. Then the von Neumann
algebra generated by $\pi_1(A_1\ten 1)$ and $\pi_2(1\ten A_2)$ is
isomorphic to the free product
$(A_1,\phi_1)\bar{\ast}(A_2,\phi_2)$ and invariant under the
action of the modular group of the free product state
$\phi\ast\phi$. By Takesaki's theorem (see e.g. \cite[Theorem
10.1]{Strat}) we deduce the existence of a normal conditional
expectation and hence it suffices to assume $A_1=A_2$,
$\phi_1=\phi_2$.

Now, we assume $A=M_n$ and that $\phi_n(x)=\sum_{k=1}^n \la_k
x_{kk}$ is the given state. We recall the notation $l(h)$ for the
creation operator on the full Fock space $\F(H)$. On
$\B(\F(\ell_2^{n^2}))\ten M_n$ we consider the state $\Phi\lel
\phi_\Om \ten \phi_n$, where $\phi_\Om$ is given by the vacuum
state. According to \cite[Theorem 5.2]{S1}, we know that the
$C^*$-algebra $C^*(L)$ generated by the  operator
 \[ L\lel \summ_{k,l=1}^n \ell(h_{kl}) \ten \sqrt{\la_k}e_{kl} \]
is free from $M_n$. We consider the semicircular operator
 \[ s\lel L+L^* \lel \summ_{k,l=1}^n
 [\sqrt{\la_k} \ell(h_{kl})+\sqrt{\la_l}\p \ell(h_{lk})^*] \ten e_{kl} \pl .\]
For $k\le l$ we obtain the generalized semicircular elements
 \[ y_{kl} \lel \sqrt{\la_k}\ell(h_{kl})+\sqrt{\la_l}\ell(h_{kl})^* \lel
  \sqrt{\la_k} (\ell(h_{kl})+\sqrt{\frac{\la_l}{\la_k}}\ell(h_{lk})^*)
 \pl
 .\]
By orthogonality we deduce  that the family $(y_{kl})_{k\le l}$ is
$^*$-free. For fixed $l\le k$ the von Neumann algebra $D_{kl}$
generated by $y_{kl}$ is isomorphic to the von Neumann algebra
$T_{\la_l/\la_k}$ introduced in \cite[section4]{S1}. Moreover, the
restriction $\Phi_{\Om}|_{D_{kl}}$ corresponds to the vacuum state
$\phi_{\la_k/\la_l}$ on $T_{\la_l/\la_k}$. Hence the von Neumann
algebra $M$ generated by all the $D_{kl}$'s is isomorphic to
$\bar{\ast}_{k\le l} (T_{\la_k/\la_l},\phi_{\la_k/\la_l})$. In
particular, $\Phi$ is faithful on the von Neumann algebra $N=M\ten
M_n$ and the isomorphism between $N$ and $\bar{\ast}_{k\le l}
(T_{\la_k/\la_l},\phi_{\la_k/\la_l}) \ten M_n$ sends $\ast_{k\le
l} \phi_{\la_k/\la_l}\ten \phi_n$ to  $\Phi$. In \cite[Theorem
2.9]{S1} Shlyakhtenko investigated an automorphism group $\al_t$
of $T_{\la_k/\la_l}$ and showed that it satisfies the KMS
condition at inverse temperature $1$ with respect to
$\phi_{\la_k/\la_l}$. In our normalization the modular group
$\si_t^{\phi}$ of an arbitrary normal faithful state $\phi$
satisfies $\phi(xy)\lel \phi(y\si_{-i}^{\phi}(x))$, i.e. the KMS
condition at inverse temperature $-1$. This implies that
$\si_t^{\phi_{\la_k/\la_l}}=\al_{-t}$. We refer to
\cite[section4]{S1} for the equation $\al_t(y_{kl})= \kla
\la_l/\la_k \mer^{-it} y_{kl}$. Thus we have
 \[ \si_t^{\phi_{\la_k/\la_l}}(y_{kl})\lel \kla \frac{\la_l}{\la_k}\mer^{it}
 y_{kl} \pl. \]
(See \cite{Ps6} for a direct argument.) We deduce from
$\si_t^{\phi_n}(e_{kl})\lel \la_k^{it}e_{kl}\la_l^{-it}$ that
 \[ \si_t^{\Phi}(y_{kl}\ten e_{kl})
 \lel \kla \frac{\la_l}{\la_k}\mer^{it}   \kla
 \frac{\la_k}{\la_l}\mer^{it}  (y_{kl}\ten e_{kl}) \lel y_{kl}\ten e_{kl}   \pl .\]
This implies that $s$ belongs to the centralizer $N^{\Phi}$ of the
von Neumann  $N$. The von Neumann algebra $W^*(s)$ generated by
$s$ is isomorphic to $L_{\infty}[0,1]$. Let $u\in W^*(s)$ be a
Haar unitary. Clearly, the algebras $uM_nu^*$ and $M_n$ are free
with respect to $\Phi$. Moreover, we deduce from
$\si_t^{\Phi}(u)=u$ that
 \[ \phi(uxu^*)\lel \phi(xu^*u)\lel \phi(x) \pl .\]
Therefore the von Neumann  algebra $M$ generated by $M_n$ and
$uM_nu^*$ is isomorphic to the free product
$(M_n,\phi_n)\bar{\ast}(M_n,\phi_n)$. Since $M_n$ and $uM_nu^*$
are invariant under $\si_t^{\phi}$, we find a normal conditional
expectation from $N=\bar{\ast}_{k\le l}
(T_{\la_k/\la_l},\phi_{\la_k/\la_l})\ten M_n$ onto $M$. According
to \cite[Lemma 2.5]{PS}, the von Neumann algebra $\bar{\ast}_{k\le
l} (T_{\la_k/\la_l},\phi_{\al_k/\la_l})$ is QWEP and the assertion
follows.\qd

\begin{rem} {\rm The preceding argument implies in particular that
 \[ M_n\big ( \bar{\ast}_{k\le l} (T_{\la_k/\la_l},\phi_{\la_k/\la_l})\big)\lel (M_n,\phi_n)\bar{\ast}L_{\infty}[0,1]\]
In the tracial situation we need $n+2\frac{(n-1)n}{2}=n^2$
semicircular random variables. This leads to the well-known
isomorphism $M_n\bar{\ast}L(\zz)=M_n(L(\ff_{n^2}))$ (see
\cite[Theorem 5.4.1]{VDN}).}
\end{rem}

Our proof of Theorem \ref{freeQWEP} requires us to show  that free
products and  ultraproducts are compatible. We will need some
notation. Let  $(A_j,\phi_j)_{j\in J}$ be a family of von Neumann
algebras $A_j$ with normal faithful state  $\phi_j$. For a free
ultrafilter $\U$ on $J$, we consider the von Neumann algebra
$B=(\prod_{j,\U} (A_j)_*)^*$ and the ultraproduct state
$\Phi(x_j)\lel \lim_{j,\U} \phi_j(x_j)$. Let us denote by $e$ the
support of $\Phi$. Then $\Phi$ is a normal faithful state on $e
Be$. We use the notation $\prod_{j,\U} [A_j,\phi_j]=e(\prod_{j,\U}
(A_j)_*)^*e$ for this von Neumann algebra. (This is the
non-tracial version of the usual von Neumann algebra ultraproduct
of $II_1$-von Neumann algebras.)\lz

\begin{lemma}\label{ultraa} Let $A_1$, $A_2$ be von Neumann algebras with normal faithful
states  $\varphi_1$ and $\varphi_2$. Let $\pi_k:A_k\to \prod_j
[A_{k,j},\phi_{k,j}]$ be a faithful state preserving
homomorphisms. Assume that the ultraproduct states
$\Phi_k((x_j))=\lim_{j,\U}\phi_{k,j}(x_j)$ satisfy
 $\si_t^{\Phi_{k}}\circ \pi_k \lel \pi_k
\circ \si_t^{\varphi_k}$ for $k=1,2$. Then there is an injective
$^*$-homomorphism
 \[ \al:(A_1,\varphi_1)\bar{\ast}(A_2,\varphi_2)\to \prodd_{\U} [(A_{1,j},\phi_{1,j})\bar{\ast}(A_{2,j},\phi_{2,j}),
 \phi_{1,j}\ast \phi_{2,j}]\]
together with a normal conditional expectation onto
$\al((A_1,\varphi_1)\bar{\ast}(A_2,\varphi_2))$.
\end{lemma}

\begin{proof}[\bf Proof:]
We need a preliminary observation. Let $N\subset M$, $E:M\to N$ be
a normal conditional expectation, and let $\phi:N\to \cz$ be a
normal state. We denote by $e$ the support of $\phi$ and by $f$
the support of $\psi=\phi\circ E$. We want to show that every
element in  $x\in eMe$ commutes with $f$. Let $E_*:N_*\to M_*$ be
the predual map such that $(E_*)^*=E$. Note that $E_*$ is a
$N$-bimodule map and $E_*(\phi)= \phi\circ E=\psi$. Given an
analytic element $x\in eNe$, we see that
 \begin{align*}
  (1-f)x\psi &= (1-f)xE_*(\phi) \lel   (1-f)E_*(x\phi)\lel (1-f)E_*(\phi\si^{\phi}_{-i}(x))\\
  &= (1-f)E_*(\phi) \si^{\phi}_{-i}(x)\lel (1-f)\psi
 \si^{\phi}_{-i}(x)\lel 0 \pl .
 \end{align*}
Thus, for analytic elements, $(1-f)xf=0$ and $(1-f)x^*f=0$. By the
density of the analytic elements, this implies that $xf=fx$ for
all $x\in eNe$.

Let $(A_1,\phi_1)$, $(A_2,\phi_2)$ and $(A_{k,j},\phi_{k,j})$
satisfy the assumptions. Let us denote by $\psi_j=\phi_{1,j} \ast
\phi_{2,j}$ the free product state and
$B_j=(A_{1,j},\phi_{1,j})\bar{\ast}(A_{2,j},\phi_{2,j})$. We
denote by $\pi_{k,j}$ the natural inclusion map of $A_{k,j}$ in
$B_j$ using the first (second) component for $k=1$ ($k=2$,
respectively). Then we find a mapping $\hat{\rho}_k: \prod_{j,\U}
A_{k,j}\to \prod_{j,\U} B_j$ given by $\hat{\rho}_k(a_j)\lel
(\pi_{k,j}(a_j))$. By density with respect to the strong operator
topology, we may extend $\hat{\rho}_k$ to a $^*$-homomorphism from
$(\prod_{j,\U} (A_{k,j})_*)^*$ to $(\prod_{j,\U} (B_j)_*)^*$,
still denoted by $\hat{\rho}$. For every $j$ and $k=1,2$ we have a
normal conditional expectation $E_{k,j}:B_j\to A_{k,j}$ such that
$\psi_j\lel \phi_{k,j}\circ E_{k,j}$. The predual maps provide
contractions $(E_{k,j})_*:(A_{k,j})_*\to (B_{j})_*$ such that
$((E_{k,j})_*)^*=E_{k,j}$. This induces (complete) contractions
$(E_k)_*:\prod_{j,\U} (A_{k,j})_*\to \prod_{j,\U} (B_{j})_*$ such
that $((E_k)_*)^*\hat{\rho}_k=id$. We follow Raynaud's notation
$(\pl )^{\bullet}$ (see \cite{Ray})   for equivalence classes in
ultraproducts and in $(\prod_{j,\U} A_*)^*$. In this notation we
have
 \[ (E_k)_*(\phi_{k,j})^{\bullet}\lel (\psi_j)^{\bullet} \pl .\]
Let us denote by $\Psi=(\psi_j)^{\bullet}$ the ultraproduct state.
The support of $\Psi$, $\Phi_k$  is denoted by $f$, $e_{k}$,
respectively. By our preliminary observation we see that elements
in $\hat{\rho}_k(e_{k}(\prod_{j,\U} (A_{k,j})_*)^*e_{k})$ commute
with $f$. Therefore, $\hat{\rho}_k$ induces a $^*$-homomorphism
$\rho_k(e_kxe_k)=f\hat{\rho}_k(e_kxe_k)f$ between the von Neumann
algebras $\prod_{j,\U}[A_{k,j},\phi_{k,j}]$ and
$\prod_{j,\U}[B_{j},\psi_{j}]$ for $k=1,2$.  For the modular group
we use Raynaud's \cite{Ray} characterization
 \begin{equation} \label{ra} \si_t^{\Psi}(f (y_j)^{\bullet}f)
 \lel f\big(\si_t^{\psi_j}(y_j)\big)^{\bullet} f \pl .
 \end{equation}
A  similar formula holds for the  $\Phi_{k}$'s.  We certainly have
$f\le \rho_k(e_{k})$. Thus we get
\begin{align*}
  &\si_t^{\Psi}\big(f \hat{\rho}_k\big(e_{k}(x_j)^{\bullet}e_{k}\big)f\big)
  \lel \si_t^{\Psi}\big(f\hat{\rho}_k\big((x_j)^{\bullet}\big)f\big)
  \lel \si_t^{\Psi}\big(f\big(\pi_{k,j}(x_j)\big)^{\bullet} f\big)
  \lel f\big (\si_t^{\psi_j}(\pi_{k,j}(x_j))\big)^{\bullet}f\\
  &= f\big(\pi_{k,j}(\si_t^{\phi_{k,j}}(x_j)) \big)^{\bullet} f
  \lel f\hat{\rho}_k\big(e_{k}(\si_t^{\phi_{k,j}}(x_j))^{\bullet}e_k\big)f \lel
  \rho_k\big( \si_t^{\Phi_{k}}(e_k(x_j)^{\bullet}e_k)\big) \pl .
 \end{align*}
By density we find that $\si_t^{\Psi}\circ \rho_k \lel \rho_k
\circ \si_t^{\Phi_k}$ for $k=1,2$. Thus we have found embeddings
$\al_k\lel \rho_k\pi_k$ of $A_k$ in $\prod_{j,\U}[B_j,\psi_j]$
satisfying $\Psi\circ \al_k\lel \varphi_k$.  Using the assumption
$\si_t^{\Phi_{k}}\circ \pi_k \lel \pi_k \circ \si_t^{\varphi_k}$,
we deduce that
 \begin{equation} \label{com}
 \si_t^{\Psi}\circ \al_k\lel \al_k\circ \si_t^{\varphi_k} \pl.
 \end{equation}
If we can show that $\al_1(A_1)$ and $\al_2(A_2)$ are free with
respect to $\Psi$, then the von Neumann algebra $C$ generated by
$\al_1(A_1)$ and $\al_2(A_2)$ is isomorphic to the free product
$(A_1,\phi_2)\bar{\ast}(A_2,\phi_2)$ and admits a normal
conditional expectation as guaranteed by   Takesaki's theorem (see
e.g. \cite[Theorem 10.1]{Strat}). Let us show freeness. Consider
 $a_r\in {\rm \AA}_{i_r}$ and assume  $i_1\neq i_2 \cdots
\neq i_m$. We may replace the $a_r$'s by their analytic
approximations
 \[ T_l^{\varphi_{i_r}}(a_r)\lel \int_{\rz} \si_t^{\varphi_{i_r}} f_l(t) dt \pl,\]
where $f_l(t)\lel \sqrt{\frac{l}{\pi}}e^{-lt^2}$. Note that
$T_l^{\varphi_{i_r}}(a_r)$ still has mean $0$. Since
$e_k(\prod_{j} A_{j,k})e_k$ is strongly dense in
$\prod_{j,\U}[A_{j,k},\phi_{k,j}]$ we may apply Kaplansky's
density theorem to approximate $\al_{i_r}(a_r)$ by a bounded net
of elements $(a_{i_r,j,s})_{s\in S}$ in the strong operator
topology. We observe that $\lim_{s}
\phi_{i_r}((a_{i_r,j,s})^{\bullet})\lel \varphi_{i_r}(a_r)=0$.
Therefore
$((a_{i_r,j,s}-\phi_{i_r,j}(a_{r,j,s})1)^{\bullet})_{s\in S}$ also
provides a bounded net converging in the strong operator topology
to  $\al_{i_r}(a_r)$. We use the notation
$\circled{a}_{i_r,j,s}=a_{i_r,j,s}-\phi_{i_r,j}(a)1$. Using the
strong continuity of the modular group it is easy to show that
$T_l^{\Psi}\big(f\big(\pi_{i_r,j}(\circled{a}_{i_r,j,s})\big)^{\bullet}f\big)$
converges to $T_l^{\Psi}(\al_{i_r}(a_r))$ in the strong operator
topology. We deduce from \eqref{ra} that
 \[ T_l^{\Psi}\big(f\big(\pi_{i_r,j}(\circled{a}_{i_r,j,s})\big)^{\bullet}f\big )
 \lel
 f\big(T_l^{\psi_j}(\pi_{i_r,j}(\circled{a}_{i_r,j,s})) \big)^{\bullet} f
 \pl .\]
For an arbitrary family $(y_j)$ we observe that
 \begin{align*}
 & (1-f) \big(T_l^{\psi_j}(y_j)\big)^{\bullet} \Psi \lel
 (1-f)  \big(T_l^{\psi_j}(y_j)\psi_j\big)^{\bullet}
 \lel (1-f)\big(\psi_j \si_{-i}^{\psi_j}(T_l(y_j)\big)^{\bullet} \\
 &= (1-f)\Psi (\si_{-i}^{\psi_j}(T_l(y_j))^{\bullet}
 \lel 0
 \pl .\end{align*}
In the last line we use the fact that the norm of
$\si_{-i}^{\psi_j}(T_l(y_j)$ can be estimated as a function of
$\noo y_j\rrm$ and $l$. Applying this argument  also for adjoints,
we deduce that
 \begin{equation}\label{comm}  f\big(T_l^{\psi_j}(y_j)\big)^{\bullet} \lel
 \big(T_l^{\psi_j}(y_j)\big)^{\bullet}f \pl .
 \end{equation}
Using \eqref{com}, \eqref{comm} and the component-wise freeness we
get
 \begin{align*}
 & \Psi\bigg(\al_{i_1}(T_{l_1}^{\varphi_{i_1}}(a_1))\al_{i_2}(T_{l_2}^{\varphi_{i_2}}(a_2))
 \cdots \al_{i_m}(T_{l_m}^{\varphi_{i_m}}(a_m))\bigg)\\
 &\pll =  \Psi\bigg(T_{l_1}^{\Psi}(\al_{i_1}(a_1))\al_{i_2}(T_{l_2}^{\varphi_{i_2}}(a_2))
 \cdots \al_{i_m}(T_{l_m}^{\varphi_{i_m}}(a_m))\bigg)\\
 &\pll = \lim_{s_1} \Psi\bigg(T_{l_1}^{\Psi}\big(f\big(
 (\pi_{i_1,j}(\circled{a}_{i_1,j,s_1})\big)^{\bullet}f\big)
 \al_{i_2}(T_{l_2}^{\varphi_{i_2}}(a_2))
 \cdots \al_{i_m}(T_{l_m}^{\varphi_{i_m}}(a_m))\bigg)\\
 &\pll = \lim_{s_1}
\Psi\bigg(f\big(T_{l_1}^{\psi_j}(\pi_{i_1,j}(\circled{a}_{i_1,j,s_1})\big)^{\bullet}f
 \al_{i_2}(T_{l_2}^{\varphi_{i_2}}(a_2))
 \cdots \al_{i_m}(T_{l_m}^{\varphi_{i_m}}(a_m))\bigg)\\
 &\pll = \lim_{s_1}\lim_{s_2}\cdots \lim_{s_m} \Psi\bigg(
 \big(T_{l_1}^{\psi_j}\pi_{i_1,j}(\circled{a}_{i_1,j,s_1})\big)^{\bullet}f
 \big(T_{l_2}^{\psi_j}\pi_{i_2,j}(\circled{a}_{i_2,j,s_2})\big)^{\bullet}f\cdots f
 \big(T_{l_m}^{\psi_j}\pi_{i_m,j}(\circled{a}_{i_m,j,s_m})\big)^{\bullet}\bigg) \\
 &\pll = \lim_{s_1}\lim_{s_2}\cdots \lim_{s_m} \Psi\bigg(
 \big(T_{l_1}^{\psi_j}\pi_{i_1,j}(\circled{a}_{i_1,j,s_1})\big)^{\bullet}
 \big(T_{l_2}^{\psi_j}\pi_{i_2,j}(\circled{a}_{i_2,j,s_2})\big)^{\bullet} \cdots
 \big(T_{l_m}^{\psi_j}\pi_{i_m,j}(\circled{a}_{i_m,j,s_m})\big)^{\bullet}\bigg) \\
 &\lel \lim_{s_1}\lim_{s_2}\cdots \lim_{s_m}  \lim_{j,\U}
 \psi_j(T_{l_1}^{\psi_j}\pi_{i_1,j}(\circled{a}_{i_1,j,s_1})T_{l_2}^{\psi_j}\pi_{i_2,j}(\circled{a}_{i_2,j,s_2}))\cdots
 T_{l_m}^{\psi_j}\pi_{i_m,j}(\circled{a}_{i_m,j,s_m})) \lel 0 \pl .
 \end{align*}
Taking the limit for $l_1,...,l_m\to \infty$ yields the
assertion.\qd

\begin{proof}[\bf Proof of Theorem \ref{freeQWEP}]
Let $A_1$ and $A_2$ be von Neumann algebras with QWEP. According
to \cite{Fub} (and a slight perturbation argument), we find state
preserving  embeddings $\pi_1:A_1\to \prod_{i,\U}
[M_{m(i)},\phi_i]$ and $\pi_2:A_2\to \prod_{j,\U'}
[M_{m(j)},\phi_j]$ together with state preserving conditional
expectations $E_1: \prod_{i,\U} [M_{m(i)},\phi_i]\to A_1$, $E_2:
\prod_{j,\U'} [M_{m(j)},\phi_j]\to A_1$. This implies that
$\pi_1(A_1)$, $\pi_2(A_2)$ are  invariant under the modular group
of the ultraproduct states  $\Phi_1=(\phi_i)^{\bullet}$,
$\Phi_2=(\phi_j)^{\bullet}$, respectively. We consider the index
set $I\times J$ with the ultrafilter $\U''$ defined as follows:
$B\in \U''$ if and only $\{i\p|\p \{j\p|\p (i,j)\in B\}\in
\U'\}\in\U'$. We define $(A_{1,ij},\phi_{1,ij})\lel
(M_{m(i)},\phi_i)$ and $(A_{2,ij},\phi_{1,ij})\lel
(M_{m(j)},\phi_j)$. According to Lemma \ref{ultraa}, we deduce
that $(A_1,\varphi_1)\bar{\ast}(A_2,\varphi_2)$ embeds in
 \[
\prodd_{(i,j),\U''} [(M_m(i),\phi_i)\bar{\ast} (M_m(j),\phi_j),
 \phi_i\ast \phi_j] \pl \]
and admits a normal conditional expectation. Therefore Lemma
\ref{fidim} and Lemma \ref{ultraqwep} imply the assertion.\qd

\lz { AMS Subject Classification 2000:} 47L25 \llz

\end{document}